\documentclass{article}
\pdfoutput=1
\usepackage[fleqn]{amsmath}
\usepackage{amssymb}
\usepackage{bbm}
\usepackage[letterpaper, portrait, margin=1in]{geometry}
\usepackage{graphicx}
\graphicspath{ {img/}}
\usepackage{algorithm}
\usepackage[noend]{algpseudocode}
\usepackage{subcaption}
\usepackage{array}
\usepackage[dvipsnames]{xcolor}

\newtheorem{definition}{Definition}

\newcolumntype{L}[1]{>{\raggedright\let\newline\\\arraybackslash\hspace{0pt}}m{#1}}
\newcolumntype{C}[1]{>{\centering\let\newline\\\arraybackslash\hspace{0pt}}m{#1}}
\newcolumntype{R}[1]{>{\raggedleft\let\newline\\\arraybackslash\hspace{0pt}}m{#1}}
\usepackage{multirow}
\usepackage{float}
\floatstyle{plaintop}
\restylefloat{table}
\usepackage{enumitem}
\newcommand{\etal}{\textit{et al}. }
\newcommand{\ie}{i.e., }
\newcommand{\eg}{e.g., }

\newcommand{\sNS}{\mathcal{N}}
\newcommand{\sES}{\mathcal{E}}
\newcommand{\sTS}{\mathcal{T}}

\newcommand{\sAS}{\mathcal{A}}
\newcommand{\sAC}{\mathcal{A}_{c}}
\newcommand{\sACd}{\hat{\mathcal{A}}_{c}}
\newcommand{\sN}{\mathcal{N}^\text{x}}
\newcommand{\sE}{\mathcal{E}^\text{x}}

\newcommand{\sA}{\mathcal{A}^\text{x}}

\newcommand{\flow}[2][]{\varphi_{{#2},{#1}}}

\newcommand{\bincf}[1]{y_{{#1}}}
\newcommand{\bmcf}[2]{x_{{#1},{#2}}}

\newcommand{\nin}[1]{\delta^-(#1)}
\newcommand{\nout}[1]{\delta^+(#1)}

\newcommand{\st}{\mbox{s.t.}\quad }

\newcommand{\eu}[2][]{\ensuremath{u_{{#2}}^{{#1}}}}

\newcommand{\GS}{\ensuremath{\mathcal{G}}}

\newcommand{\footremember}[2]{%
    \footnote{#2}
    \newcounter{#1}
    \setcounter{#1}{\value{footnote}}%
}
 
\title{Large-Scale Zone-Based Evacuation Planning: \\
Models, Algorithms, and Evaluation}
\author{Mohd. Hafiz Hasan\footremember{UM}{University of Michigan, Ann Arbor. University of Michigan, Ann Arbor, Michigan 48109, Email: hasanm@umich.edu.} \and Pascal Van Hentenryck\footremember{GT}{Georgia Institute of Technology, Atlanta, Georgia 30332, Email: pvh@isye.gatech.edu.}}
\date{}

\begin{document}

\maketitle

\begin{abstract}
  In zone-based evacuation planning, the region to evacuate is divided
  into zones and each zone must be assigned a path to safety and
  departure times along the path. Zone-based evacuations are highly
  desirable in practice because they allow emergency services to
  communicate evacuation orders and to control the evacuation more
  accurately.  Zone-based evacuations may also be combined with
  contraflows (to maximize the network capacities) and may impose
  additional constraints on the evacuation path (e,g,. path
  convergence) and the departure times (e.g., non-preemption).

  This paper presents a systematic study large-scale zone-based
  evacuation planning, both from an effectiveness and a computational
  standpoint.  It reviews existing optimization algorithms, and
  presents new ones, and evaluates them, on a real, large-scale case
  study, both from a macroscopic standpoint and through microscopic
  simulations under a variety of assumptions. The results provide some
  unique perspectives on the strengths and weaknesses of each approach
  and the implications of evacuation functionalities. The paper also
  suggests new directions for future research in zone-based evacuation
  and beyond in order to address the fundamental challenges by
  emergency services around the world.
\end{abstract}

\section{Introduction}

Large-scale evacuations are often necessary and critical to the
preservation of safety and lives of residents in regions threatened by
man-made or natural disasters like floods, hurricanes, and
wildfires. According to a report by the International Federation of
Red Cross and Red Crescent Societies \cite{cross2010}, the first
decade of the 21\textsuperscript{st} century witnessed 7184 disasters
around the world which affected a total of 2.55 billion people,
accounted for the deaths of more than 1 million people, and incurred
\$986 billion in economic losses. There is now an evacuation of 1,000
or more people every two or three weeks in the United States alone.

Effective disaster management requires, among others, evacuation plans
that ensure resources like transportation network capacity and time
are not completely overwhelmed by evacuee demand. Traffic congestion
and associated delays which usually result from self-evacuations, in
which individuals are given the freedom to choose their own evacuation
routes, destinations, and times, significantly increase the risk of
casualties being stranded in disaster affected areas. Therefore, it is
crucial for emergency authorities to be equipped with centralized
disaster management tools capable of generating and prescribing plans
that guarantee optimal utilization of evacuation resources to attain
specific goals like maximizing the number of evacuees reaching safety
or minimizing overall evacuation time. Evacuation planning algorithms
fulfill this need by producing prescriptive evacuation plans, i.e., a
set of operational instructions for authorities to manage and
orchestrate large-scale evacuations through the specification of
directions for evacuees, including routes from their homes to
designated safe destinations and departure times to be followed, as
well as identification of roads that need to be closed to facilitate
traffic flow. This contrasts to self-evacuations that are more
difficult to control and may produce significant congestion.

Hamacher and Tjandra \cite{hamacher2002} distinguish between
microscopic and macroscopic approaches to evacuation modeling.
Microscopic approaches model individual characteristics of evacuees,
their interactions with each other, and how these factors influence
their movement. In contrast, macroscopic approaches aggregate evacuees
and model their movements as flow in a network, making them much more
amenable to optimization. Macroscopic models are often defined in
terms of {\em flows over time} in order to capture capacity
constraints more accurately.  In particular, they typically use the
concept of {\em time-expanded graphs} pioneered by Ford and Fulkerson
\cite{ford1962}.

This study is concerned with macroscopic approaches to prescriptive
evacuation planning, although all results are validated using
microscopic traffic simulations. Moreover, it focuses on {\em
  zone-based evacuations} in which all evacuees from the same
residential zone are assigned a single evacuation route. Most
emergency services rely on zone-based evacuations which facilitate the
communication of evacuation plans, reduce confusion, increase
compliance, and allow for a more reliable control of the evacuation
Indeed, zone-based evacuations are probably the only practical method
for communicating instructions precisely to the population at
risk. They are however much more computationally challenging to plan
and finding scalable algorithms has been one of the foci of recent
research.

The core Zone-based Evacuation Planning Problem (ZEPP) considered in
this study consists in assigning an evacuation path, as well as
departure times, to each zone in the region. However, even when
restricting attention to zone-based evacuations, several important
design decisions remain to be taken: They include, but are not limited
to, contraflows, convergent plans, and preemption.

\begin{itemize}
\item {\em Contraflows}, also known as lane reversals, are the idea of using
inbound lanes for outbound traffic in evacuations. Several studies have
suggested contraflow procedures as a viable method for increasing
network capacity (e.g., \cite{wolshon2001,theodoulou2004}).

\item {\em Convergent plans} ensure that each evacuee coming to an
  intersection follows the same path subsequently. The rationale for
  convergent plans is that they eliminate forks from all evacuation
  paths and hence reduce driver hesitation at road intersections,
  which has been shown to be a significant source of delays
  \cite{townsend2006}. Convergent paths also allow roads which are not
  part of the evacuation paths to be blocked, facilitating vehicular
  guidance and enforcement of the evacuation plans.

\item {\em Non-preemptive} evacuations ensure that the evacuation of a
  zone, once it starts, proceeds without interruptions. Non-preemptive
  evacuations are also preferred by emergency services since they are
  easier to enforce. 
\end{itemize}

\noindent
Each of these decisions has a significant impact, not only on the
effectiveness of an evacuation (e.g., the number of evacuees reaching
safety), but also on the computational properties of the optimization
model and its potential solution techniques.

The purpose of this paper is to provide a systematic comparative study
of zone-based evacuation planning, their design choices, their
fidelity, and their computational performance. It synthesizes existing
algorithms, complements them with some new ones to fill some of the
gaps in the design space, and compares them on a real-life case study
both at macroscopic and microscopic scales. The case study concerns
the the Hawkesbury-Nepean (HN) region located north-west of Sydney,
Australia; The associated time-expanded graph has 30,000 nodes and
75,000 arcs.

The benefits of this study are threefold:
\begin{enumerate}
\item It systematically evaluates, on a large case study and for the
  first time, a variety of zone-based evacuation planning algorithms
  both from macroscopic and microscopic viewpoints;
  
\item It quantifies, for the first time, the benefits and limitations
  of contraflows, convergent plans, and preemption, providing unique
  perspectives on how to deploy these algorithms in practice;

\item It highlights the approaches that are best suited to capture
  each of these design features and the computational burden they
  impose.

\item It suggests potential directions for future research.
\end{enumerate} 
Finally, the paper addresses some other avenues for future research
that are not satisfactorily addressed in the literature yet but are
critical in practice. 

The rest of this paper is organized as follows. Section
\ref{sec:problem_formulation} outlines the key concepts in zone-based
evacuation planning, including time-expanded graphs, contraflows, and
convergent evacuations. Section \ref{sec:MIP} provides a Mixed Integer
Program (MIP) for non-convergent, preemptive evacuation planning with
contraflows. This MIP introduces the main decision variables that
appear in most approaches discussed in the paper. Section
\ref{sec:BN} presents a Benders Decomposition approach
for non-convergent, preemptive ZEPP and Section \ref{sec:BC} reviews
the Benders decomposition approach for the convergent preemptive ZEEP
originally proposed in \cite{romanski2016}. Section \ref{sec:CPG}
presents the conflict-based path-generation heuristic algorithm
proposed in \cite{pillac2016} for the non-convergent, preemptive ZEEP,
while Section \ref{sec:CG} reviews the column-generation approach to
the non-convergent, non-preemptive ZEEP which was originally proposed
in \cite{pillac2015} and improved upon in \cite{Hafiz2017}. Section
\ref{sec:casestudy} presents the case study, the experimental setting,
and some initial observations. Section \ref{sec:macroscopic} presents
the experimental results from a macroscopic standpoint and Section
\ref{sec:microscopic} provides the details of the microscopic
simulations.  Section \ref{sec:perspective} gives some perspectives on
directions for future research and knowledge gaps that needs to be
filled. Section \ref{sec:related} discusses some related work to give
more context to all the results in the paper. Finally, Section
\ref{sec:conclusion} concludes the paper.

\section{Problem Formulation}
\label{sec:problem_formulation}

\begin{figure}[!tbhp]
\centering
\includegraphics[width=0.5\linewidth]{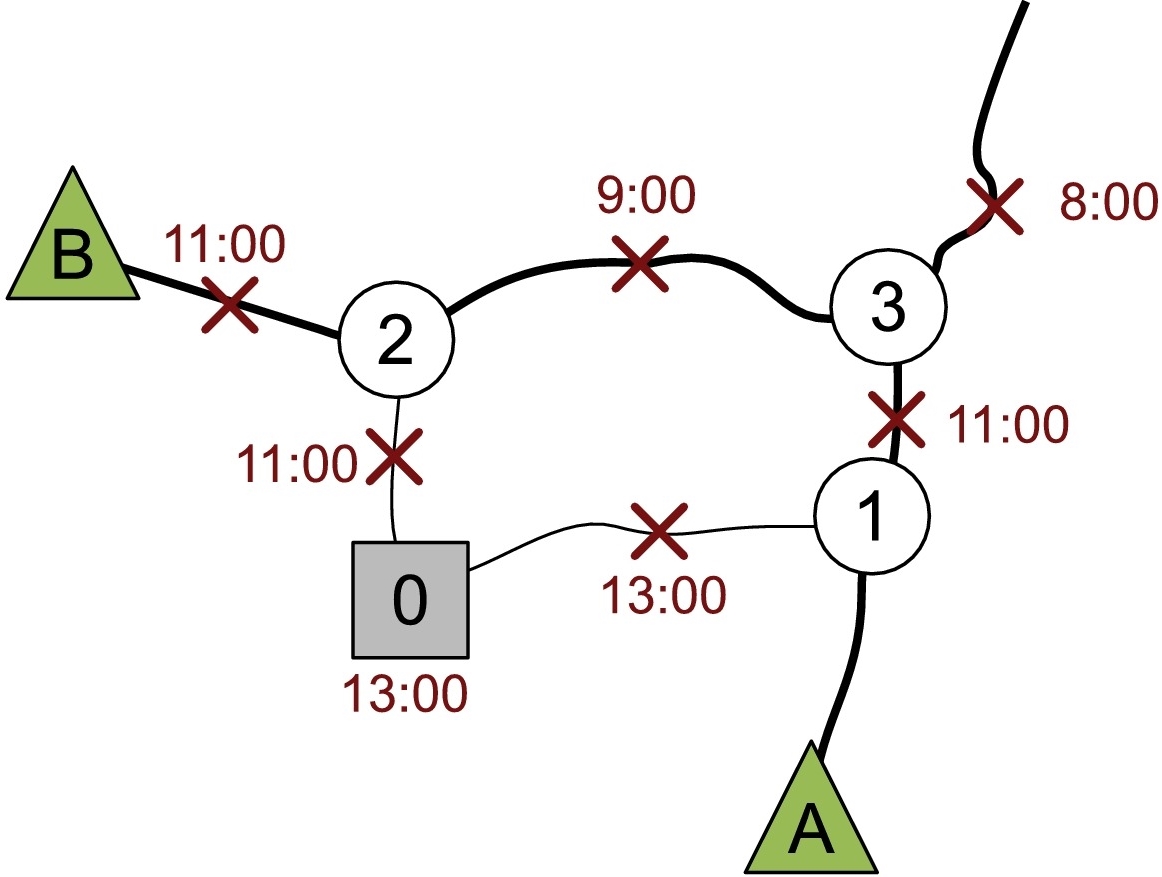}
\caption{A Sample Evacuation Scenario.}
\label{fig:evac_scenario}
\end{figure}

Figure \ref{fig:evac_scenario} shows an example of the evacuation
scenario that is addressed by all methods presented in this
work. Square 0 represents an evacuation node (e.g., a residential zone), triangles A and B represent safe nodes (final evacuation
destinations), circles 1-3 represent transit nodes (road
intersections), and arcs represent roads connecting the nodes. Times
on each arc indicate when each road will become unavailable (\eg due
to being flooded) and the time on the evacuation node indicates the
final deadline by which it must be evacuated. In this example, the
evacuation deadline for node 0 is 13:00 since its last outgoing arc
will be blocked at that time.

\begin{figure}[!tbhp]
\centering
\includegraphics[width=0.5\linewidth]{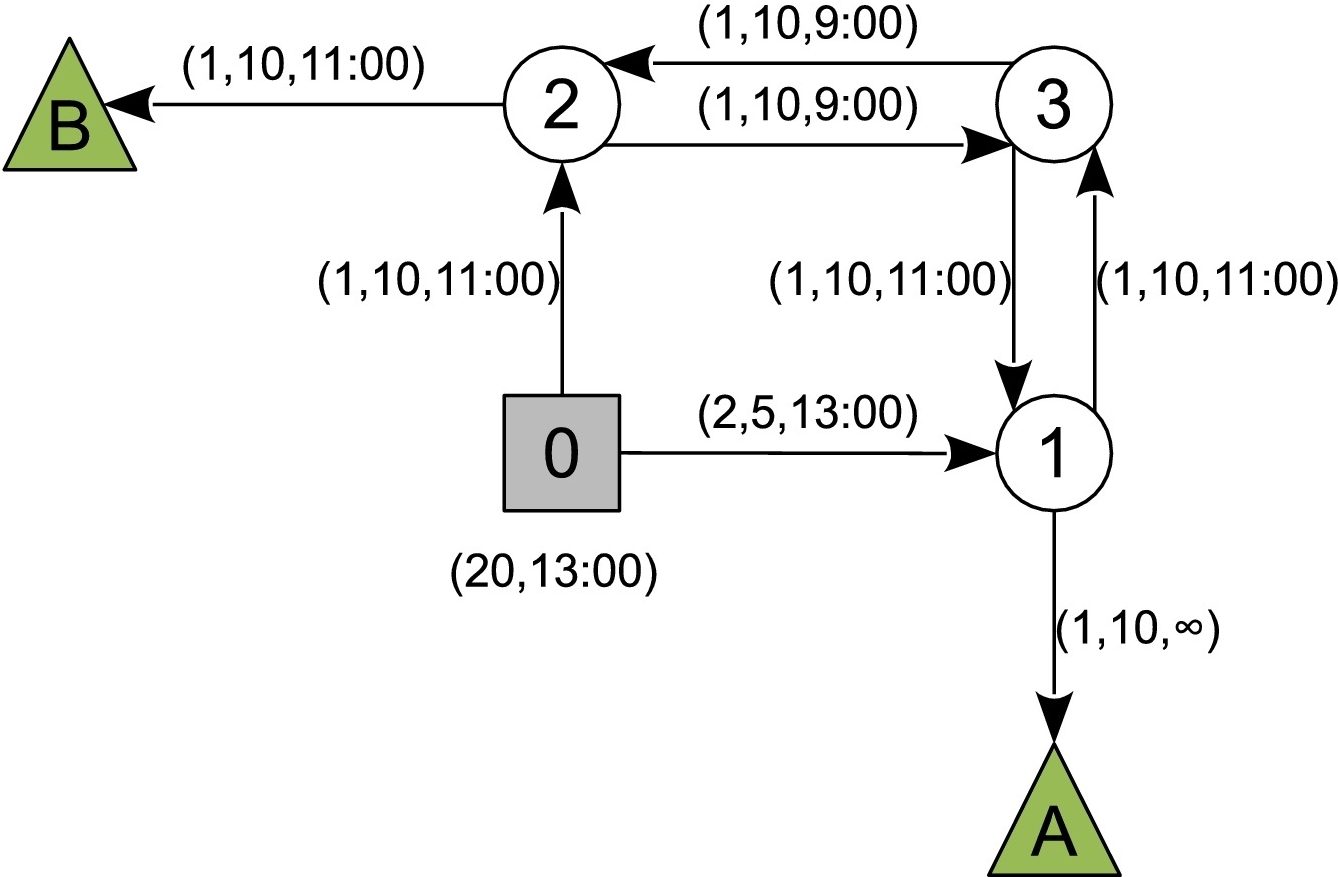}
\caption{The Static Evacuation Graph for the Scenario in Figure \ref{fig:evac_scenario}.}
\label{fig:evac_graph}
\end{figure}

The evacuation scenario can be abstracted by a static evacuation graph
$\mathcal{G}=(\mathcal{N}=\mathcal{E}\cup\mathcal{T}\cup\mathcal{S},\mathcal{A})$
where $\mathcal{E}$, $\mathcal{T}$, and $\mathcal{S}$ are the set of
evacuation, transit, and safe nodes respectively, and $\mathcal{A}$ is
the set of all arcs. Each evacuation node $k\in\mathcal{E}$ has
associated with its demand $d_k$ representing the number of vehicles
to be evacuated and its evacuation deadline $f_k$, while each arc
$e\in\mathcal{A}$ has associated with its travel time $s_e$, its
capacity $u_e$ in terms of vehicles per unit time, and its block time
$f_e$ at which the road becomes unavailable. Figure
\ref{fig:evac_graph} shows the static evacuation graph for the
scenario of Figure \ref{fig:evac_scenario}. The evacuation node is
labeled with its demand and evacuation deadline while the arcs are
labeled with their travel time, capacity, and block time. Also note
that the evacuation node has no incoming arcs and the safe nodes have
no outgoing arcs.

\begin{figure}[!tbhp]
\centering
\includegraphics[width=0.5\linewidth]{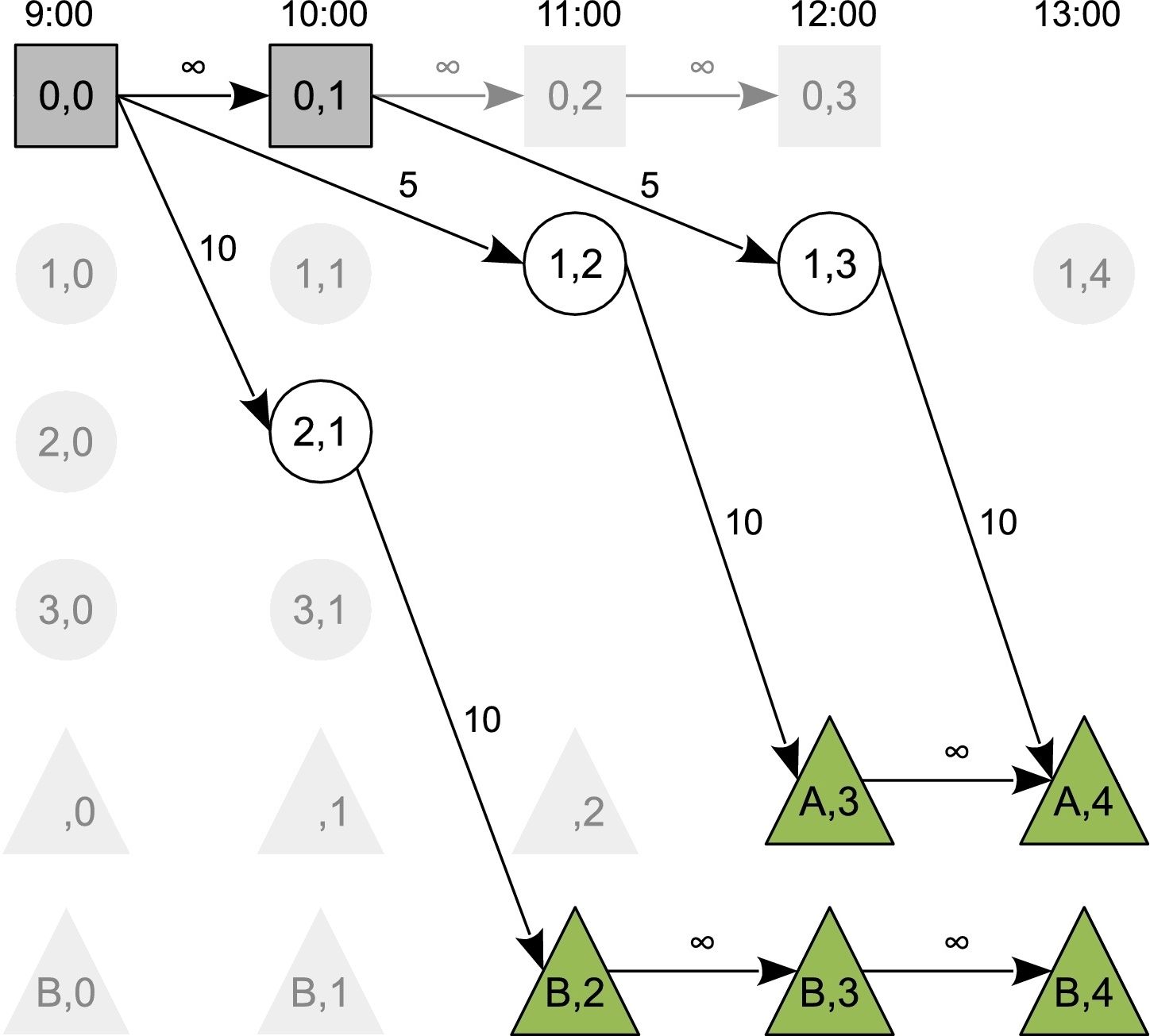}
\caption{The Time-expanded Graph of Static Graph in Figure \ref{fig:evac_graph}.}
\label{fig:time_expanded_graph}
\end{figure}

In order to reason about traffic flows over time, the static graph is
converted into a time-expanded graph
$\mathcal{G}^x=(\mathcal{N}^x=\mathcal{E}^x\cup\mathcal{T}^x\cup\mathcal{S}^x,\mathcal{A}^x)$. The
conversion is performed by first discretizing the time horizon
$\mathcal{H}$ into time steps of identical length $t\in\mathcal{H}$,
creating a copy of all nodes at each time step, and replacing each arc
$e=(i,j)$ with corresponding arcs $e_t=(i_t,j_{t+s_e})$ for each time
step that $e$ is available. Figure \ref{fig:time_expanded_graph} shows
the time-expanded graph constructed from the static graph of Figure
\ref{fig:evac_graph}, where each arc is labeled with its
capacity. Infinite capacity arcs are introduced connecting the
evacuation and safe nodes at each time step to allow evacuees to wait
at those nodes. Nodes that cannot be reached from either the
evacuation or safe nodes within the time horizon are removed from the
graph (they are greyed out in Figure \ref{fig:time_expanded_graph}).

An evacuation plan can then be defined to contain the following two
components: (a) a set of evacuation paths, each represented by a
sequence of connected nodes in the static graph from an evacuation
node to a safe node, specifying the route to be taken by residents of
each evacuation node to reach safety, and (b) an evacuation schedule
indicating the number of vehicles that need to depart from each
evacuation node at each time step $t\in\mathcal{H}$. The Zone-Based
Evacuation Planning Problem (ZEPP) can now be defined as follows.

\begin{definition}
  Given an evacuation graph $\GS$, the Zone-Based Evacuation Planning
  Problem (ZEPP) consists of finding an evacuation path from each
  evacuation zone to a safe node that maximize the flow of evacuees to
  safe nodes, while satisfying the problem constraints.
\end{definition}

\paragraph{Contraflows} Contraflows are an important tool in evacuation planning and
scheduling. To capture their benefits, this study assumes the
existence of a subset $\mathcal{A}_c\subseteq\mathcal{A}$ of arcs in
the static graph that may be used in contraflows during
evacuations. The unique arc that goes in the opposite direction of arc
$e\in\mathcal{A}_c$ is denoted by $\bar{e}$. The set $\mathcal{A}_c$
can then be partitioned into $\hat{\mathcal{A}}_c$ and
$\check{\mathcal{A}}_c$ such that $\forall e\in\hat{\mathcal{A}}_c:
\bar{e}\in\check{\mathcal{A}}_c$. Finally, $e_0\in\sAS$ denotes the
static edge associated with edge $e\in\sA$ and $\delta^{-}(i)$ and
$\delta^{+}(i)$ denote the set of incoming and outgoing arcs of node
$i$ respectively.

\paragraph{Convergent Evacuations} Convergent paths reduce confusion and hesitation
during an evacuation. They can be formally defined by the following definitions. 

\begin{definition}
A graph $\GS = (\sNS, \sAS)$ is \emph{connected} if, for all $k \in \sES$, there exists a path from $k$ to a safe node.
\end{definition}

\begin{definition}
A graph $\GS = (\sNS, \sAS)$ is \emph{convergent} if, for all $i \in \sES \cup \sTS$, the outdegree of $i$ is $1$.
\end{definition}

\noindent
As stated by \cite{even2015}, any connected evacuation graph
$\GS$ contains a connected and convergent subgraph $\mathcal{G}'$. If
an evacuation graph is connected and convergent, each evacuation node
has a unique path to a safe node.  The Convergent Zone-Based
Evacuation Planning Problem (C-ZEPP) is defined as follows:

\begin{definition}
  Given a connected evacuation graph $\GS$, the Convergent Zone-Based
  Evacuation Planning Problem (C-ZEPP) consists of finding a
  convergent subgraph $\mathcal{G}'$ of $\GS$ and a set of evacuee
  departure times that maximize the flow from evacuation nodes to safe
  nodes, while satisfying the problem constraints.
\end{definition}

\begin{figure}[!tbhp]
	\centering
	\includegraphics[width=0.7\linewidth]{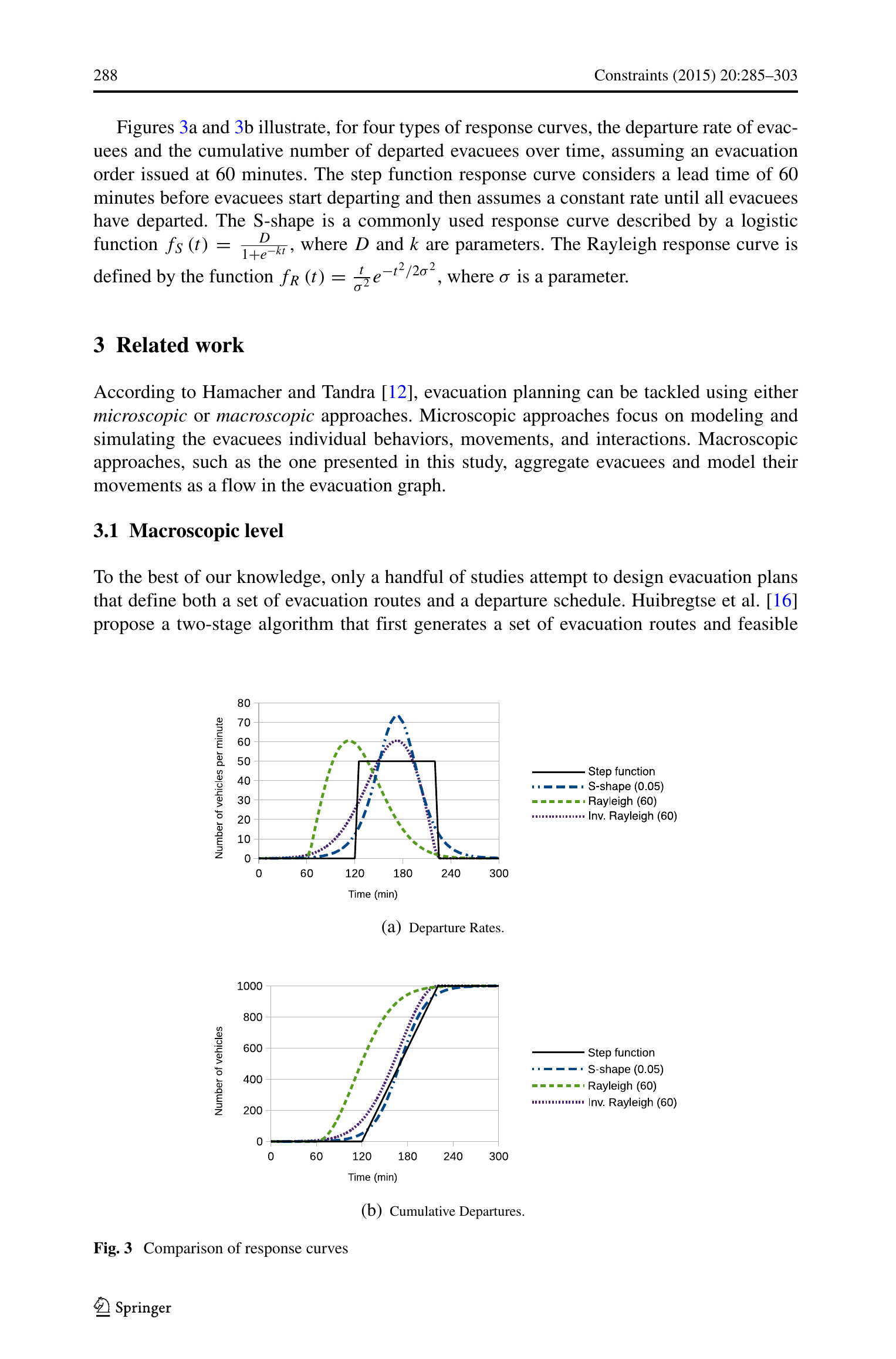}
	\caption{The number of evacuees departing an evacuation node $k$ as a function of time, $D_k(t)$, using four different types of response curves (from \cite{pillac2015}).}
	\label{fig:response_curves}
\end{figure}

\paragraph{Non-Preemption and Response Curves}
Non-preemptive evacuation plans are typically organized around the
concept of response curves \cite{pel2012}. A response curve $f$ is a
function that models the number of evacuees departing an evacuation
node over time after an evacuation start time $t_0\in\mathcal{H}$. The
number $D_k(t)$ of evacuees departing an evacuation node $k$ at time
$t$ is defined using a selected response curve $f$ as follows:
\begin{equation}
D_k(t)=
\begin{cases}
0&\text{if}\ t<t_0\\
f(t-t_0)&\text{if}\ t\geq t_0
\end{cases}
\end{equation}
\noindent $D_k(t)$ can be used to precisely specify a non-preemptive
evacuation schedule for evacuation node $k$. Figure
\ref{fig:response_curves} shows $D_k(t)$ utilizing four different
types of response curves. The S-shape, Rayleigh and inverse Rayleigh
curves use $t_0$ = 60 minutes, while the step function uses $t_0$ =
120 minutes. The step response curve, where evacuees depart at a
constant rate after $t_0$ until a region is completely evacuated, is
the response curse considered in this paper.

\begin{definition}
  Given an evacuation graph $\GS$, the Non-Preemptive Zone-Based
  Evacuation Planning Problem (NP-ZEPP) consists of finding for each
  evacuation zone an evacuation path to a safe node, a departure time,
  and a response curve that maximize the flow of evacuees to safe
  nodes, while satisfying the problem constraints.
\end{definition}

\noindent
The paper considers the ZEPP, C-ZEPP, and NP-ZEPP problems and
algorithms to solve them.  Each of these problems is considered with
and without contraflows.

\section{The Basic MIP Model for the ZEPP} 
\label{sec:MIP}

This section presents a Mixed Integer Program (MIP) model for the
ZEPP. Model (\ref{lp:od-obj}-\ref{lp:od-end}) in Figure \ref{fig:MIP}
provides the intuition that serves as the basis for the more complex
models presented subsequently. The decision variables of the model are
as follows: Binary variable $\bmcf{e}{k}$ is equal to 1 if and only if
edge $e\in\sAS$ belongs to the evacuation path for evacuation node $k$,
and $\flow[k]{e}$ is a continuous variable equal to the flow of
evacuees from evacuation node $k$ on edge $e\in\sA$. To indicate which
road should be used in contraflows, binary variable $y_e$ represents
whether arc $e$ is used in its normal direction ($y_e=1$) or in
contraflow ($y_e=0$). Each road segment $(e,\bar{e})$ with
$e\in\hat{\mathcal{A}}_c,\bar{e}\in\check{\mathcal{A}}_c$ can then be
utilized in one of three possible configurations: (a)
$(y_e,y_{\bar{e}})=(1,1)$ where both arcs are used in their normal
directions, (b) $(y_e,y_{\bar{e}})=(1,0)$ where arc $\bar{e}$ is used
in contraflow, or (c) $(y_e,y_{\bar{e}})=(0,1)$ where arc $e$ is used
in contraflow.

\begin{figure}[!t]
\begin{align}
\max\quad & \sum_{k\in\sE}\sum_{e\in\nin{v_t}} \flow[k]{e} \label{lp:od-obj}\\
\st& \sum_{e\in\nout{k}}\bmcf{e}{k} = 1 & \forall k \in \sES \label{lp:od-st-onepath} \\
& \sum_{e\in\nin{i}}\bmcf{e}{k} - \sum_{e\in\nout{i}}\bmcf{e}{k} = 0 & \forall k \in \sES, i \in \sTS \label{lp:od-st-flow} \\ 
& \sum_{e\in\nin{i}}\flow[k]{e} - \sum_{e\in\nout{i}}\flow[k]{e} = 0 & \forall i \in\sN\setminus\{v_s,v_t\},  k \in \sES 
\label{lp:od-flow}\\
& \sum_{k\in\sES}\flow[k]{e} \leq u_{e} 		& \forall e\in\sA\setminus\sAC \label{lp:od-ub}\\
& \sum_{k\in\sES}\flow[k]{e} \leq \bincf{e_0}\eu{e} + (1-\bincf{\bar{e_0}})\eu{\bar{e}}   &  \forall e\in\sAC \label{lp:od-ub-cf} \\
& \flow[k]{e} \leq u_{e}*\bmcf{e}{k} 		& \forall e\in\sA,  k \in \sES \label{lp:od-ubp}\\
& \bincf{e} + \bincf{\bar{e}} \geq 1 &  \forall e \in \sACd \label{lp:od-bin-cf}\\
& \bincf{e} \in \{0,1\} & \forall e\in\sAC \\
& \flow[k]{e} \geq 0, \bmcf{e}{k} \in \{0,1\} & \forall e\in\sA,  k \in \sES \label{lp:od-end}
\end{align}
\caption{The MIP Model for the ZEEP.}
\label{fig:MIP}
\end{figure}

Constraints \eqref{lp:od-st-onepath} ensure that exactly one path is
used to route the flow coming from evacuation nodes, while constraints
\eqref{lp:od-st-flow} ensure the continuity of the path.  Constraints
\eqref{lp:od-flow} ensure  flow conservation through the
time-expanded graph.  Constraints \eqref{lp:od-ub} enforce the
capacity of each edge in the time-expanded graph.  Constraints
\eqref{lp:od-ub-cf} enforce the capacity on edges that allow
contraflow: They allocate to $e$ the capacity of edge $\bar{e}$
whenever $\bar{e}$ is used in contraflow, and forbid any flow on $e$
when it is used in contraflow.  Constraints \eqref{lp:od-ubp} ensure
that there is no flow of evacuees coming from an evacuation node $k$ if
edge $e$ is not part of the evacuation path for $k$.  Constraints
\eqref{lp:od-bin-cf} prohibit the simultaneous use of $e$ and
$\bar{e}$ in contraflow.  The objective \eqref{lp:od-obj} maximizes
the number of evacuees reaching safety.

Model (\ref{lp:od-obj}-\ref{lp:od-end}) is intractable for the case
study used in this paper which has approximately 30,000 nodes and
75,000 arcs. The computational difficulty comes from two
interconnected components: The selection of the paths is a design
component, while the scheduling of the evacuation is multi-commodity
flow. {\em The algorithms described in this paper address this
  computational challenge by separating these two aspects.} Observe
also that the temporal aspects (i.e., when to schedule evacuees along
a path) are an important and difficult aspect of the ZEPP. Finally,
it is interesting to mention that additional requirements, such as
convergent evacuations and non-preemption, may lead to elegant
computational contributions that would not be possible otherwise.

%

\section{Benders Decomposition for the ZEEP}
\label{sec:BN}

This section presents an (approximate) Benders decomposition for the
ZEEP. This Benders decomposition is referred to as the Benders
Non-convergent (BN) method in the rest of the paper. The Restricted
Master Problem (RMP) of the Benders decomposition selects evacuation
paths that are then used in the subproblem (SP) for scheduling the
flows of evacuees over time along these paths. It is important to note
that the subproblem is a multi-commodity flow and hence it is not
totally unimodular. As a result, the Benders decomposition in the BN
method solves a relaxation of the ZEEP where the integrality
constraints on flow variables are relaxed. A final MIP is solved to
obtain an integer solution to the subproblem.

As is traditional in Benders decomposition, the objective values of
the RMP and SP provide upper and lower bounds on the optimal solution
to Model (\ref{lp:od-obj}-\ref{lp:od-end}) without integrality
constraints on the flow variables. When they converge, evacuation
paths from the RMP and the evacuation schedule from the SP form an
optimal solution. Otherwise, a Benders cut is generated from the
solution of the SP and introduced into the RMP as an additional
constraint to remove the current evacuation paths from the RMP's feasible
region, after which the entire process is repeated.

\paragraph{The Restricted Master Problem}

The RMP, depicted in Figure \ref{fig:BN:RMP}, finds evacuation paths
for each evacuation zone. It operates on the static graph and its main
decision variables are the binary variables $x_{e,k}$ of Model
(\ref{lp:od-obj}-\ref{lp:od-end}). {\em In addition to the Benders
  cuts, the RMP also reasons about aggregate flows and aggregated
  capacities, an idea that was proposed in \cite{romanski2016} to obtain
  reasonable evacuation paths early on.} In particular, variable
$\psi_{e,k}$ represents the aggregate flow of evacuees from evacuation
node $k$ along arc $e$ and arc capacities are aggregated over the time
horizon in all of the RMP's capacity constraints. Finally,
$z_{\text{RMP}}$ is the RMP's objective value and represents the
number of evacuees reaching safety.

\begin{figure}[!tbhp]
\begin{equation}\label{eqn:bn_obj1}
\max z_{\text{RMP}}
\end{equation}
\begin{equation*}
\text{subject to}
\end{equation*}
\begin{equation}\label{eqn:bn_obj2}
z_{\text{RMP}}\leq\sum_{k\in\mathcal{E}}\sum_{e\in\delta^{+}(k)}\psi_{e,k}
\end{equation}
\begin{equation}\label{eqn:bn_onepathperevac}
\sum_{e\in\delta^{+}(k)}x_{e,k}= 1 \qquad \forall k\in\mathcal{E}
\end{equation}
\begin{equation}\label{eqn:bn_onepathpreserve}
\sum_{e\in\delta^{+}(i)}x_{e,k}\leq 1\qquad\forall i\in\mathcal{T},\forall k\in\mathcal{E}
\end{equation}
\begin{equation}\label{eqn:bn_flowconserve}
\sum_{e\in\delta^{-}(i)}\psi_{e,k}-\sum_{e\in\delta^{+}(i)}\psi_{e,k}=0\qquad\forall i\in\mathcal{T},\forall k\in\mathcal{E}
\end{equation}
\begin{equation}\label{eqn:bn_demand}
\sum_{e\in\delta^{+}(k)}\psi_{e,k}\leq d_k\qquad\forall k\in\mathcal{E}
\end{equation}
\begin{equation}\label{eqn:bn_flowonselect1}
\psi_{e,k}\leq x_{e,k}\sum_{t\in\mathcal{H}}u_{e_t}\qquad \forall e\in\mathcal{A}\setminus\mathcal{A}_c,\forall k\in\mathcal{E}
\end{equation}
\begin{equation}\label{eqn:bn_flowonselect2}
\psi_{e,k}\leq x_{e,k}\sum_{t\in\mathcal{H}}(u_{e_t}+u_{\bar{e}_t})\qquad \forall e\in\mathcal{A}_c,\forall k\in\mathcal{E}
\end{equation}
\begin{equation}\label{eqn:bn_totalcapacity1}
\sum_{k\in\mathcal{E}}\psi_{e,k}\leq\sum_{t\in\mathcal{H}}u_{e_t}\qquad \forall e\in\mathcal{A}\setminus\mathcal{A}_c
\end{equation}
\begin{equation}\label{eqn:bn_totalcapacity2}
\sum_{k\in\mathcal{E}}\psi_{e,k}\leq y_e\sum_{t\in\mathcal{H}}u_{e_t}+(1-y_{\bar{e}})\sum_{t\in\mathcal{H}}u_{\bar{e}_t}\qquad \forall e\in\mathcal{A}_c
\end{equation}
\begin{equation}\label{eqn:bn_contraflow}
y_e+y_{\bar{e}}\geq 1\qquad\forall e\in\mathcal{A}_c
\end{equation}
\begin{equation}
\psi_{e,k}\geq 0\qquad\forall e\in\mathcal{A},\forall k\in\mathcal{E}
\end{equation}
\begin{equation}
x_{e,k}\in\{0,1\}\qquad\forall e\in\mathcal{A},\forall k\in\mathcal{E}
\end{equation}
\begin{equation}
y_e\in\{0,1\}\qquad\forall e\in\mathcal{A}_c
\end{equation}
\caption{The Restricted Master Problem for the BN method.}
\label{fig:BN:RMP}
\end{figure}

Constraints \eqref{eqn:bn_obj2}, together with objective function
\eqref{eqn:bn_obj1}, maximize the flow of evacuees from all evacuation
nodes. Constraints \eqref{eqn:bn_onepathperevac} specify that exactly
one path is generated for each evacuation node, and constraints
\eqref{eqn:bn_onepathpreserve} ensure that the one path requirement
is preserved throughout the graph. Constraints
\eqref{eqn:bn_flowconserve} ensure that flow is conserved throughout
the graph, and constraints \eqref{eqn:bn_demand} make sure that total
flow from each evacuation node does not exceed its demand. Constraints
\eqref{eqn:bn_flowonselect1} and \eqref{eqn:bn_flowonselect2} permit
evacuee flow from evacuation node $k$ on an arc only if the arc is
selected for the evacuation path of $k$. Constraints
\eqref{eqn:bn_totalcapacity1} ensure that the flow from all evacuation
nodes along an arc does not exceed the aggregate capacity for arcs
that may not be used in contraflow, while constraints
\eqref{eqn:bn_totalcapacity2} do the same for arcs that may be used in
contraflow. Finally, constraints \eqref{eqn:bn_contraflow} indicate
that at most one arc in road segment $(e,\bar{e})$ with
$e\in\hat{\mathcal{A}}_c,\bar{e}\in\check{\mathcal{A}}_c$ can be used
in contraflow. To generate evacuation plans with contraflows,
constraints \eqref{eqn:bn_contraflow} can be replaced with constraints
\begin{equation*}
  \label{eqn:bn_nocontraflow}
y_e+y_{\bar{e}}=2\qquad\forall e\in\mathcal{A}_c
\end{equation*}
which forces all arcs to be used in their normal directions.

\paragraph{The Benders Subproblem}

The SP, depicted in Figure \ref{fig:BN:SP}, utilizes paths generated
from the RMP together with the time-expanded graph $\mathcal{G}^x$ to
generate an evacuation schedule that maximizes the number of evacuees
reaching safety along those paths. The paths are specified by the
values $\{\bar{x}_{e,k}\}$ and $\{\bar{y}_e\}$ for variables
$\{x_{e,k}\}$ and $\{y_e\}$ in the RMP. The SP uses variable
$\varphi_{e_t,k}$ to represent the flow of evacuees from evacuation
node $k$ along arc $e_t$ in $\mathcal{G}^x$, and $z_{\text{SP}}$ is
the SP's objective value.

\begin{figure}[!t]
\begin{equation}\label{eqn:bn_sp_obj}
  \max z_{\text{SP}}=\sum_{k\in\mathcal{E}}\sum_{e_t\in\delta^{+}(k)}\varphi_{e_t,k}
\end{equation}
\begin{equation*}
\text{subject to}
\end{equation*}
\begin{equation}\label{eqn:bn_sp_flowconserve}
\sum_{e_t\in\delta^{-}(i)}\varphi_{e_t,k}-\sum_{e_t\in\delta^{+}(i)}\varphi_{e_t,k}=0\qquad\forall i\in\mathcal{T}^{x},\forall k\in\mathcal{E}
\end{equation}
\begin{equation}\label{eqn:bn_sp_demand}
\sum_{t\in\mathcal{H}}\sum_{e_t\in\delta^{+}(k)}\varphi_{e_t,k}\leq d_k\qquad\forall k\in\mathcal{E}
\end{equation}
\begin{equation}\label{eqn:bn_sp_flowonselect1}
\varphi_{e_t,k}\leq \bar{x}_{e,k}\cdot u_{e_t}\qquad\forall e\in\mathcal{A}\setminus\mathcal{A}_c,\forall t\in\mathcal{H},\forall k\in\mathcal{E}
\end{equation}
\begin{equation}
\label{eqn:bn_sp_flowonselect2}
\varphi_{e_t,k}\leq \bar{x}_{e,k}\cdot(u_{e_t}+u_{\bar{e}_t})\qquad\forall e\in\mathcal{A}_c,\forall t\in\mathcal{H},\forall k\in\mathcal{E}
\end{equation}
\begin{equation}\label{eqn:bn_sp_capacity1}
\sum_{k\in\mathcal{E}}\varphi_{e_t,k}\leq u_{e_t}\qquad\forall e\in\mathcal{A}\setminus\mathcal{A}_c,\forall t\in\mathcal{H}
\end{equation}
\begin{equation}\label{eqn:bn_sp_capacity2}
\sum_{k\in\mathcal{E}}\varphi_{e_t,k}\leq \bar{y}_e\cdot u_{e_t}+(1-\bar{y}_{\bar{e}})\cdot u_{\bar{e}_t}\qquad\forall e\in\mathcal{A}_c,\forall t\in\mathcal{H}
\end{equation}
\begin{equation}
\varphi_{e_t,k}\geq 0\qquad\forall e_t\in\mathcal{A}^{x},\forall k\in\mathcal{E}
\end{equation}
\caption{The Bender Subproblem for the BN method.}
\label{fig:BN:SP}
\end{figure}

The objective function \eqref{eqn:bn_sp_obj} maximizes the flow of
evacuees from all evacuation nodes. Constraints
\eqref{eqn:bn_sp_flowconserve} enforce flow conservation throughout
$\mathcal{G}^x$, while constraints \eqref{eqn:bn_sp_demand} ensure
that the total flow from each evacuation node does not exceed its
demand. Constraints \eqref{eqn:bn_sp_flowonselect1} and
\eqref{eqn:bn_sp_flowonselect2} permit flow only on the selected arcs
for each evacuation node. Constraints \eqref{eqn:bn_sp_capacity1}
ensure that the total flow from all evacuation nodes along an arc does
not exceed its capacity for arcs that may not be used in contraflow,
while constraints \eqref{eqn:bn_sp_capacity2} do the same for arcs
that may be used in contraflow.

\paragraph{The Benders Cuts}

A Benders optimality cut is generated from the solution of the SP and
added to the RMP as long as the objective values of the RMP and SP do
not converge. The cut is of the form
\begin{equation}
\begin{split}\label{eqn:bn_benderscut}
z_{\text{RMP}}\leq&\sum_{k\in\mathcal{E}}d_k\cdot\pi_k+\sum_{e\in\mathcal{A}\setminus\mathcal{A}_c}\sum_{t\in\mathcal{H}}u_{e_t}\sum_{k\in\mathcal{E}}x_{e,k}\cdot\pi_{e_t,k}+\\
&\sum_{e\in\mathcal{A}_c}\sum_{t\in\mathcal{H}}(u_{e_t}+u_{\bar{e}_t})\sum_{k\in\mathcal{E}}x_{e,k}\cdot\pi_{e_{t,c},k}+\sum_{e\in\mathcal{A}\setminus\mathcal{A}_c}\sum_{t\in\mathcal{H}}u_{e_t}\cdot\pi_{e_t}+\\
&\sum_{e\in\mathcal{A}_c}\sum_{t\in\mathcal{H}}[y_e\cdot u_{e_t}+(1-y_{\bar{e}})\cdot u_{\bar{e}_t}]\cdot\pi_{e_{t,c}}
\end{split}
\end{equation}
where $\{\pi_k\}$, $\{\pi_{e_{t},k}\}$, $\{\pi_{e_{t,c},k}\}$,
$\{\pi_{e_t}\}$, and $\{\pi_{e_{t,c}}\}$ are the dual variables of
constraints \eqref{eqn:bn_sp_demand}, \eqref{eqn:bn_sp_flowonselect1},
\eqref{eqn:bn_sp_flowonselect2}, \eqref{eqn:bn_sp_capacity1}, and
\eqref{eqn:bn_sp_capacity2} respectively.  Since the SP is always
feasible, Benders feasibility cuts are never generated by this
algorithm.

\paragraph{The Benders Non-convergent Algorithm}

Algorithm \ref{alg:bendersnonconv} summarizes the entire BN
algorithm which uses $\text{RMP}(\mathcal{G},\mathcal{H})$ to denote an
optimal solution obtained from solving the RMP given static graph
$\mathcal{G}$ and time horizon $\mathcal{H}$ as inputs,
$\text{SP}(\Psi,\mathcal{H})$ to denote an optimal solution of the SP
given a solution to the RMP, $\Psi$, and time horizon $\mathcal{H}$ as
inputs, and $z(\sigma)$ to denote the objective value of a solution
$\sigma$.

\begin{algorithm*}[!t]
	\caption{Benders Non-convergent}\label{alg:bendersnonconv}
	\begin{algorithmic}[1]
		\State $t^*\leftarrow\min\{t\in\mathcal{H}\,|\,z(\text{RMP}(\mathcal{G},[0..t]))=z(\text{RMP}(\mathcal{G},\mathcal{H}))\}$
		\State $z_{\text{RMP}}\leftarrow z(\text{RMP}(\mathcal{G},[0..t^*]))$
		\State $z_{\text{SP}}\leftarrow z(\text{SP}(\text{RMP}(\mathcal{G},[0..t^*]),\mathcal{H}))$
		\State $z_{\text{SP,max}}\leftarrow z_{\text{SP}}$
		\While{$z_{\text{RMP}}-z_{\text{SP,max}}>\varepsilon$}
			\State \text{Generate Benders cut from solution of SP and add it to RMP}
			\State $z_{\text{RMP}}\leftarrow z(\text{RMP}(\mathcal{G},\mathcal{H}))$
			\State $z_{\text{SP}}\leftarrow z(\text{SP}(\text{RMP}(\mathcal{G},\mathcal{H}),\mathcal{H}))$
			\State $z_{\text{SP,max}}\leftarrow\max\{z_{\text{SP,max}},z_{\text{SP}}\}$
		\EndWhile
		\State Solve $\text{SP}(\text{RMP}(\mathcal{G},\mathcal{H}),\mathcal{H})$ with $\varphi_{e_t,k}\ \text{integer}\ \forall e_t\in\mathcal{A}^{x},\forall k\in\mathcal{E}$
		\State \textbf{return} Evacuation paths from solution of RMP and evacuation schedule from solution of SP
	\end{algorithmic}
\end{algorithm*}

The BN algorithm begins with a procedure that searches for the
tightest time horizon $t^*$ that can preserve the optimal solution to
the RMP, $z(\text{RMP}(\mathcal{G},\mathcal{H}))$. This step was
originally proposed by Even \etal \cite{even2015} who found that a
tighter time horizon produces better evacuation paths for the flow
scheduling problem of their two-stage approach. The BN
method adopts a similar strategy to seed the Benders
decomposition. The procedure is implemented using a simple sequential
search which solves $\text{RMP}(\mathcal{G},\mathcal{H})$ with
progressively smaller values of $\mathcal{H}$ in search of $t^*$.

After this step, the algorithm proceeds to first solve the RMP to
generate evacuation paths, and then the SP using the generated paths
as input to generate an evacuation schedule. The minimum objective
value $z_{\text{RMP}}$ of the RMP is then compared to the maximum
objective value $z_{\text{SP,max}}$ of the SP. If they do not converge
(if their difference is larger than a convergence criterion
$\varepsilon$, which we set to 0), a Benders cut is generated
utilizing dual variables from the solution of the SP and added to the
RMP to remove the current evacuation paths from its feasible region. The
process of solving the RMP and SP is then repeated until convergence.
Since the Benders decomposition relaxes the integrality constraints on
the flow variables, the subproblem is solved one more time as a MIP
after convergence to obtain an integer solution to the subproblem.

\section{Benders Decomposition for Convergent Evacuation Planning}
\label{sec:BC}

This section present the Benders decomposition of \cite{romanski2016}
for the C-ZEEP, i.e., for the convergent preemptive zone-based
evacuation planning. This Benders decomposition is referred to as the
Benders convergent (BC) method in this paper. The BC method shares a
lot of similarities with the BN method. The BC method however imposes
that the evacuation paths form a convergent graph and hence that the
outdegree at each transit node is at most 1. {\em This convergence property
has some fundamental consequences: (1) the BC method is exact since
the subproblem becomes totally submodular; (2) it trivially supports
contraflows; and (3) its computational performance is strong compared to
all the other algorithms.}

\paragraph{The Restricted Master Problem}

\begin{figure}[!t]
\begin{equation}\label{eqn:bc_obj1}
\max z_{\text{RMP}}
\end{equation}
\begin{equation*}
\text{subject to}
\end{equation*}
\begin{equation}\label{eqn:bc_obj2}
z_{\text{RMP}}\leq\sum_{k\in\mathcal{E}}\sum_{e\in\delta^{+}(k)}\psi_e
\end{equation}
\begin{equation}\label{eqn:bc_flowconserve}
\sum_{e\in\delta^{-}(i)}\psi_e-\sum_{e\in\delta^{+}(i)}\psi_e=0\qquad\forall i\in\mathcal{T}
\end{equation}
\begin{equation}\label{eqn:bc_convergence}
\sum_{e\in\delta^{+}(i)}x_e\leq 1\qquad \forall i\in\mathcal{E}\cup\mathcal{T}
\end{equation}
\begin{equation}\label{eqn:bc_capacity}
\psi_e\leq x_e\sum_{t\in\mathcal{H}}u_{e_t}\qquad \forall e\in\mathcal{A}
\end{equation}
\begin{equation}\label{eqn:bc_demand}
\sum_{e\in\delta^{+}(k)}\psi_e\leq d_k\qquad \forall k\in\mathcal{E}
\end{equation}
\begin{equation}
\psi_e\geq 0\qquad\forall e\in\mathcal{A}
\end{equation}
\begin{equation}
x_e\in\{0,1\}\qquad\forall e\in\mathcal{A}
\end{equation}
\caption{The Restricted Master Problem for the BC method.}
\label{fig:BC:RMP}
\end{figure}

The RMP for the BC method is presented in Figure
\ref{fig:BC:RMP}. {\em Because the paths are convergent, the model is
considerably simpler: There is no need to track the origin of the flow
(i.e., the evaluation zone) and have different flow conservation
constraints for each evacuation zone.} The model still uses a binary
variable $x_e$ to indicate if an arc $e$ is to be part of an
evacuation path. But it uses a single variable $\psi_e$ to represent
the aggregate flow of evacuees along arc $e$ over the time horizon. Arc
capacities are again aggregated over the time horizon in the capacity
constraints. Constraint \eqref{eqn:bc_obj2} combined with objective function
\eqref{eqn:bc_obj1} maximizes the flow of evacuees leaving all
evacuation nodes. Constraint \eqref{eqn:bc_flowconserve} enforces flow
conservation, while constraint \eqref{eqn:bc_convergence} enforces the
convergence of arcs selected by the evacuation paths. Constraints
\eqref{eqn:bc_capacity} permit flows only on selected arcs and ensures
aggregate flow along them do not exceed their aggregate
capacity. Finally, constraint \eqref{eqn:bc_demand} ensures that the
total flow from each evacuation node does not exceed its demand.

\paragraph{The Benders Subproblem}

The Benders problem for method BC, depicted in Figure \ref{fig:BC:SP},
is again simpler due to path convergence and uses a a variable
$\varphi_{e_t}$ to represent flow of evacuees along arc $e_t$ in
$\mathcal{G}^x$. Objective function \eqref{eqn:bc_sp_obj} maximizes
flow of evacuees across all evacuation nodes. Constraints
\eqref{eqn:bc_sp_flowconserve} enforce flow conservation, constraints
\eqref{eqn:bc_sp_capacity} permit flow only on arcs selected for
evacuation paths and ensure the flow does not exceed the arc's
capacity, and constraints \eqref{eqn:bc_sp_demand} ensure that the
total flow leaving each evacuation node does not exceed its demand.

\begin{figure}[!t]
\begin{equation}\label{eqn:bc_sp_obj}
\max z_{\text{SP}}=\sum_{k\in\mathcal{E}}\sum_{e_t\in\delta^{+}(k)}\varphi_{e_t}
\end{equation}
\begin{equation*}
\text{subject to}
\end{equation*}
\begin{equation}\label{eqn:bc_sp_flowconserve}
\sum_{e_t\in\delta^{-}(i)}\varphi_{e_t}-\sum_{e_t\in\delta^{+}(i)}\varphi_{e_t}=0\qquad \forall i\in\mathcal{T}^{x}
\end{equation}
\begin{equation}\label{eqn:bc_sp_capacity}
\varphi_{e_t}\leq \bar{x}_e\cdot u_{e_t}\qquad \forall e\in\mathcal{A},\forall t\in\mathcal{H}
\end{equation}
\begin{equation}\label{eqn:bc_sp_demand}
\sum_{e_t\in\delta^{+}(k)}\varphi_{e_t}\leq d_k\qquad \forall k\in\mathcal{E}
\end{equation}
\begin{equation}
\varphi_{e_t}\geq 0\qquad \forall e_t\in\mathcal{A}^{x}
\end{equation}
\caption{The Bender Subproblem for the BN method.}
\label{fig:BC:SP}
\end{figure}

\paragraph{The Benders Cuts}

The Benders optimality cuts are of the form
\begin{equation}\label{eqn:bc_benderscut}
z_{\text{RMP}}\leq\sum_{e\in\mathcal{A}}x_e\sum_{t\in\mathcal{H}}u_{e_t}\cdot \pi_{e_t}+\sum_{k\in\mathcal{E}}d_k\cdot \pi_k
\end{equation}
and use dual variables $\{\pi_{e_t}\}$ and $\{\pi_k\}$ associated with
constraints \eqref{eqn:bc_sp_capacity} and \eqref{eqn:bc_sp_demand} of
the SP respectively. Again, Benders feasibility cuts are never
generated in this algorithm because the SP is always feasible.

\paragraph{Pareto-Optimal Cuts}

The convergence of Benders decomposition can be accelerated through
utilization of Pareto-optimal cuts \cite{romanski2016}, i.e., cuts
that are not dominated by any other Benders cut. The Magnanti-Wong
method \cite{magnanti1981} is utilized to generate these stronger
cuts. The method requires a core point, i.e., a point located within
the relative interior of the convex hull of the feasibility domain of
the RMP's first-stage variable $\{x_e\}$. For this formulation, the
core point utilized is simply $x_e^0=\frac{1}{|\delta^{+}(i)+1|}$ for
each arc $e=(i,j)$. The dual of the Magnanti-Wong problem (DMWP) which
utilizes this core point and the optimal objective value of the SP,
$z_\text{SP}$, is solved to generate a Pareto-optimal cut.

\begin{equation}
\max\sum_{k\in\mathcal{E}}\sum_{e_t\in\delta^{+}(k)}\varphi_{e_t}+\xi\cdot z_\text{SP}
\end{equation}
\begin{equation*}
\text{subject to}
\end{equation*}
\begin{equation}
\sum_{e_t\in\delta^{-}(i)}\varphi_{e_t}-\sum_{e_t\in\delta^{+}(i)}\varphi_{e_t}=0\qquad \forall i\in\mathcal{T}^x
\end{equation}
\begin{equation}\label{eqn:mw_capacity}
\varphi_{e_t}+x_e\cdot u_{e_t}\cdot \xi\leq x_e^0\cdot u_{e_t}\qquad \forall e\in\mathcal{A},\forall t\in\mathcal{H}
\end{equation}
\begin{equation}\label{eqn:mw_demand}
\sum_{e_t\in\delta^{+}(k)}\varphi_{e_t}+d_k\cdot\xi\leq d_k\qquad \forall k\in\mathcal{E}
\end{equation}
\begin{equation}
\varphi_{e_t}\geq 0\qquad \forall e_t\in\mathcal{A}^{x}
\end{equation}

\noindent In order to generate the Pareto-optimal cut, coefficients
$\{\pi_{e_t}\}$ and $\{\pi_k\}$ in cut \eqref{eqn:bc_benderscut} are
taken from the dual variables of constraints \eqref{eqn:mw_capacity}
and \eqref{eqn:mw_demand} respectively instead of those from
constraints of the SP.

\paragraph{Contraflow Extension}

The BC algorithm proposed in \cite{romanski2016} did not consider
contraflows but it can be easily extended to support this
functionality. {\em In fact, convergent evacuations make contraflows
  very easy}, as their tree structure guarantees that, for any road
segment $(e,\bar{e})$ with
$e\in\hat{\mathcal{A}}_c,\bar{e}\in\check{\mathcal{A}}_c$, if $x_e =
1$, then $x_{\bar{e}}=0$. In other words, if an arc
$e\in\hat{\mathcal{A}}_c$ is in an evacuation path, the corresponding
unique arc in the opposite direction $\bar{e}$ is not. This makes it
possible to use arc $\bar{e}$ in contraflow if arc
$e\in\hat{\mathcal{A}}_c$ is being used in an evacuation plan, since
arc $\bar{e}$ is guaranteed not to be part of any other evacuation
path by the convergence constraint. As a consequence, the BC algorithm
is extended to allow contraflow as follows. Before the algorithm is
executed, the capacities of all arcs $e\in\mathcal{A}_c$ is replaced
with new capacities $u_{e_t,\text{new}}=u_{e_t}+u_{\bar{e}_t}$. To
identify where to use contraflows, it suffices to identify arcs
$e\in\mathcal{A}_c$ with flows $\varphi_{e_t}>u_{e_t}$, meaning that
the extra capacity afforded by using arc $\bar{e}$ in contraflow is
necessary to achieve optimality.

\section{The Conflict-based Path Generation Method}
\label{sec:CPG}

This section summarizes the heuristic algorithm originally presented
in \cite{pillac2016} to solve the ZEPP and called the Conflict-based
Path Generation (CPG) method. The CPG method originated from an
attempt to design a column-generation algorithm for the ZEPP. However,
each new path creates a collection of variables, i.e., the path
variables and the associated flow variables, and these variables are
linked like in constraints \eqref{lp:od-ubp} of the MIP model. Since
the duals of these constraints are not readily available, it did not
appear easy to derive a column-generation algorithm at the
time. Hence, the CPG mimics the behavior of a column-generation
algorithm but its pricing subproblem is a heuristic. More precisely,
the CPG breaks down the evacuation planning problem into a subproblem
(SP) responsible for generating evacuation paths and a restricted
master problem (RMP) responsible for selecting paths and scheduling
the evacuation.  The method maintains a subset of critical evacuation
nodes $\mathcal{E}'\subseteq\mathcal{E}$, i.e., evacuation nodes that have not
been fully evacuated, and it alternates execution of the SP and RMP
until $\mathcal{E}'$ is empty.

\paragraph{The Restricted Master Problem}

The RMP of the CPG method, shown in Figure \ref{fig:CPG:RMP}, selects
an evacuation path for each evacuation node and schedules the evacuees
over them to maximize the number of evacuees reaching safety. The
paths are selected from a set of evacuation paths $\Omega'$ generated
by the SP. The RMP uses a binary variable $x_p$ to indicate whether a
path $p\in\Omega'$ is selected for the evacuation plan, a continuous
variable $\varphi_p^t$ to represent the number of evacuees departing
along path $p$ at departure time $t$, and a continuous variable
$\bar{\varphi}_k$ to represent the number of evacuees that cannot be
evacuated at evacuation node $k$. In addition to these,
$\Omega_k\subset\Omega'$ is the subset of evacuation paths for
evacuation node $k$, $\omega(e)\subseteq\Omega'$ is the subset of
paths that contain arc $e$, $\mathcal{H}_p\subseteq\mathcal{H}$ is the
subset of time steps over which path $p$ is usable, $\tau_p^e$ is the
number of time steps required to reach arc $e$ when traversing path
$p$, and $u_p$ is the capacity of path $p$.

\begin{figure}[!t]
\begin{equation}\label{eqn:cpg_mp_obj}
\max\sum_{p\in\Omega}\sum_{t\in\mathcal{H}_p}\varphi_p^t
\end{equation}
\begin{equation*}
\text{subject to}
\end{equation*}
\begin{equation}\label{eqn:cpg_mp_onepath}
\sum_{p\in\Omega_k}x_p=1\qquad\forall k\in\mathcal{E}
\end{equation}
\begin{equation}\label{eqn:cpg_mp_demand}
\sum_{p\in\Omega_k}\sum_{t\in\mathcal{H}_p}\varphi_p^t+\bar{\varphi}_k=d_k\qquad\forall k\in\mathcal{E}
\end{equation}
\begin{equation}\label{eqn:cpg_mp_capacity1}
\sum_{\substack{p\in\omega(e)\\t-\tau_p^e\in\mathcal{H}_p}}\varphi_p^{t-\tau_p^e}\leq u_{e_t}\qquad\forall e\in\mathcal{A}\setminus\mathcal{A}_c,\forall t\in\mathcal{H}
\end{equation}
\begin{equation}\label{eqn:cpg_mp_capacity2}
\sum_{\substack{p\in\omega(e)\\t-\tau_p^e\in\mathcal{H}_p}}\varphi_p^{t-\tau_p^e}\leq y_e\cdot u_{e_t}+(1-y_{\bar{e}})\cdot u_{\bar{e}_t}\qquad\forall e\in\mathcal{A}_c,\forall t\in\mathcal{H}
\end{equation}
\begin{equation}\label{eqn:cpg_mp_contraflow1}
y_e+y_{\bar{e}}\geq 1\qquad\forall e\in\mathcal{A}_c
\end{equation}
\begin{equation}\label{eqn:cpg_mp_pathcap}
\sum_{t\in\mathcal{H}_p}\varphi_p^t\leq |\mathcal{H}_p|\cdot x_p\cdot u_p\qquad\forall p\in\Omega'
\end{equation}
\begin{equation}
\varphi_p^t\geq 0\qquad\forall p\in\Omega',\forall t\in\mathcal{H}_p
\end{equation}
\begin{equation}
\bar{\varphi}_k\geq 0\qquad\forall k\in\mathcal{E}
\end{equation}
\begin{equation}
y_e\in\{0,1\}\qquad\forall e\in\mathcal{A}_c
\end{equation}
\begin{equation}
x_p\in\{0,1\}\qquad\forall p\in\Omega'
\end{equation}
\caption{The Restricted Master Problem for the CPG Method.}
\label{fig:CPG:RMP}
\end{figure}

Objective function \eqref{eqn:cpg_mp_obj} maximizes the flow of
evacuees over all paths. Constraints \eqref{eqn:cpg_mp_onepath} allow
only one path from being selected per evacuation node, while
constraints \eqref{eqn:cpg_mp_demand} ensures that the sum of evacuees
who reach or who do not reach safety is equal to the demand at each
evacuation node. Constraints \eqref{eqn:cpg_mp_capacity1} and
\eqref{eqn:cpg_mp_capacity2} enforce the capacity of arcs that may not
and may be used in contraflow respectively. Constraints
\eqref{eqn:cpg_mp_contraflow1} prohibit the simultaneous use of arcs
$e$ and $\bar{e}$ in contraflow for road segment $(e,\bar{e})$ with
$e\in\hat{\mathcal{A}}_c,\bar{e}\in\check{\mathcal{A}}_c$. Finally,
constraints \eqref{eqn:cpg_mp_pathcap} allows for flows only on
selected paths. To generate an evacuation plan that does not permit
contraflow, Constraint \eqref{eqn:cpg_mp_contraflow1} is replaced with
Constraint \eqref{eqn:cpg_mp_contraflow2} to ensure all arcs are only
used in their normal directions.
\begin{equation}\label{eqn:cpg_mp_contraflow2}
y_e+y_{\bar{e}}=2\qquad\forall e\in\mathcal{A}_c
\end{equation}
Observe Constraints \eqref{eqn:cpg_mp_pathcap} that feature both
variables $\varphi_p^t$ and $x_p$. These constraints must be generated
every time a new path is available and they make it difficult to
obtain a traditional pricing subproblem since their duals are not
available.

\paragraph{The Path Generation Subproblem}

The SP utilizes a conflict-based path generation heuristic to generate
evacuation paths that could potentially improve the objective value of
the RMP. These paths are generated by solving the following
multiple-origins, multiple-destinations shortest path problem:
\begin{equation}\label{eqn:cpg_sp_obj}
\min\sum_{k\in\mathcal{E}'}\sum_{e\in\mathcal{A}}c_e\cdot y_{e,k}
\end{equation}
\begin{equation*}
\text{subject to}
\end{equation*}
\begin{equation}\label{eqn:cpg_sp_pathcontinuity}
\sum_{e\in\delta^{-}(i)}y_{e,k}-\sum_{e\in\delta^{+}(i)}y_{e,k}=0\qquad\forall i\in\mathcal{T},\forall k\in\mathcal{E}'
\end{equation}
\begin{equation}\label{eqn:cpg_sp_onepath}
\sum_{e\in\delta^{+}(k)}y_{e,k}=1\qquad\forall k\in\mathcal{E}'
\end{equation}
\begin{equation}
y_{e,k}\in\{0,1\}\qquad\forall e\in\mathcal{A},\forall k\in\mathcal{E}'
\end{equation}

\noindent
The problem formulation utilizes a binary variable $y_{e,k}$ to
indicate whether arc $e$ belongs to the path generated for evacuation
node $k$. Constraints \eqref{eqn:cpg_sp_pathcontinuity} enforce
path continuity while Constraints \eqref{eqn:cpg_sp_onepath}
ensure only one path is generated for each critical node. Objective
function \eqref{eqn:cpg_sp_obj} minimizes the total cost of all paths. Arc cost $c_e$ is defined as a linear combination of an arc's travel time $s_e$, the number of times arc $e$ is utilized in the current set of
paths $\Omega'$, and the utilization of arc $e$ in the current solution:
\begin{equation}\label{eqn:cpg_sp_edgecost}
c_e=\alpha_t\frac{s_e\cdot r}{\max_{e\in\mathcal{A}}s_e}+\alpha_c\frac{\sum_{\substack{p\in\Omega'\\e\in p}}1}{|\Omega'|}+\alpha_u\frac{\sum_{\substack{p\in\Omega'\\e\in p}}\sum_{t\in\mathcal{H}_p}\varphi_p^t}{\sum_{t\in\mathcal{H}}u_{e_t}}
\end{equation}
\noindent In Equation \eqref{eqn:cpg_sp_edgecost}, $\alpha_t$,
$\alpha_c$, and $\alpha_u$ are positive weights which sum to 1, and
$r$ is a random noise factor that is initialized to 1 and subsequently
modified to $r\in[1-\epsilon,1+\epsilon]$ depending on the number of
iterations in which the objective value of the RMP did not
improve. The value $\epsilon$ is set to 0.50 in this study.

\paragraph{The Conflict-based Path Generation Algorithm}
\label{sec:cpg_algorithm}

This CPG algorithm is summarized in Algorithm
\ref{alg:conflictpathgen}. $\text{PathGenerationSP}$
$(\mathcal{E}',\Omega',\Lambda)$ denotes a subroutine that solves the
SP with a set of critical evacuation nodes $\mathcal{E}'$, a set of
evacuation paths $\Omega'$, and an evacuation schedule $\Lambda$
obtained from the solution of the RMP as
inputs. $\text{EvacuationSchedulingMP}(\Omega')$ denotes a subroutine
that solves the RMP using a set of evacuation paths $\Omega'$ as
input. $\text{FindCriticalEvacuationNodes}(\Lambda)$, as its name
suggests, is a subroutine that identifies evacuation nodes that have
not been fully evacuated using an evacuation schedule $\Lambda$ as
input.

\begin{algorithm*}[!t]
	\caption{Conflict-based Path Generation}\label{alg:conflictpathgen}
	\begin{algorithmic}[1]
		\State $\Omega'\leftarrow\text{PathGenerationSP}(\mathcal{E},\text{\O},\text{\O})$
		\State $\Lambda\leftarrow\text{EvacuationSchedulingMP}(\Omega')$
		\State $\mathcal{E}'\leftarrow\text{FindCriticalEvacuationNodes}(\Lambda)$
		\While{$\mathcal{E}'\neq\text{\O}$}
			\State $\Omega'\leftarrow\Omega'\cup\text{PathGenerationSP}(\mathcal{E}',\Omega',\Lambda)$
			\State $\Lambda\leftarrow\text{EvacuationSchedulingMP}(\Omega')$
			\State $\mathcal{E}'\leftarrow\text{FindCriticalEvacuationNodes}(\Lambda)$
		\EndWhile
		\State	$\Lambda\leftarrow\text{Solve EvacuationSchedulingMP}(\Omega')$ with $\varphi_p^t\ \text{integer}\ \forall p\in\Omega',\forall t\in\mathcal{H}_p$
		\State \textbf{return} Selected evacuation paths from solution of RMP and evacuation schedule $\Lambda$
	\end{algorithmic}
\end{algorithm*}

The algorithm begins by first solving the SP to generate an evacuation
path for each evacuation node. The RMP is then solved to schedule the
flow of evacuees over these paths. Critical evacuation nodes are then
identified and stored in $\mathcal{E}'$, and as long as this set is
not empty, the process of solving the SP to generate additional
evacuation paths and the RMP to produce an evacuation schedule that
maximizes the flow of evacuees is repeated. The algorithm terminates
when $\mathcal{E}'$ is empty or when a maximum number of iterations is
reached (maximum number of iterations is set to 10 in this
study). Upon completion, the RMP is solved one last time as an IP,
where variables $\{\varphi_p^t\}$ are set to be integers, to produce
an evacuation schedule with integral flow values. In all but the
instance with the largest population, the CPG method produces
evacuation plans very quickly. 

\section{Column Generation for Evacuation Planning}
\label{sec:CG}

This section presents the column-generation algorithm (CG) introduced
in \cite{pillac2015} to solve the NP-ZEPP, i.e., the CG method
generates non-preemptive, non-convergent zone-based evacuation
paths. {\em Interestingly, forbidding preemption makes it possible
  to design an exact column generation, avoiding the difficulties
  faced by the CPG method.} The key idea underlying the CG is to
generate {\em time-response evacuation plans} of the form $p=\langle
P,f,t_0\rangle$ where
\begin{enumerate}
\item $P$ is an evacuation path for a given zone $k$;
\item $t_0$ is the starting time of the evacuation along path $P$;
\item $f\in\mathcal{F}$ is a response curve from a set $\mathcal{F}$
  of predefined response curves.
\end{enumerate}
The CG method also features a multi-objective function that minimizes
the overall evacuation time in addition to maximizing the number of
evacuees reaching safety.

\paragraph{The Restricted Master Problem}

The RMP selects time-response
evacuation plans from a subset of feasible plans $\Omega'$ to maximize
the number of evacuees reaching safety and minimize the overall evacuation time. The formulation uses a number of constants associated with the
plans. In particular, $c_p$ denotes the cost for selecting plan $p$,
$\Omega_k\subseteq\Omega'$ is the subset of plans for evacuation node
$k$, $\omega(e)\subseteq\Omega'$ is the subset of plans that utilize
arc $e$, and $a_{p,e_t}$ denotes the flow of evacuees along arc $e$ at
time $t$ induced by plan $p$ (as prescribed by the response curve and
the departure time). The cost $c_p$ of plan $p$ is defined as a
function that applies a linear penalty on the arrival time of evacuees
at the safe node and heavily penalizes the number of evacuees that
cannot reach safety. More precisely, $c_p$ is defined as:
\begin{equation}\label{eqn:cg_mp_plancost}
c_p=\sum_{e\in p}\sum_{t\in\mathcal{H}}c_{e_t}\cdot a_{p,e_t}+\bar{c}\cdot \bar{a}_p
\end{equation}
\noindent where $\bar{a}_p$ denotes the number of evacuees not reaching safety
when executing of plan $p$, and $c_{e_t}$ and $\bar{c}$
are defined as follows:
\begin{equation}
c_{e_t}=c_{(i,j)_t}=
\begin{cases}
\frac{t}{|\mathcal{H}|}&\text{if }j\in\mathcal{S}\\
0&\text{otherwise}
\end{cases}
\end{equation}
\begin{equation}
\bar{c}=100\max_{e\in\mathcal{A},t\in\mathcal{H}}\{c_{e_t}\}\cdot\max_{k\in\mathcal{E}}\{d_k\}
\end{equation}

\begin{figure}[!t]
\begin{equation}\label{eqn:cg_mp_obj}
\min\sum_{p\in\Omega'}x_p\cdot c_p
\end{equation}
\begin{equation*}
\text{subject to}
\end{equation*}
\begin{equation}\label{eqn:cg_mp_onepath}
\sum_{p\in\Omega_k}x_p=1\qquad\forall k\in\mathcal{E}
\end{equation}
\begin{equation}\label{eqn:cg_mp_capacity}
\sum_{p\in\omega(e)}a_{p,e_t}\cdot x_p\leq u_{e_t}\qquad\forall e\in\mathcal{A},\forall t\in\mathcal{H}
\end{equation}
\begin{equation}
x_p\geq 0\qquad\forall p\in\Omega'
\end{equation}
\caption{The Restricted Master Problem for the CG Method.}
\label{fig:CG:RMP}
\end{figure}

\noindent
The RMP uses a binary variable $x_p$ to indicate whether plan
$p\in\Omega'$ is selected and is shown in
Figure \ref{fig:CG:RMP}. It is essentially a set-covering problem
with constraints on the arc capacities. Constraints
\eqref{eqn:cg_mp_onepath} ensure that only one plan is selected for
each evacuation node $k$ and Constraint \eqref{eqn:cg_mp_capacity}
enforce all arc capacities. The RMP minimizes the overall cost which
essentially causes the objective function \eqref{eqn:cg_mp_obj} to be
multi-objective and lexicographic: The RMP first maximizes of the
number of evacuees reaching safety and then minimizes the overall
evacuation time. The formulation is a linear relaxation of the
original RMP; After completion of the column-generation phase, the RMP
will be solved as a MIP. 

\paragraph{The Reduced Cost Formulation}

To find a time-response evacuation plan $p$ that can improve the current RMP, its
reduced cost $r_p$ must be negative, i.e.,
\begin{equation}\label{eqn:cg_reducedcost1}
r_p=c_p-\mathbf{a}_p^{\intercal}\boldsymbol{\pi}<0
\end{equation}
\noindent where $\mathbf{a}_p$ is the column of constraint
coefficients of $x_p$ and $\boldsymbol{\pi}$ is the vector of dual
values from the optimal solution of the RMP. Letting $\{\pi_k\}$, and
$\{\pi_{e_t}\}$ denote the dual variables of constraints
\eqref{eqn:cg_mp_onepath} and \eqref{eqn:cg_mp_capacity}, and
respectively, and substituting \eqref{eqn:cg_mp_plancost} into
\eqref{eqn:cg_reducedcost1}, the reduced cost can be formulated as:
\begin{align}\label{eqn:cg_reducedcost2}
r_p&=\sum_{e\in p}\sum_{t\in\mathcal{H}}c_{e_t}\cdot
a_{p,e_t}+\bar{c}\cdot\bar{a}_p-\pi_k-\sum_{e\in
  p}\sum_{t\in\mathcal{H}}\pi_{e_t}\cdot
a_{p,e_t}\nonumber\\ &=-\pi_k+\bar{c}\cdot\bar{a}_p+\sum_{e\in
  p}\sum_{t\in\mathcal{H}}(c_{e_t}-\pi_{e_t})\cdot a_{p,e_t}
\end{align}

\paragraph{The Pricing Subproblem}

The PSP is responsible for identifying a new time-response evacuation
plan that satisfies the negative reduced cost criteria. The
formulation of the PSP exploits some key characteristics of the
reduced cost. First, {\em the reduced cost contains terms that are
  specific to a single time-response evacuation plan $p$.} Since the
time-response evacuation plans are independent of each other, the PSP
can also be solved independently for each evacuation node
$k\in\mathcal{E}$ and for each predefined response curve
$f\in\mathcal{F}$, allowing multiple PSPs to be solved concurrently in
parallel. Moreover, since $\pi_k$ does not depend on the path, finding
a plan $p$ with negative reduced cost is equivalent to finding an
evacuation path $P$ and a start time $t_0$ that minimize the last two
terms of \eqref{eqn:cg_reducedcost2} for each $k\in\mathcal{E}$ and
for each $f\in\mathcal{F}$. Denote the last two terms of
\eqref{eqn:cg_reducedcost2} as $\text{Cost}(P,t_0)$:
\begin{equation}\label{eqn:cost_of_path0}
\text{Cost}(P,t_0)=\bar{c}\cdot\bar{a}_p+\sum_{e\in p}\sum_{t\in\mathcal{H}}(c_{e_t}-\pi_{e_t})\cdot a_{p,e_t}
\end{equation}

\begin{figure}[!t]
	\centering
	\includegraphics[width=0.6\linewidth]{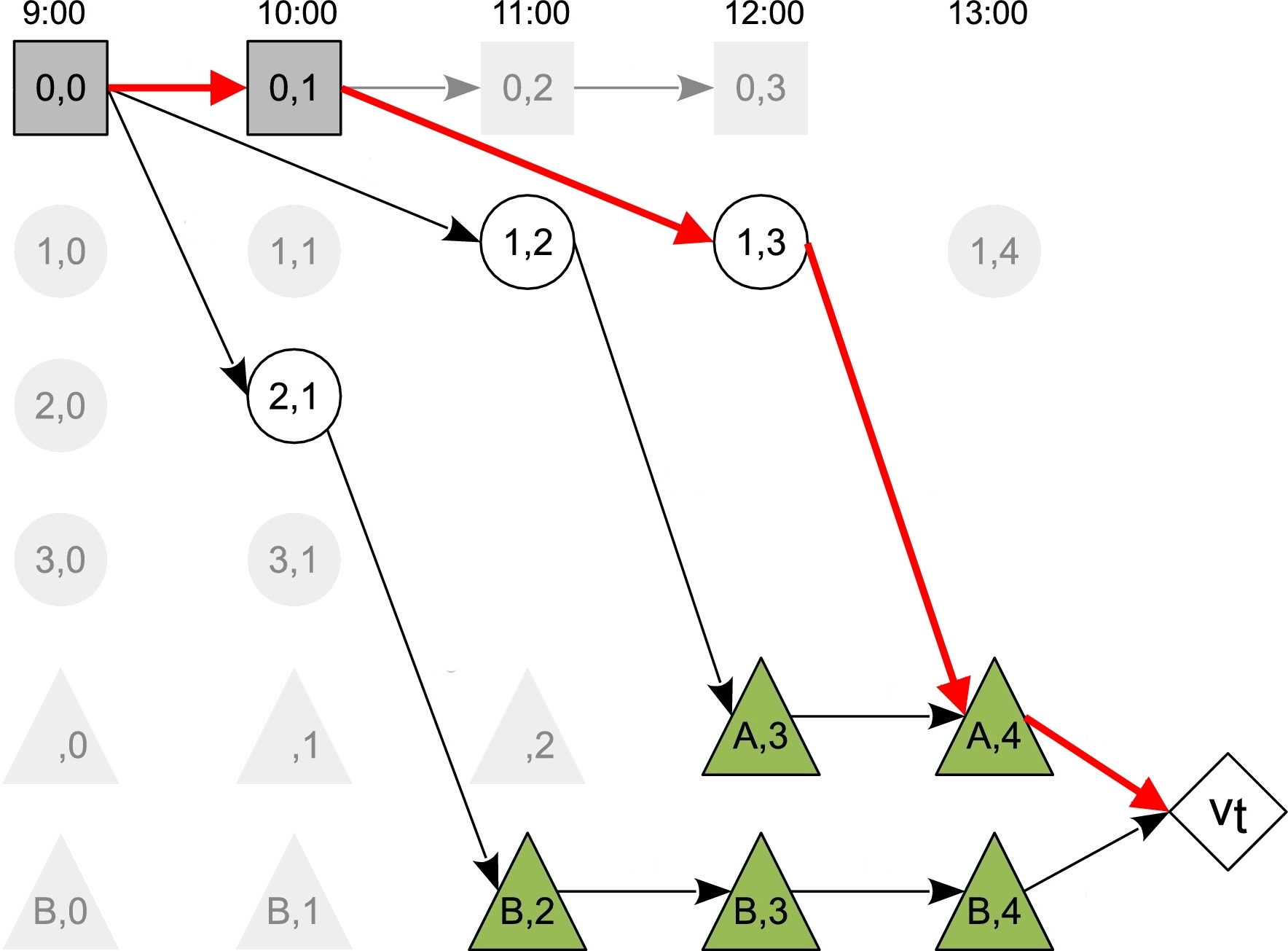}
	\caption{Path $P^x$ in the extended time-expanded graph.}
	\label{fig:path_in_time_expanded_graph}
\end{figure}

\noindent {\em The second key observation is that the path $P$ and
  evacuation start time $t_0$ that minimizes $\text{Cost}(P,t_0)$ can
  be obtained by applying a least-cost path algorithm on an extended
  time-expanded graph $\mathcal{G}^x$ with carefully defined arc
  costs. In particular, the formulation recognizes that time-response
  evacuation plan follows the same path at each time step and hence
  the arc costs can be aggregated together.} The extension to the
time-expanded graph involves the introduction of a virtual super-sink,
$v_t$, which all safe nodes $s\in\mathcal{S}^x$ are connected to with
arcs $e_t\in\mathcal{A}_s^x=\{(s,v_t)\ |\ s\in\mathcal{S}^x\}$. Now
denote by $\mathcal{A}_w^x$ the set of all infinite capacity arcs used
to model evacuees waiting at the evacuation nodes. For a given
evacuation node $k$ and response curve $f$, a path $P^x$ in
$\mathcal{G}^x$ from evacuation node $k_0$ (evacuation node $k$ and
time 0) to $v_t$ corresponds to a time-response evacuation plan
$p=\langle P,f,t_0\rangle$, where $P$ is given by the sequence of
nodes visited by $P^x$ excluding $v_t$ and $t_0$ is given by the time
of the first non-waiting arc leaving $\mathcal{E}^x$. For instance,
path $P^x$ represented by the red colored arcs in Figure
\ref{fig:path_in_time_expanded_graph} corresponds to path $P=\langle
0,1,\text{A}\rangle$ and evacuation start time $t_0$ = 10:00.

The cost of the combination (path,start time), $\text{Cost}(P,t_0)$,
can be calculated by first assigning arc costs $c_{e_t}^{\text{sp}}$
to each arc $e_t\in\mathcal{A}^x$ as follows:
\begin{equation}\label{eqn:psp_arccost1}
c_{e_t}^{\text{sp}}=\sum_{t'=t}^{|\mathcal{H}|}(c_{e_{t'}}-\pi_{e_{t'}})\cdot f(t'-t)\qquad\forall e_t\in\mathcal{A}^x\setminus(\mathcal{A}_w^x\cup\mathcal{A}_s^x)
\end{equation}
\begin{equation}\label{eqn:psp_arccost2}
c_{e_t}^{\text{sp}}=\bar{c}\cdot (d_k-F(|\mathcal{H}|-t))\qquad\forall e_t\in\mathcal{A}_s^x
\end{equation}
\begin{equation}
c_{e_t}^{\text{sp}}=0\qquad\forall e_t\in\mathcal{A}_w^x
\end{equation}
\noindent Equation \eqref{eqn:psp_arccost1} aggregates future costs of
arc $e_t$, should it be selected for a time-response evacuation plan,
while Equation \eqref{eqn:psp_arccost2} accounts for the cost of
evacuees not reaching safety for time-response evacuation plans which
end with that arc.

With these arc cost definitions, $\text{Cost}(P,t_0)$ for a path $P$
and an evacuation start time $t_0$ that corresponds to a path $P^x$
can be calculated using Equation \eqref{eqn:cost_of_path1}:
\begin{align}
\text{Cost}(P,t_0)&=\sum_{e_t\in P^x}c_{e_t}^{\text{sp}}\label{eqn:cost_of_path1}\\
&=\bar{c}\cdot (d_k-F(|\mathcal{H}|-t))+\sum_{e_t\in P^x\setminus\mathcal{A}_s^x}\sum_{t'=t}^{|\mathcal{H}|}(c_{e_{t'}}-\pi_{e_{t'}})\cdot f(t'-t)\label{eqn:cost_of_path2}\\
&=\bar{c}\cdot\bar{a}_p+\sum_{e\in p}\sum_{t\in\mathcal{H}}(c_{e_t}-\pi_{e_t})\cdot a_{p,e_t}\label{eqn:cost_of_path3}
\end{align}
\noindent Equations \eqref{eqn:cost_of_path2} and
\eqref{eqn:cost_of_path3} show that the expansion of
\eqref{eqn:cost_of_path1} will eventually lead to the original
equation for $\text{Cost}(P,t_0)$ in \eqref{eqn:cost_of_path0}.

With this formulation, the goal of the PSP, which is to find a path
$P$ and evacuation start time $t_0$ combination that minimizes
$\text{Cost}(P,t_0)$, can be accomplished by finding a shortest path
from $k_0$ to $v_t$ in the extended time-expanded graph for each
$k\in\mathcal{E}$ and $f\in\mathcal{F}$. This allows a shortest path
algorithm such as the Bellman-Ford algorithm to be applied to solve
the PSP in polynomial time.

\paragraph{The Contraflow Extension}

The CG method can also be extended to produce an evacuation plan that
allows for contraflows. It suffices to 
replace constraints \eqref{eqn:cg_mp_capacity} in the RMP with Constraints
\eqref{eqn:cg_mp_capacity1} and \eqref{eqn:cg_mp_capacity2}, and
to introduce additional Constraints \eqref{eqn:cg_mp_contraflow} and
\eqref{eqn:cg_mp_relaxy}.
\begin{equation}\label{eqn:cg_mp_capacity1}
\sum_{p\in\omega(e)}a_{p,e_t}\cdot x_p\leq u_{e_t}\qquad\forall e\in\mathcal{A}\setminus\mathcal{A}_c,\forall t\in\mathcal{H}
\end{equation}
\begin{equation}\label{eqn:cg_mp_capacity2}
\sum_{p\in\omega(e)}a_{p,e_t}\cdot x_p\leq y_e\cdot u_{e_t}+(1-y_{\bar{e}})\cdot u_{\bar{e}_t}\qquad\forall e\in\mathcal{A}_c,\forall t\in\mathcal{H}
\end{equation}
\begin{equation}\label{eqn:cg_mp_contraflow}
y_e+y_{\bar{e}}\geq 1\qquad\forall e\in\mathcal{A}_c
\end{equation}
\begin{equation}\label{eqn:cg_mp_relaxy}
y_e\in[0,1]\qquad\forall e\in\mathcal{A}_c
\end{equation}

\noindent
Constraints \eqref{eqn:cg_mp_capacity1} enforce capacity on arcs that
may not be used in contraflow, while Constraints
\eqref{eqn:cg_mp_capacity2} do the same for arcs that may. Constraints
\eqref{eqn:cg_mp_contraflow} prevents arcs $e$ and $\bar{e}$ from
being used in contraflow simultaneously for road segment $(e,\bar{e})$
with $e\in\hat{\mathcal{A}}_c,\bar{e}\in\check{\mathcal{A}}_c$, and
Constraints \eqref{eqn:cg_mp_relaxy} apply a linear relaxation on
variable $y_e$. Once the column generation procedure has terminated,
variable $y_e$ is made binary, the final RMP is solved as a MIP, and
the rest of the CG method remains unchanged.

\paragraph{Elementary Paths}
\label{sec:cycle_free_path_extension}

\begin{figure}[!t]
\centering
\includegraphics[width=0.5\linewidth]{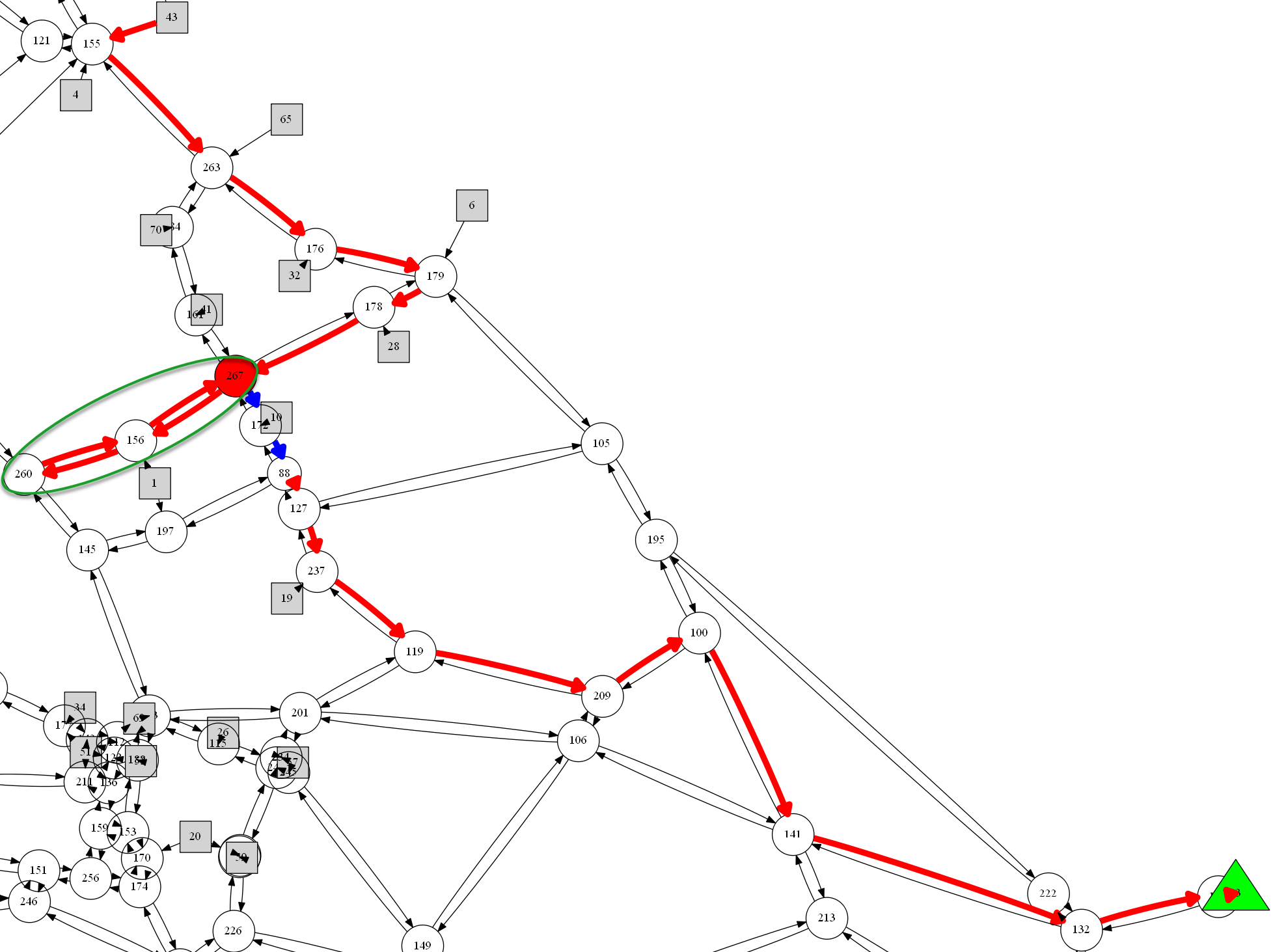}
\caption{An Example of Non-Elementary Path.}
\label{fig:NEP}
\end{figure}

The time-expanded graph $\mathcal{G}^x$ is by construction acyclic as
its arcs only connect nodes at different time steps. As such, the
shortest paths identified in the PSP are also acyclic. However, this
fact does not preclude the PSP from generating paths that visit the
same transit node in $\mathcal{G}^x$ at different time steps, as there
are no restrictions enforced in the shortest path algorithms
preventing such paths from being generated. While such paths are
acyclic in $\mathcal{G}^x$, their corresponding counterparts in the
static graph $\mathcal{G}$ contain cycles, as they visit the same
transit node more than once.

These cyclic paths are called non-elementary (they visit the same node multiple
times) and an example of such a non-elementary path is shown in Figure
\ref{fig:NEP}. Non-elementary paths in the static graph are
undesirable in real evacuations, as they give evacuees the impression
that the evacuation plans are sub-optimal and reduce trust in
emergency services. However, when the CG algorithm is applied to the
real case study, about 44\% of the generated evacuation paths are not
elementary.

\begin{figure}[!t]
\begin{equation}\label{eqn:feas_sp_obj}
\min\sum_{e_t\in\mathcal{A}^x}c_{e_t}^{\text{sp}}\cdot x_{e_t}
\end{equation}
\begin{equation*}
\text{subject to}
\end{equation*}
\begin{equation}\label{eqn:feas_sp_source}
\sum_{e_t\in\delta^{+}(k_0)}x_{e_t}=1
\end{equation}
\begin{equation}\label{eqn:feas_sp_flowconserve}
\sum_{e_t\in\delta^{-}(i)}x_{e_t}-\sum_{e_t\in\delta^{+}(i)}x_{e_t}=0\qquad\forall i\in\mathcal{N}^x\setminus\{k_0,v_t\}
\end{equation}
\begin{equation}\label{eqn:feas_sp_sink}
\sum_{e_t\in\delta^{-}(v_t)}x_{e_t}=1
\end{equation}
\begin{equation}\label{eqn:feas_sp_resource}
\sum_{i_t\in \mathrm{\Lambda}(i)}\sum_{e_t\in\delta^{+}(i_t)} x_{e_t}\leq 1\qquad\forall i\in\mathcal{T}
\end{equation}
\begin{equation}\label{eqn:feas_sp_vars}
x_{e_t}\in\{0,1\}\qquad\forall e_t\in\mathcal{A}^x
\end{equation}
\caption{The Pricing Subproblem With Elementary Paths.}
\label{fig:PSPEP}
\end{figure}

This section outlines the pricing subproblem proposed in
\cite{Hafiz2017} that only generates time-response evacuation plans
with elementary paths. Let $\mathrm{\Lambda}(i)$ denote the set of
time-expanded nodes in $\mathcal{G}^x$ for a node $i\in\mathcal{T}$,
i.e., $\mathrm{\Lambda}(i)=\{i_t\,|\,t\in\mathcal{H}\}$. A path $P^x$
in $\mathcal{G}^x$ corresponds to an elementary path $P$ in
$\mathcal{G}$ if and only if $P^x$ visits at most a single node in
$\mathrm{\Lambda}(i)$ for each node $i \in\mathcal{T}$. As a result,
instead of finding a least-cost path, the revised PSP must find a
least-cost path that is also an elementary path in the static graph.
Figure \ref{fig:PSPEP} depicts the new formulation of the pricing
problem.  The formulation uses binary decision variable $x_{e_t}$ to
indicate whether edge $e_t$ should be selected as part of the shortest
path. Objective function \eqref{eqn:feas_sp_obj} minimizes the total
cost of the path. Constraint \eqref{eqn:feas_sp_source} specifies that
exactly one path should emanate from source node $k_0$, while
constraint \eqref{eqn:feas_sp_sink} ensures the path ends at super
sink node $v_t$. Constraints \eqref{eqn:feas_sp_flowconserve} enforce
path continuity at every node other than the source and super
sink. Finally, Constraints \eqref{eqn:feas_sp_resource} guarantee that
each transit node is visited by the path at most once throughout the
entire time horizon.

This version of the PSP is a shortest path problem with resource
constraints \cite{irnich2005} which is known to be NP-hard
\cite{gary1979}. In this particular formulation, the resources are
simply the unit ``visited'' resources associated with each transit
node in $\mathcal{G}^x$. The set of all time steps of a particular
transit node, $\{i_t\,|\,i\in\mathcal{T},\forall t\in\mathcal{H}\}$,
is allocated only one unit of this ``visited'' resource, and the
resource is completely consumed if this node were to be visited by a
path. For the case study in this paper, this constrained shortest-path
problem must be solved repeatedly for a graph with approximately 30000
nodes and 75000 arcs.

While solving formulation
\eqref{eqn:feas_sp_obj}-\eqref{eqn:feas_sp_vars} using a MIP solver
will result in the desired shortest elementary path, the hybrid strategy
proposed in \cite{Hafiz2017} is capable to obtain these paths
faster. The hybrid strategy combines the above formulation with a
$k$-shortest-path-based algorithm based on the implementation of
Jimenez and Marzal's Recursive Enumeration Algorithm (REA)
\cite{jimenez1999}.  This algorithm incrementally generates a
$k$\textsuperscript{th}-shortest path based on information of the
$(k-1)$ shortest paths. It can be used to find the shortest elementary
path by first generating the shortest path (setting $k=1$). If the
path is elementary, the algorithm terminates.  Otherwise, the next
shortest path is generated by the REA (by incrementing $k$ by 1) and
the elementary check is applied on this path. This process is repeated
until an elementary path is obtained.

Computational experiments on the case study show that the
$k$-shortest-path-based algorithm is extremely fast at finding
shortest elementary paths when the required $k$ values are relatively
small ($k<10^5$). Unfortunately, in some rare instances, the value of
$k$ required to obtain an elementary path is extremely large (in the
millions), and under these circumstances, solving the MIP formulation
produces faster results. Therefore, the hybrid strategy combines both
methods by first utilizing the $k$-shortest-path-based algorithm to
find shortest elementary paths up to a threshold value for $k$ (in this
study, the threshold is set to $10^5$). If this $k$-threshold is
reached and a elementary path is yet to be found, a MIP formulation is
solved.  This hybrid strategy exploits the strengths of both methods
and is extremely effective at identifying shortest elementary paths
quickly in almost all cases.

\section{Clearance Time Minimization}
\label{sec:mincleartime}

Of the four methods presented, only the CG method has a
multi-objective function which minimizes total evacuation time in
addition to maximizing number of evacuees reaching safety. The BN, BC,
and CPG methods only optimize for the latter goal. However, evacuation
authorities are also deeply interested in the minimum clearance time,
\ie the smallest amount of time to evacuate an entire region. A
precise definition of minimum clearance time, $h^*$, is as follows:
\begin{equation}\label{eqn:mincleartime}
h^*=\min\Bigg\{t\in\mathcal{H}\ \bigg|\ z(\text{EPP}(\mathcal{G},[0..t]))=\sum_{k\in\mathcal{E}}d_k\Bigg\}
\end{equation}
\noindent where $\text{EPP}(G,\mathcal{H})$ denotes the optimal
solution obtained from an EPP formulation given static graph
$\mathcal{G}$ and time horizon $\mathcal{H}$ as inputs. This section
shows how to obtain the minimum clearance time for each method.

\paragraph{Benders Non-convergent and Convergent Methods}
\label{sec:mct_bn_bc}

The BN and BC methods each consists of an RMP and an SP which generate
upper and lower bounds to the objective value. As proposed in
\cite{romanski2016}, a lower bound on the minimum clearance time,
$h^\dagger$, can first be obtained by performing a binary search over
the time horizon using just the RMP. Next, a sequential search using
the full BN or BC method can be used to find $h^*$, beginning from its
lower bound $h^\dagger$. This approach seems preferable over a binary
search for the second stage as $h^*$ is very likely closer to the
lower bound $h^\dagger$ and hence the second part of the algorithm
will converge faster by starting a sequential search from that
time. Algorithm \ref{alg:mincleartime} summarizes the entire approach.

\begin{algorithm*}[!tbhp]
	\caption{Clearance Time Minimization for BN and BC Methods}\label{alg:mincleartime}
	\begin{algorithmic}[1]
		\State $h^\dagger\leftarrow\min\{t\in\mathcal{H}\,|\,z(\text{RMP}(\mathcal{G},[0..t]))=\sum_{k\in\mathcal{E}}d_k\}$
		\State $h^*\leftarrow\min\{t\in[h^\dagger..\mathcal{H}]\,|\,z(\text{EPP}(\mathcal{G},[0..t]))=\sum_{k\in\mathcal{E}}d_k\}$
		\State \textbf{return} $h^*$, evacuation paths from solution of RMP, and evacuation schedule from solution of SP
	\end{algorithmic}
\end{algorithm*}

\paragraph{Conflict-based Path Generation Method}\label{sec:mct_cpg}

Since the CPG method does not maintain upper and lower bounds to the
objective value of the EPP, the clearance time can be performed by
a binary search over the time horizon the full CPG method.

\paragraph{Column Generation Method}\label{sec:mct_cg}

Even though the CG method's multi-objective function minimizes total
evacuation time, the quantity is not equivalent to clearance
time. Clearance time is equivalent to the time at which the last
evacuee arrives at its safe node, and this is not the quantity being
minimized in the objective function. Minimizing total evacuation time
might result in minimal clearance time, but it might also produce
suboptimal clearance times as the penalty incurred by late arrival of
the last evacuee could possibly be diluted by early arrival costs.
Application of a binary search over the time horizon is an option,
however this paper did not resort to this approach due to the
significant run times of the CG method. Therefore, the minimum
clearance time experiments only report the arrival time of the last
evacuee produced by the CG method while fully recognizing that it
might be suboptimal as the method's objective function does not
explicitly minimize clearance times.

\section{The Case Study}
\label{sec:casestudy}

This section presents a case study to investigate the effectiveness
and run times of the four methods on a real-world evacuation
scenario. It also discusses some preliminary observations on various
properties of the evacuation algorithms.

\begin{figure}[!t]
\centering
\includegraphics[width=1.0\linewidth]{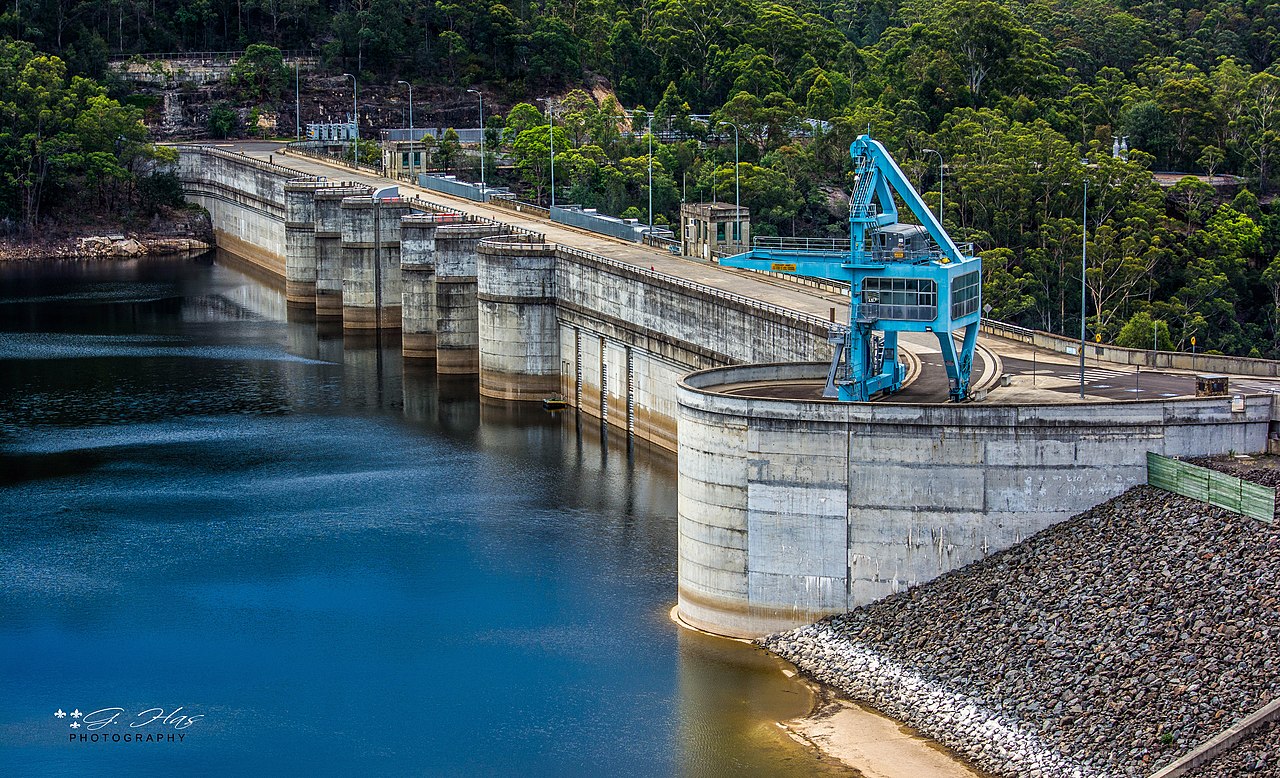}
\caption{The Warragamba Dam in New South Wales.}
\label{fig:warragamba}
\end{figure}

\begin{figure}[!tbhp]
	\centering
	\begin{subfigure}{.55\textwidth}
		\centering
		\includegraphics[width=\linewidth]{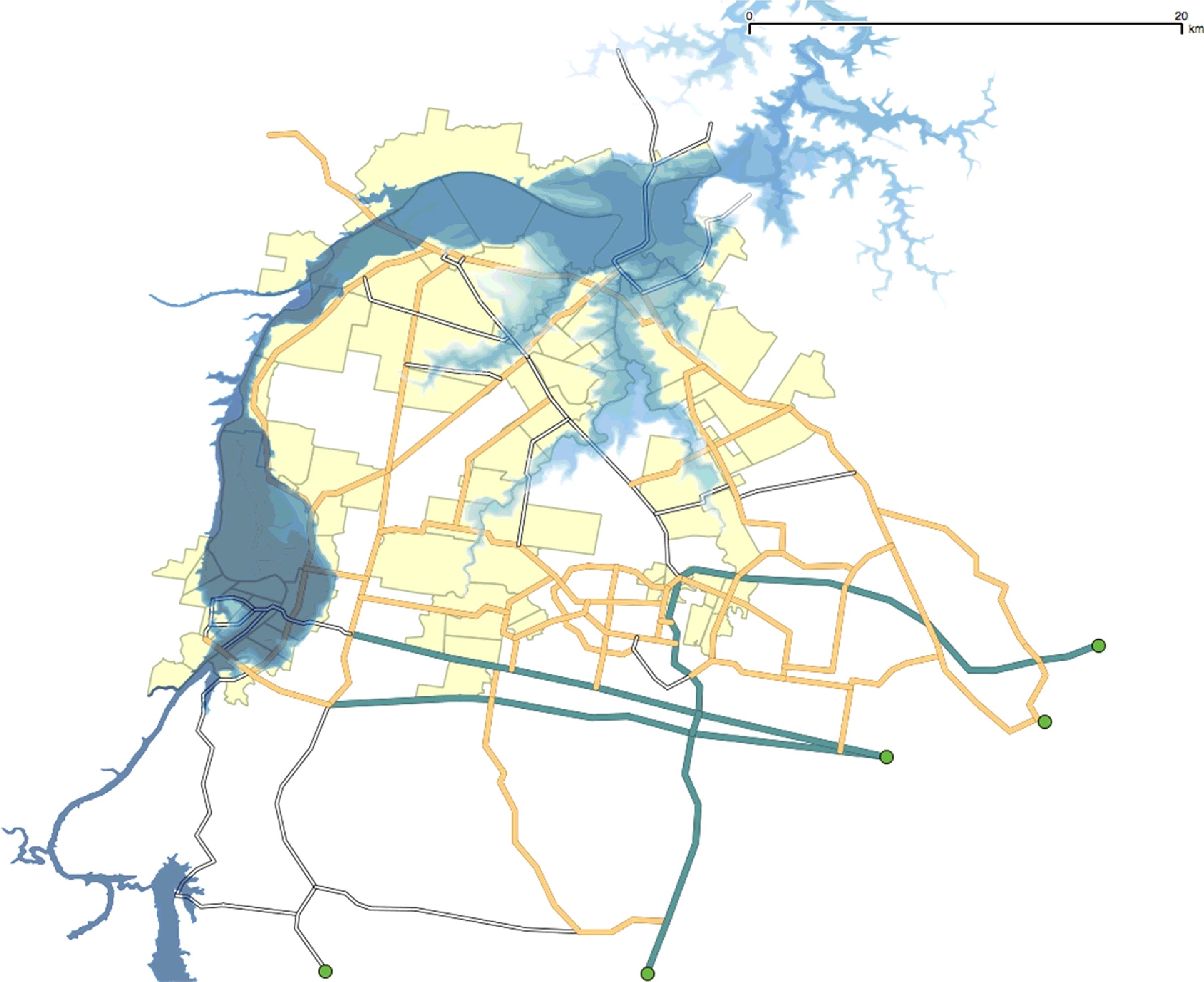}
		\caption{HN map (from \cite{pillac2016})}
		\label{fig:hn_map}
	\end{subfigure}%
	\begin{subfigure}{.45\textwidth}
		\centering
		\includegraphics[width=\linewidth]{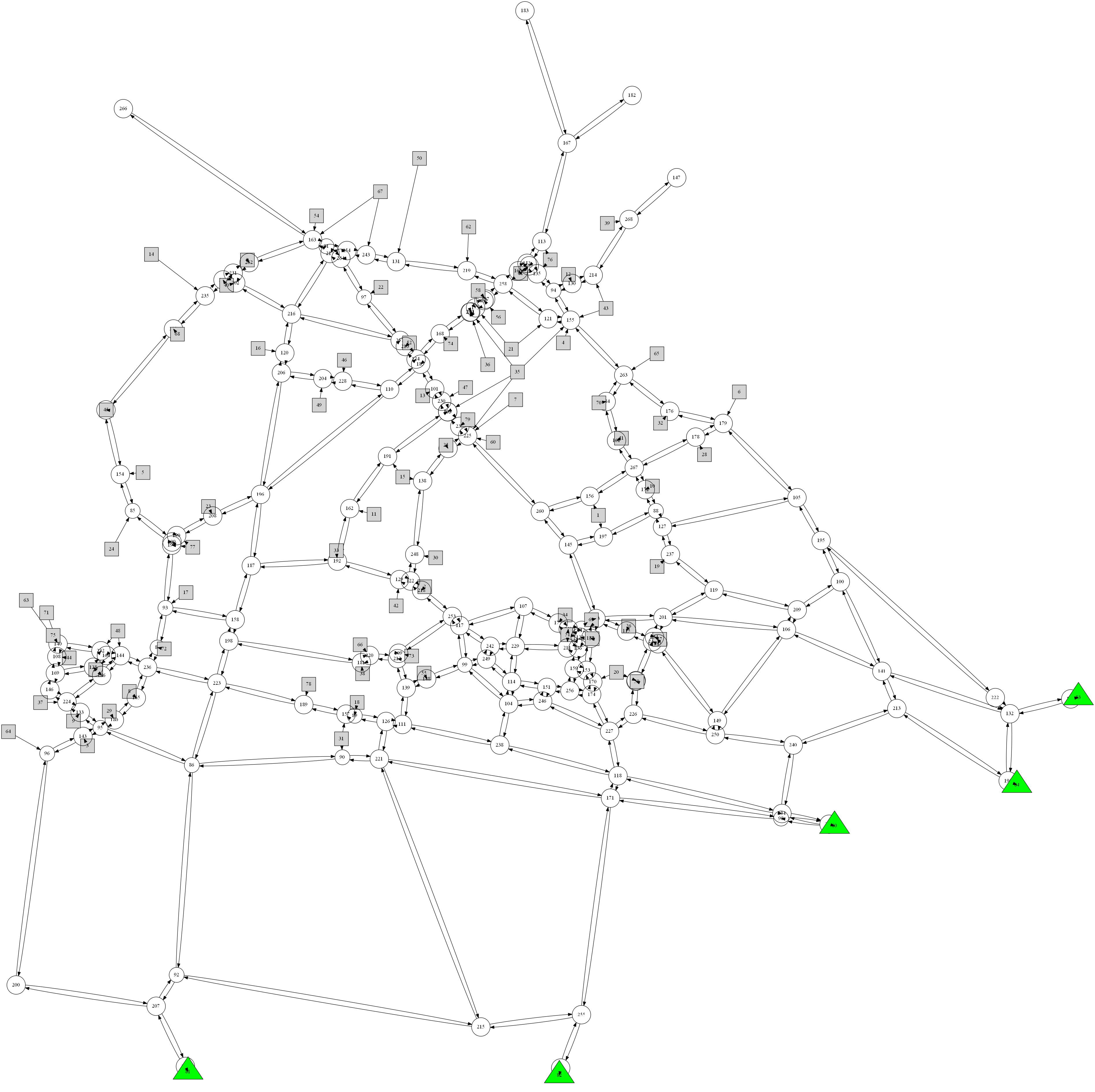}
		\caption{HN evacuation graph}
		\label{fig:hn_graph}
	\end{subfigure}
	\caption{The Map and Evacuation Graph of the Hawkesbury-Nepean Region.}
	\label{fig:hn_map_n_graph}
\end{figure}

\paragraph{The Case Study}
The case study is the Hawkesbury-Nepean (HN) region located north-west
of Sydney, Australia, which is separated from the Blue Mountains, a
catchment area, by the Warragamba dam (See Figure
\ref{fig:warragamba}. This dam often spills over and, if it breaks, it
would create a massive flooding event for West Sydney. The region's
evacuation graph consists of 80 evacuation nodes, 184 transit nodes,
and 5 safe nodes. Evacuation node deadlines and road block times for
the region were obtained from a hydro-dynamic simulation of a 1-in-100
years flood event. The region has a total of 38,343 vehicles to be
evacuated in its base instance; However the results also consider the
effect of increasing the population size by linearly scaling the base
instance by a factor $x \in[1.0,3.0]$, as West Sidney has sustained
significant population growth. Figure \ref{fig:hn_map} shows a
bird's-eye view of the entire region, while Figure \ref{fig:hn_graph}
shows its corresponding evacuation graph, with squares, circles, and
triangles representing evacuation, transit, and safe nodes
respectively. Each instance (or evacuation scenario) is referred to as
HN80-Ix from this point forth, where $x$ is the population scaling
factor.

\paragraph{Experimental Settings}
The performance of each method is evaluated under two settings: (a) a
deadline setting where the maximum number of evacuees reaching safety
is determined within a fixed time horizon $\mathcal{H}$ = 10 hours,
and (b) a minimum clearance time setting where the smallest amount of
time needed to evacuate the entire region is determined. Under both
settings, the time horizon is discretized into 5 minute time
steps. For the CG method, the set of predefined response curves
$\mathcal{F}$ was populated with step response curves with flow rates
$\gamma\in\{2,6,10,25,50\}$ evacuees per time step. All methods were
implemented in C++, with multi-threaded components being handled using
OpenMP, used in conjunction with \textsc{Gurobi} 6.5.2 to solve all
mathematical programs. All experiments were conducted on the
University of Michigan's Flux high-performance computing cluster, with
each utilizing 8 cores of a 2.5 GHz Intel Xeon E5-2680v3 processor and
64 GB of RAM.

\begin{figure}[!tbhp]
	\centering
	\begin{subfigure}{1\textwidth}
		\centering
		\includegraphics[width=0.6\linewidth]{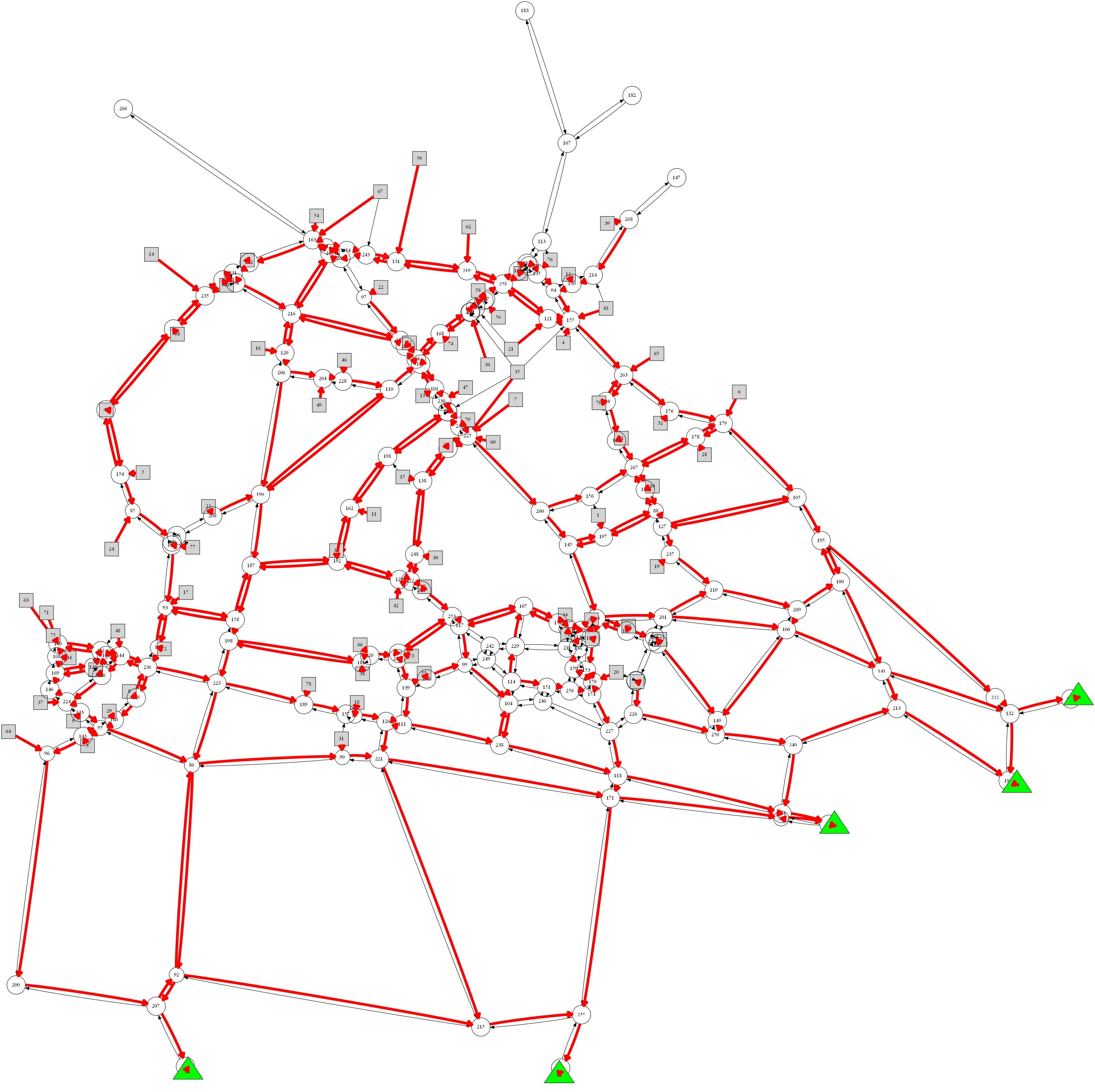}
		\caption{Non-Convergent Evacuation Paths.}
		\label{fig:nonconv_path}
	\end{subfigure}
\newline
	\begin{subfigure}{1\textwidth}
		\centering
		\includegraphics[width=0.6\linewidth]{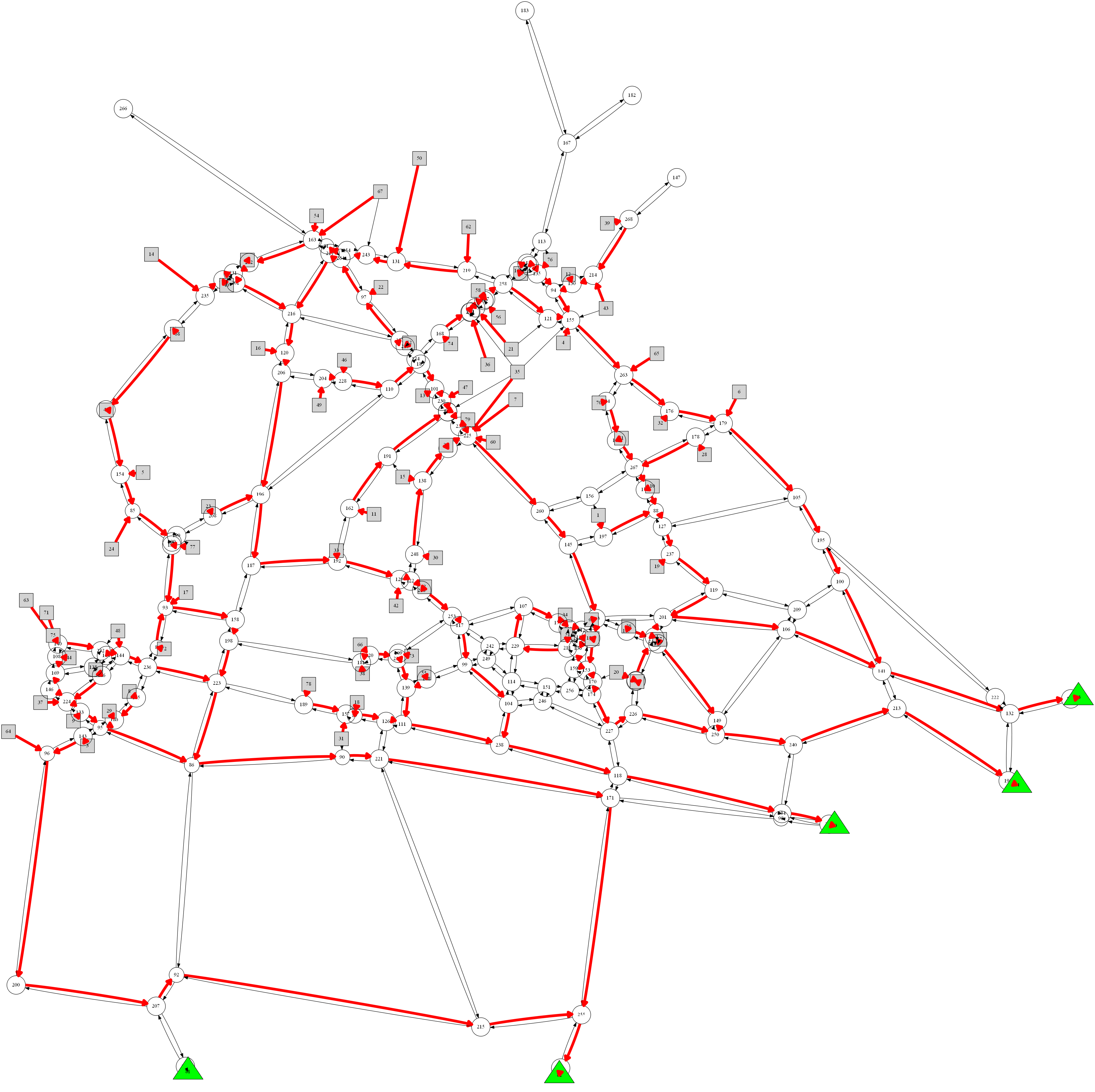}
		\caption{Convergent Evacuation Paths.}
		\label{fig:conv_path}
	\end{subfigure}
	\caption{Non-Convergent and Convergent Evacuation Paths
          without Contraflow Generated by the BN and BC Methods.}
	\label{fig:nonconv_n_conv_paths}
\end{figure}

\paragraph{Convergent versus Non-Convergent Paths}
Figures \ref{fig:nonconv_path} and \ref{fig:conv_path} illustrate
examples of generated non-convergent and convergent evacuation paths
that do not use contraflow. The paths, whose arcs are represented by
red colored arrows, are overlaid on top of the evacuation graph in
these figures. It can be seen in Figure \ref{fig:conv_path} that
convergent paths form a tree with leaves at the evacuation nodes
(squares in the graph) and rooted at the safe nodes (triangles in the
graph). The non-convergent paths in Figure \ref{fig:nonconv_path} do
not share this property; however not being constrained by this
property allows more arcs to be utilized for evacuation (at the
expense introducing potential delays caused by driver hesitation when
multiple paths fork at road intersections).

\begin{figure}[!t]
	\centering
        \includegraphics[width=0.7\linewidth]{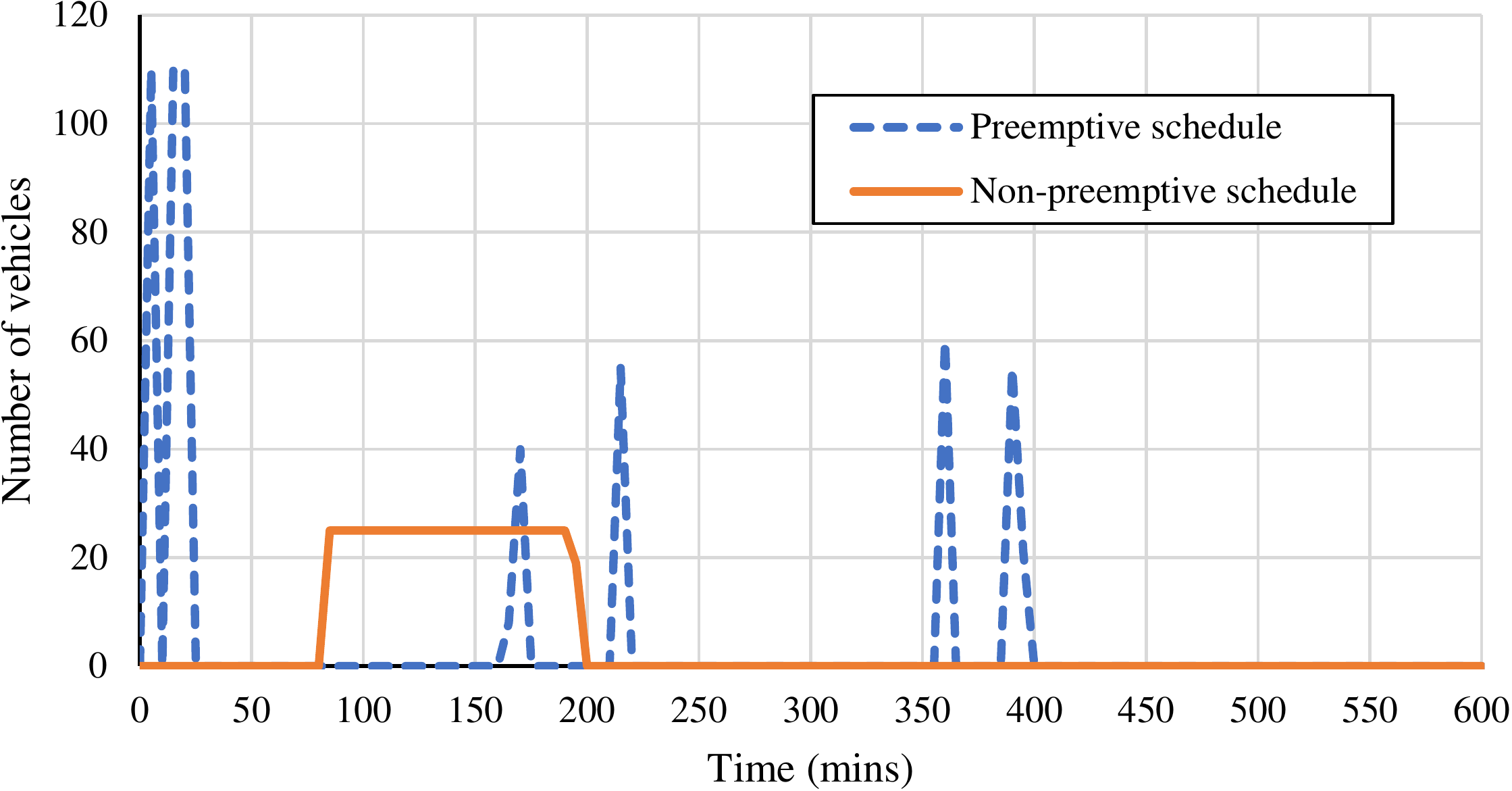}
	\caption{Preemptive and Non-Preemptive Evacuation Schedules for an Evacuation Node with 569 Evacuees (Generated by the BN and CG Methods Respectively).}
	\label{fig:evac_schedule}
\end{figure}

\paragraph{Preemptive Versus Non-Preemptive Schedules}
Figure \ref{fig:evac_schedule} highlights the difference between
preemptive and non-preemptive evacuation schedules generated for an
evacuation node with 569 evacuees. The preemptive schedule is
characterized by multiple spikes followed by interruptions in evacuee
departure rates over time, which may lead to some challenges in
the enforcement of the schedule. This is contrasted with the
non-preemptive schedule which uses a step response curve with a flow
rate of 25 evacuees every 5 minutes. Evacuation is started at the
85\textsuperscript{th} minute, and a constant departure rate is
maintained until the node has been completely evacuated, making the
schedule arguably easier to enforce compared to the preemptive one.

\paragraph{Convergence of the Benders Decomposition}

\begin{figure}[!tbhp]
	\centering
	\begin{subfigure}{.5\textwidth}
		\centering
		\includegraphics[width=.9\linewidth]{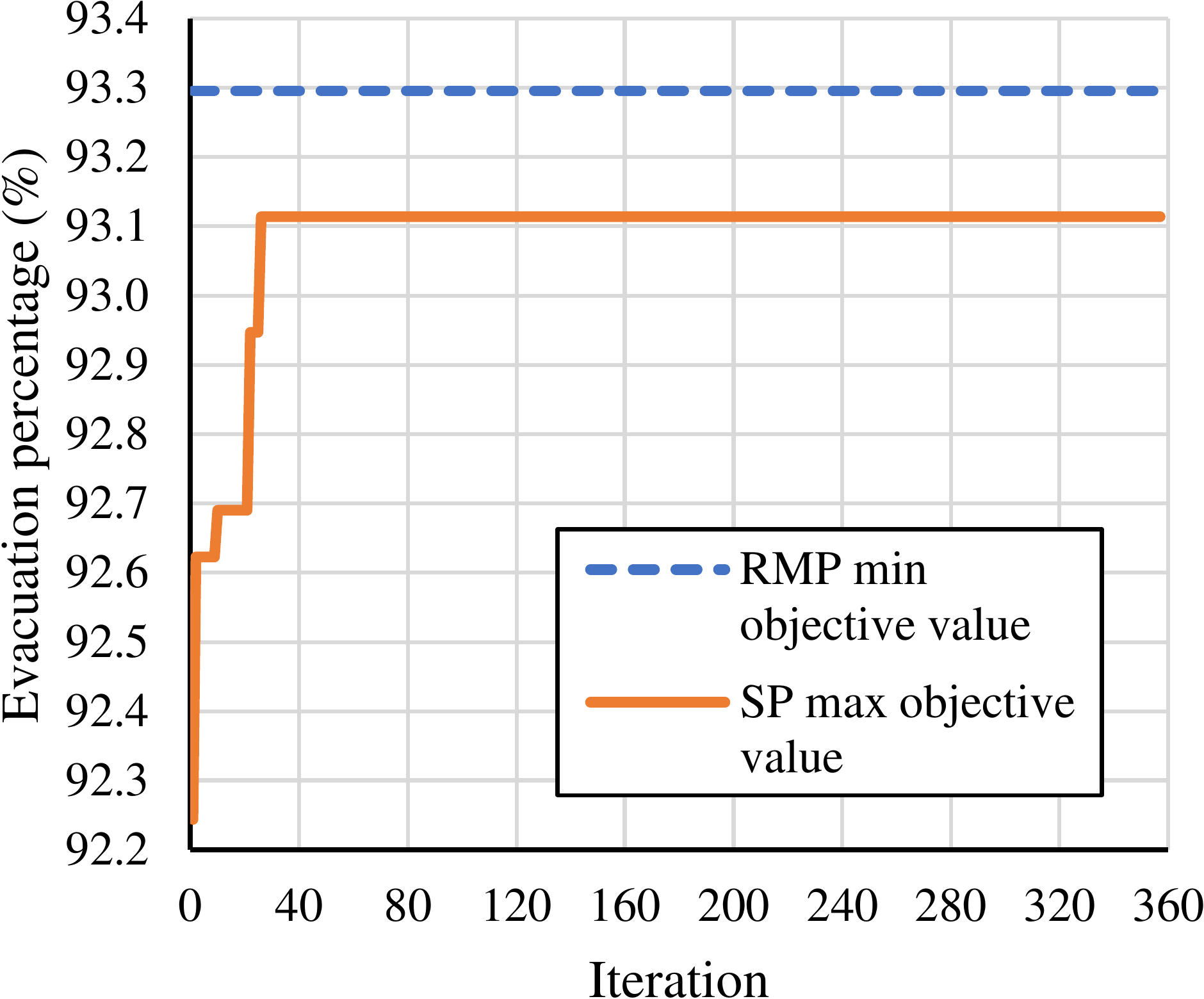}
		\caption{Convergence of the BN Method}
		\label{fig:bn_convergence}
	\end{subfigure}%
	\begin{subfigure}{.5\textwidth}
		\centering
		\includegraphics[width=.9\linewidth]{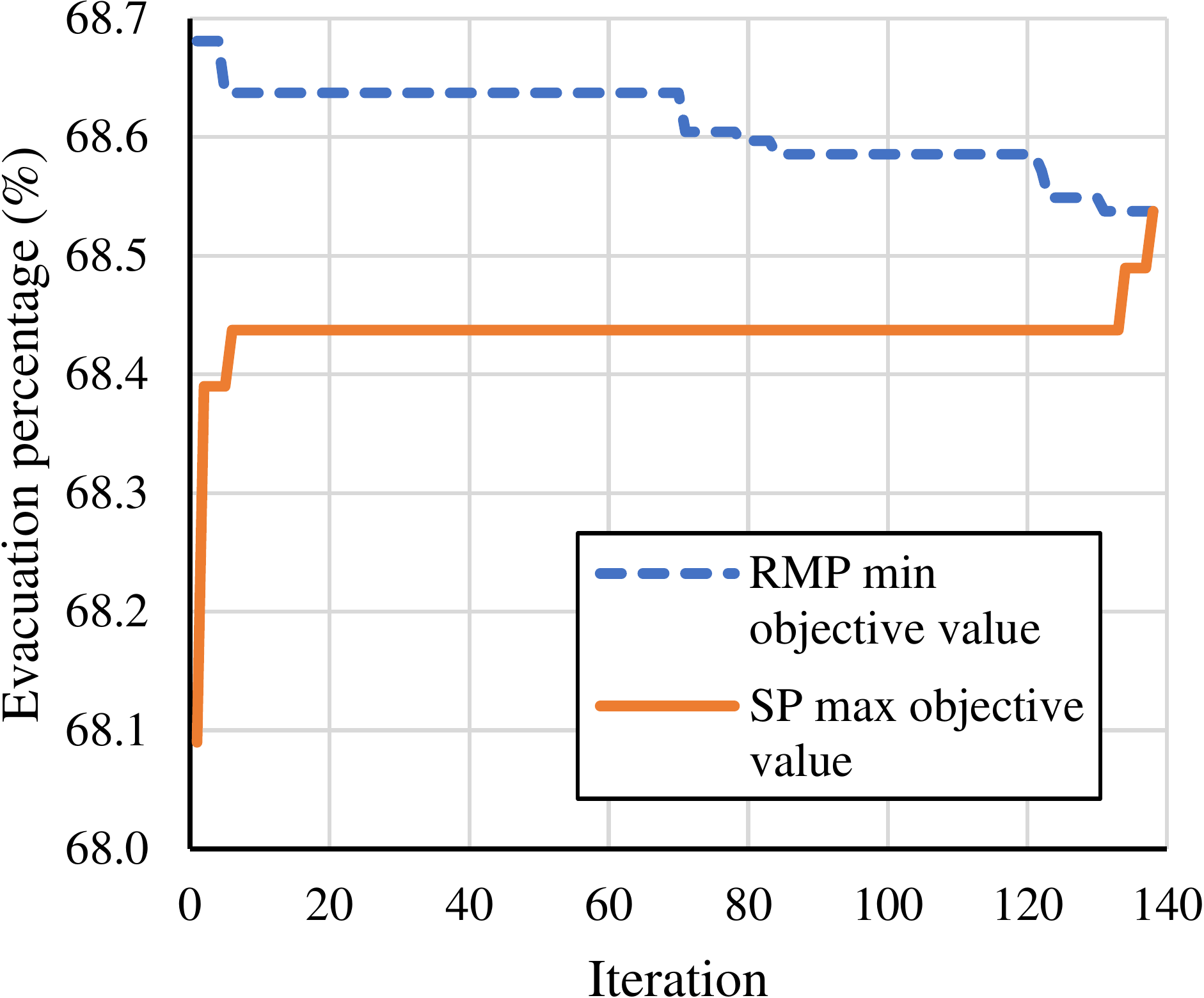}
		\caption{Convergence of BC Method}
		\label{fig:bc_convergence}
	\end{subfigure}
	\caption{Convergence characteristics of the BN and BC methods under the deadline setting for the HN80-I3.0 instance.}
	\label{fig:bn_n_bc_convergence}
\end{figure}

Figures \ref{fig:bn_convergence} and \ref{fig:bc_convergence} reveal
how the upper and lower bounds of the BN and BC methods converge over
time for a particular experiment in the deadline setting. For this
instance, the BC method converged in less than 140
iterations, while the BN method did not, even after 360 iterations, at
which point the algorithm was terminated as it exceeded a set time
limit (time limits for each method are elaborated further in Section
\ref{sec:macroscopic}). Nevertheless, it can be seen that
final optimality gap is very small ($\approx 0.2\%$) and this gap was
attained in less than 40 iterations.

\paragraph{Elementary and Non-Elementary Paths}
Table \ref{tab:compare_ori_new_psp} compares results of the column
generation phase of the CG method without and with the elementary path
restriction. The key takeaway is that the two formulations produce the
same optimal values. Minor differences in the last two instances may
be due to the column generation phase being terminated before
convergence as a CPU time limit of 5760 minutes was reached (CPU time
limits applied to all experiments are detailed further in Section
\ref{sec:deadline_setting}). Restricting attention to elementary paths
increases the CPU times, which is not surprising, since finding
shortest paths subject to resource constraints is an NP-hard problem:
Even though the hybrid strategy employed for finding elementary paths
is highly effective, it still cannot compete with the polynomial time
Bellman-Ford algorithm used in the original formulation. Nevertheless,
it is interesting to observe that the CPU time advantage of the
original formulation diminishes as the population size increases.
Another intriguing observation is that CG with elementary paths
reaches optimality in fewer iterations and it has fewer columns in its
final RMP in almost all instances.

\begin{table}[!tbhp]
	\begin{center}
		\begin{tabular}{|c|r|r|r|r|r|r|r|r|}
			\hline \multirow{2}{*}{\textbf{Instance}} &
                        \multicolumn{4}{c|}{\textbf{Original Shortest
                            Path PSP}} &
                        \multicolumn{4}{C{60mm}|}{\textbf{New Elementary
                            Shortest Path PSP}} \\ \cline{2-9} &
                        \multicolumn{1}{C{7mm}|}{\textbf{Iter \#}} &
                        \multicolumn{1}{C{14mm}|}{\textbf{Column \#}}
                        & \multicolumn{1}{C{11mm}|}{\textbf{CPU Time
                            (mins)}} &
                        \multicolumn{1}{C{15mm}|}{\textbf{Optimal Obj
                            Val}} &
                        \multicolumn{1}{C{7mm}|}{\textbf{Iter \#}} &
                        \multicolumn{1}{C{14mm}|}{\textbf{Column \#}}
                        & \multicolumn{1}{C{11mm}|}{\textbf{CPU Time
                            (mins)}} &
                        \multicolumn{1}{C{15mm}|}{\textbf{Optimal Obj
                            Val}} \\ \hline HN80-I1.0 & 79 & 12251 &
                        39 & 8816 & 79 & 11678 & 124 & 8816 \\ \hline
                        HN80-I1.1 & 104 & 15072 & 136 & 10405 & 95 &
                        13853 & 218 & 10405 \\ \hline HN80-I1.2 & 229
                        & 22571 & 799 & 12116 & 190 & 19543 & 834 &
                        12116 \\ \hline HN80-I1.4 & 152 & 20184 & 690
                        & 15935 & 108 & 17048 & 404 & 15935 \\ \hline
                        HN80-I1.7 & 178 & 21871 & 1312 & 22635 & 120 &
                        19883 & 2760 & 22635 \\ \hline HN80-I2.0 & 197
                        & 31418 & 5760 & 30490 & 145 & 25051 & 5760 &
                        30490 \\ \hline HN80-I2.5 & 121 & 22806 & 5760
                        & 46189 & 129 & 23513 & 5760 & 46188 \\ \hline
                        HN80-I3.0 & 132 & 31726 & 5760 &
                        \multicolumn{1}{R{15mm}|}{$1.960\times10^9$} &
                        87 & 21233 & 5760 &
                        \multicolumn{1}{R{15mm}|}{$1.961\times10^9$}
                        \\ \hline
		\end{tabular}
	\end{center}
	\caption{Results of column generation phase of CG method using original shortest path and new elementary shortest path PSP formulations (when no contraflow is allowed).}
	\label{tab:compare_ori_new_psp}
\end{table}

\begin{figure}[!tbhp]
	\centering
	\includegraphics[width=0.7\linewidth]{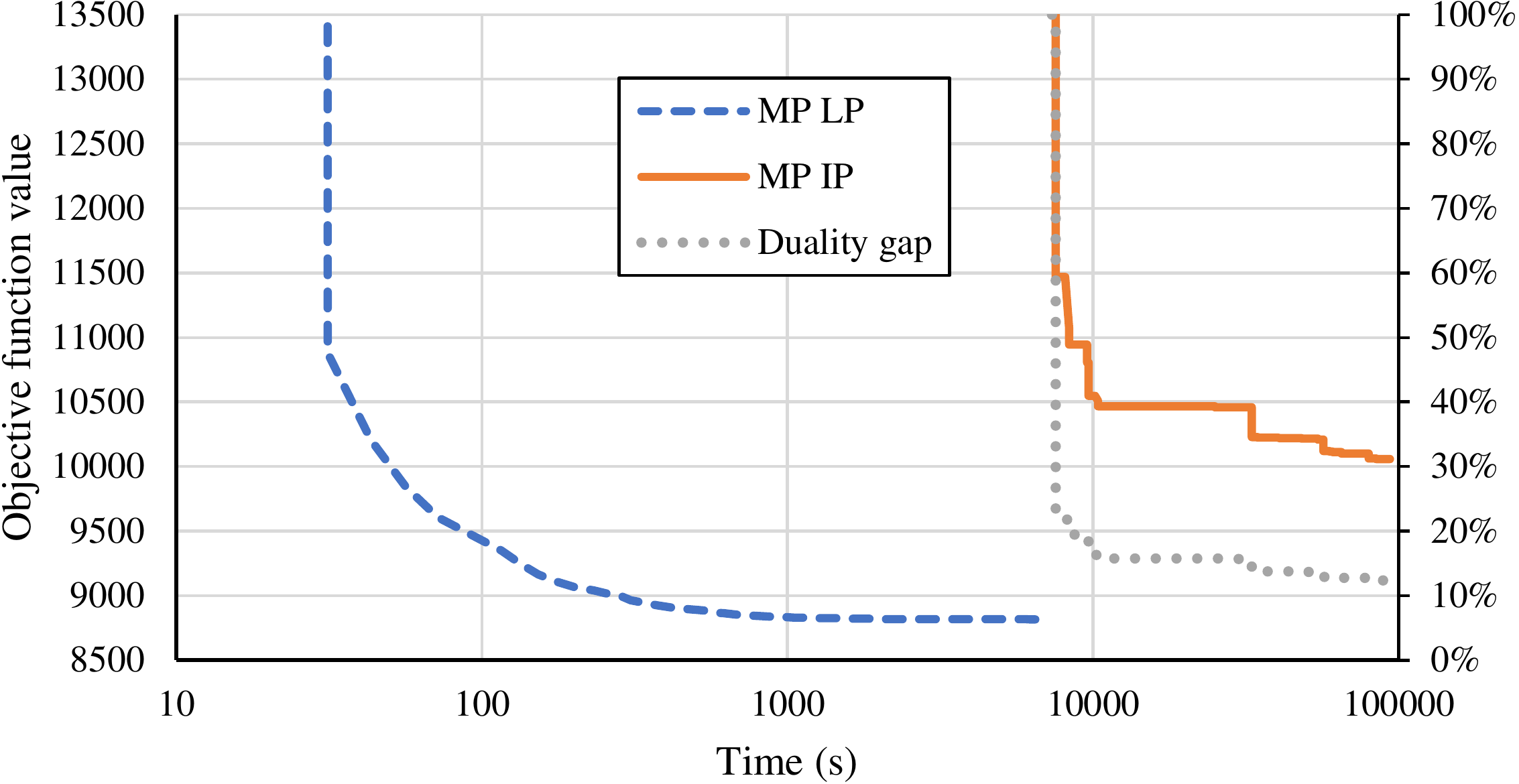}
	\caption{Evolution of the solution quality of the CG method over time for the HN80-I1.0 instance.}
	\label{fig:cg_convergence}
\end{figure}

\paragraph{Convergence of the Column Generation} 
Figure \ref{fig:cg_convergence} demonstrates how the objective value
of the RMP of the CG method evolves over time during its
column-generation phase. It also shows the evolution of the objective
value of the best incumbent solution found for the RMP when it is
solved as an IP in its last iteration, together with the progression
of its optimality gap over time. It can be seen that there is a steep
decline in the objective value of the restricted master problem within
the first 100 seconds, after which the value slowly approaches a
minimum. The same trend is observed when the restricted master problem
is solved as a MIP, with the optimality gap of the best incumbent
solution settling to a value of $\approx 13.5\%$ when the algorithm
reached its time limit.

\section{Macroscopic Evaluation}
\label{sec:macroscopic}

This section provides a summary of the results obtained from all four methods under the deadline and minimum clearance time settings, together with the specific conditions under which each method is applied.

\subsection{The Deadline Setting}
\label{sec:deadline_setting}

Under the deadline setting, each method maximizes the number of
evacuees reaching safety for the HN80-Ix instances (with x
$\in[1.0,3.0]$) within a fixed time horizon $\mathcal{H}$ = 10
hours.

\paragraph{The BN Method}

The results of the BN method are summarized in Tables
\ref{tab:BN_deadline_noCF} and \ref{tab:BN_deadline_withCF} without
and with contraflow respectively. As shown in Algorithm
\ref{alg:bendersnonconv}, BN first searches for the tightest time
horizon $t^*$ that can preserve the optimal solution of the RMP,
$z(\text{RMP}(\mathcal{G},\mathcal{H}))$. The tables show the CPU time
taken for this first phase in column ``$t^*$ CPU Time''. Each RMP
instance in this procedure is given a CPU time limit of 10
minutes. The tables also show the total number of iterations required
to complete the entire method as well as the corresponding total CPU
time taken. The complete method is allocated a CPU time limit of 24
hours. Finally, the tables show the minimum objective value of the RMP
$z_{\text{RMP}}$ and the maximum objective value of the SP,
$z_{\text{SP,max}}$ at termination in terms of evacuation percentage,
as well as the optimality gap, i.e., the percentage difference
between the upper and lower bounds of the solution given by
$\frac{z_{\text{RMP}}-z_{\text{SP,max}}}{z_{\text{SP,max}}}$.

There are three key observations from Tables
\ref{tab:BN_deadline_noCF} and \ref{tab:BN_deadline_withCF}: (a)
Without contraflows, the method converges to an evacuation percentage
of 100\% for all instances except for HN80-I3.0, for which the method
produces a 93.1\% evacuation percentage. This was also the only
instance where the method did not converge within the allocated CPU
time limit; however the final optimality gap of 0.20\% assures that
the obtained solution is very close to being optimal. (b) The method
converges after only 1 iteration for all instances except HN80-I2.5
and HN80-I3.0 without contraflow. For instances that took 1 iteration,
the bulk of the CPU time is spent on the search for $t^*$. (c) When
using contraflows, the method converges faster on all instances: The
increased network capacity provided by contraflows makes the
evacuation instances easier to solve.

\begin{table}[!t]
	\begin{center}
		\begin{tabular}{ |C{16mm}|R{12mm}|R{24mm}|R{30mm}|R{15mm}|R{12mm}|R{22mm}| }
			\hline
			\multicolumn{1}{|C{16mm}|}{\textbf{Instance}} & \multicolumn{1}{C{12mm}|}{\textbf{Iter \#}} & \multicolumn{1}{C{24mm}|}{$\boldsymbol{t^*}$ \textbf{CPU Time (mins)}} & \multicolumn{1}{C{30mm}|}{\textbf{Total CPU Time (mins)}} & \multicolumn{1}{C{15mm}|}{$\boldsymbol{z_{\text{RMP}}}$ \textbf{(\%)}} & \multicolumn{1}{C{12mm}|}{$\boldsymbol{z_{\text{SP,max}}}$ \textbf{(\%)}} & \multicolumn{1}{C{22mm}|}{\textbf{Optimality Gap (\%)}}\\
			\hline
			HN80-I1.0 & 1 & 131 & 135 & 100.0 & 100.0 & 0.00 \\
			\hline
			HN80-I1.1 & 1 & 111 & 117 & 100.0 & 100.0 & 0.00 \\
			\hline
			HN80-I1.2 & 1 & 104 & 110 & 100.0 & 100.0 & 0.00 \\
			\hline
			HN80-I1.4 & 1 & 92 & 96 & 100.0 & 100.0 & 0.00 \\
			\hline
			HN80-I1.7 & 1 & 103 & 110 & 100.0 & 100.0 & 0.00 \\
			\hline
			HN80-I2.0 & 1 & 79 & 84 & 100.0 & 100.0 & 0.00 \\
			\hline
			HN80-I2.5 & 4 & 38 & 98 & 100.0 & 100.0 & 0.00 \\
			\hline
			HN80-I3.0 & 358 & 5 & 1449 & 93.3 & 93.1 & 0.20 \\
			\hline
		\end{tabular}
	\end{center}
	\caption{Results of the BN method under the deadline setting without contraflow.}
	\label{tab:BN_deadline_noCF}
\end{table}

\begin{table}[!t]
	\begin{center}
		\begin{tabular}{ |C{16mm}|R{12mm}|R{24mm}|R{30mm}|R{15mm}|R{12mm}|R{22mm}| }
			\hline
			\multicolumn{1}{|C{16mm}|}{\textbf{Instance}} & \multicolumn{1}{C{12mm}|}{\textbf{Iter \#}} & \multicolumn{1}{C{24mm}|}{$\boldsymbol{t^*}$ \textbf{CPU Time (mins)}} & \multicolumn{1}{C{30mm}|}{\textbf{Total CPU Time (mins)}} & \multicolumn{1}{C{15mm}|}{$\boldsymbol{z_{\text{RMP}}}$ \textbf{(\%)}} & \multicolumn{1}{C{12mm}|}{$\boldsymbol{z_{\text{SP,max}}}$ \textbf{(\%)}} & \multicolumn{1}{C{22mm}|}{\textbf{Optimality Gap (\%)}}\\
			\hline
			HN80-I1.0 & 1 & 43 & 47 & 100.0 & 100.0 & 0.00 \\
			\hline
			HN80-I1.1 & 1 & 43 & 45 & 100.0 & 100.0 & 0.00 \\
			\hline
			HN80-I1.2 & 1 & 41 & 43 & 100.0 & 100.0 & 0.00 \\
			\hline
			HN80-I1.4 & 1 & 48 & 50 & 100.0 & 100.0 & 0.00 \\
			\hline
			HN80-I1.7 & 1 & 42 & 45 & 100.0 & 100.0 & 0.00 \\
			\hline
			HN80-I2.0 & 1 & 42 & 45 & 100.0 & 100.0 & 0.00 \\
			\hline
			HN80-I2.5 & 1 & 36 & 39 & 100.0 & 100.0 & 0.00 \\
			\hline
			HN80-I3.0 & 1 & 35 & 37 & 100.0 & 100.0 & 0.00 \\
			\hline
		\end{tabular}
	\end{center}
	\caption{Results of the BN method under the deadline setting with contraflow.}
	\label{tab:BN_deadline_withCF}
\end{table}

\paragraph{The BC Method} Tables \ref{tab:BC_deadline_noCF} and
\ref{tab:BC_deadline_withCF} show results of the BC method under the
deadline setting without and with contraflow respectively. The tables
are similar to Tables \ref{tab:BN_deadline_noCF} and
\ref{tab:BN_deadline_withCF} because of the similarities between the
BN and BC methods. The CPU time limits allocated for the BC method are
different however: Each RMP instance in the search for $t^*$ is given
a limit of 10 minutes, whereas the entire method is allocated only 2
hours of CPU time due its faster convergence.

The key observations from these two tables are as follows: (a) The
method evacuates everyone for all instances when contraflow is
allowed, and for instances HN80-Ix with x $\in[1.0,1.7]$ when
contraflow is not allowed. However, unlike the BN method, this method
converges to optimal solutions for all instances as evidenced by their
0\% optimality gaps. (b) For instances in which a 100\% evacuation
rate is achieved, the method converges in just 1 iteration, and the
bulk of the CPU time in these instances is spent searching for
$t^*$. (c) Except for instance HN80-I1.4, the method converges faster
when contraflow is allowed. As mentioned earlier, the BC method is
extremely effective in finding high-quality evacuations quickly.

\begin{table}[!tbhp]
	\begin{center}
		\begin{tabular}{ |C{16mm}|R{12mm}|R{24mm}|R{30mm}|R{15mm}|R{12mm}|R{22mm}| }
			\hline
			\multicolumn{1}{|C{16mm}|}{\textbf{Instance}} & \multicolumn{1}{C{12mm}|}{\textbf{Iter \#}} & \multicolumn{1}{C{24mm}|}{$\boldsymbol{t^*}$ \textbf{CPU Time (mins)}} & \multicolumn{1}{C{30mm}|}{\textbf{Total CPU Time (mins)}} & \multicolumn{1}{C{15mm}|}{$\boldsymbol{z_{\text{RMP}}}$ \textbf{(\%)}} & \multicolumn{1}{C{12mm}|}{$\boldsymbol{z_{\text{SP,max}}}$ \textbf{(\%)}} & \multicolumn{1}{C{22mm}|}{\textbf{Optimality Gap (\%)}}\\
			\hline
			HN80-I1.0 & 1   & 0.35  & 0.51  & 100.0 & 100.0 & 0.00 \\
			\hline
			HN80-I1.1 & 1   & 10.19 & 10.28 & 100.0 & 100.0 & 0.00 \\
			\hline
			HN80-I1.2 & 1   & 10.15 & 10.24 & 100.0 & 100.0 & 0.00 \\
			\hline
			HN80-I1.4 & 1   & 1.13  & 1.21  & 100.0 & 100.0 & 0.00 \\
			\hline
			HN80-I1.7 & 1   & 10.21 & 10.32 & 100.0 & 100.0 & 0.00 \\
			\hline
			HN80-I2.0 & 14  & 0.09  & 1.80  & 96.9  & 96.9  & 0.00 \\
			\hline
			HN80-I2.5 & 22  & 0.28  & 3.43  & 81.4  & 81.4  & 0.00 \\
			\hline
			HN80-I3.0 & 138 & 0.01  & 16.22 & 68.5  & 68.5  & 0.00 \\
			\hline
		\end{tabular}
	\end{center}
	\caption{Results of the BC method under the deadline setting without contraflow.}
	\label{tab:BC_deadline_noCF}
\end{table}

\begin{table}[!tbhp]
	\begin{center}
		\begin{tabular}{ |C{16mm}|R{12mm}|R{24mm}|R{30mm}|R{15mm}|R{12mm}|R{22mm}| }
			\hline
			\multicolumn{1}{|C{16mm}|}{\textbf{Instance}} & \multicolumn{1}{C{12mm}|}{\textbf{Iter \#}} & \multicolumn{1}{C{24mm}|}{$\boldsymbol{t^*}$ \textbf{CPU Time (mins)}} & \multicolumn{1}{C{30mm}|}{\textbf{Total CPU Time (mins)}} & \multicolumn{1}{C{15mm}|}{$\boldsymbol{z_{\text{RMP}}}$ \textbf{(\%)}} & \multicolumn{1}{C{12mm}|}{$\boldsymbol{z_{\text{SP,max}}}$ \textbf{(\%)}} & \multicolumn{1}{C{22mm}|}{\textbf{Optimality Gap (\%)}}\\
			\hline
			HN80-I1.0 & 1 & 0.19 & 0.27 & 100.0 & 100.0 & 0.00 \\
			\hline
			HN80-I1.1 & 1 & 0.22 & 0.30 & 100.0 & 100.0 & 0.00 \\
			\hline
			HN80-I1.2 & 1 & 0.28 & 0.37 & 100.0 & 100.0 & 0.00 \\
			\hline
			HN80-I1.4 & 1 & 2.31 & 2.42 & 100.0 & 100.0 & 0.00 \\
			\hline
			HN80-I1.7 & 1 & 0.50 & 0.60 & 100.0 & 100.0 & 0.00 \\
			\hline
			HN80-I2.0 & 1 & 0.61 & 0.70 & 100.0 & 100.0 & 0.00 \\
			\hline
			HN80-I2.5 & 1 & 0.38 & 0.48 & 100.0 & 100.0 & 0.00 \\
			\hline
			HN80-I3.0 & 1 & 0.20 & 0.30 & 100.0 & 100.0 & 0.00 \\
			\hline
		\end{tabular}
	\end{center}
	\caption{Results of the BC method under the deadline setting with contraflow.}
	\label{tab:BC_deadline_withCF}
\end{table}

\paragraph{The CPG Method}
Results for the CPG method are outlined in Tables
\ref{tab:CPG_deadline_noCF} and \ref{tab:CPG_deadline_withCF}
respectively. The tables show the number of iterations, the CPU time,
and the evacuation percentage. As mentioned in Section
\ref{sec:cpg_algorithm}, the CPG algorithm is allowed a maximum of 10
iterations, after which it was terminated even if it still had
critical nodes which are not completely evacuated. Additionally, the
RMP and SP are each allocated a CPU time limit of 1 hour.

The tables show that a 100\% evacuation rate is achieved in all
instances but HN80-I2.5 and HN80-I3.0 without contraflow. These two
instances were terminated by the iteration limit. The results also
show relatively short CPU times of less than a minute for instances
HN80-Ix with x $\in[1.0,2.0]$, and a spike in CPU time for instance
HN80-I3.0 both with and without contraflow. This is possibly caused by
the increase in difficulty of solving these instances. The positive
correlation between the number of iterations required and the value of
x for the instances provides further evidence that the heuristic
requires more iterations to evacuate more evacuees. Similar to the BN
and BC methods, the CPU times for each instance are smaller when
contraflow is allowed. On almost all instances, the CPG method
is very effective as well. 

\begin{table}[!tbhp]
	\begin{center}
		\begin{tabular}{ |c|r|r|r|}
			\hline
			\multicolumn{1}{|c|}{\textbf{Instance}} & \multicolumn{1}{c|}{\textbf{Iter \#}} & \multicolumn{1}{c|}{\textbf{CPU Time (mins)}} & \multicolumn{1}{c|}{\textbf{Evacuation Percentage (\%)}} \\
			\hline
			HN80-I1.0 & 2  & 0.05  & 100.0 \\
			\hline
			HN80-I1.1 & 3  & 0.07  & 100.0 \\
			\hline
			HN80-I1.2 & 3  & 0.08  & 100.0 \\
			\hline
			HN80-I1.4 & 3  & 0.10  & 100.0 \\
			\hline
			HN80-I1.7 & 3  & 0.15  & 100.0 \\
			\hline
			HN80-I2.0 & 4  & 0.69  & 100.0 \\
			\hline
			HN80-I2.5 & 10 & 11.04 & 96.7  \\
			\hline
			HN80-I3.0 & 10 & 96.48 & 86.3  \\
			\hline
		\end{tabular}
	\end{center}
	\caption{Results of deadline setting CPG method without contraflow.}
	\label{tab:CPG_deadline_noCF}
\end{table}

\begin{table}[!tbhp]
	\begin{center}
		\begin{tabular}{ |c|r|r|r|}
			\hline
			\multicolumn{1}{|c|}{\textbf{Instance}} & \multicolumn{1}{c|}{\textbf{Iter \#}} & \multicolumn{1}{c|}{\textbf{CPU Time (mins)}} & \multicolumn{1}{c|}{\textbf{Evacuation Percentage (\%)}} \\
			\hline
			HN80-I1.0 & 1 & 0.04  & 100.0 \\
			\hline
			HN80-I1.1 & 1 & 0.03  & 100.0 \\
			\hline

                        HN80-I1.2 & 2 & 0.04  & 100.0 \\
			\hline
			HN80-I1.4 & 2 & 0.05  & 100.0 \\
			\hline
			HN80-I1.7 & 2 & 0.06  & 100.0 \\
			\hline
			HN80-I2.0 & 3 & 0.15  & 100.0 \\
			\hline
			HN80-I2.5 & 4 & 2.44  & 100.0 \\
			\hline
			HN80-I3.0 & 4 & 21.23 & 100.0 \\
			\hline
		\end{tabular}
	\end{center}
	\caption{Results of the CPG method under the deadline setting with contraflow.}
	\label{tab:CPG_deadline_withCF}
\end{table}

\paragraph{The CG Method} 
The results of the CG method without and with contraflow are
summarized in Tables \ref{tab:CG_deadline_noCF} and
\ref{tab:CG_deadline_withCF} respectively. They show the number of
iterations and the CPU time taken by the column generation phase to
converge, the total number of columns generated by the phase, the
total CPU time taken by the entire CG method, the optimality gap of
the final integer solution of the RMP, calculated using formula
$\frac{z_\text{RMP,MIP}-z_\text{RMP,LP}}{z_\text{RMP,MIP}}$ where
$z_\text{RMP,LP}$ an d $z_\text{RMP,MIP}$ are the final objective
values of the column generation and the MIP respectively, and the
evacuation percentage achieved by the method. A couple of CPU time
limits are applied to the method: 96 hours for the column generation
phase and 24 hours for the final MIP.

The tables show that the final MIP reaches its 24 hour time limit in
all instances. The column-generation phase also reaches its 96 hour
time limit in the most challenging instances, HN80-Ix with x
$\in[2.0,3.0]$ without contraflow and HN80-I3.0 with contraflow. A
100\% evacuation rate is achieved in instances HN80-Ix with x
$\in[1.0,1.4]$ without contraflow and with x $\in[1.0,2.0]$ with
contraflow. Correspondingly, the optimality gaps of these instances
are relatively low, being $<20\%$ when contraflow is not allowed and
$<10\%$ when it is allowed. This indicates that the quality of the
solutions obtained for these instances are relatively high.  The large
optimality gap of the other instances can be explained by the large
penalty incurred in their objective values because not everyone is
evacuated safely in the final integer solutions. A steady increase in
CPU times is also observed as the population scaling factor x
increases. For the same instance, the CPU times are smaller when
contraflow is allowed, and so are the optimality gaps when 100\%
evacuation is achieved.

\begin{table}[!tbhp]
	\begin{center}
		\begin{tabular}{ |C{16mm}|R{12mm}|R{19mm}|R{30mm}|R{20mm}|R{14mm}|R{20mm}| }
			\hline
			\multicolumn{1}{|C{16mm}|}{\textbf{Instance}} & \multicolumn{1}{C{12mm}|}{\textbf{Iter \#}} & \multicolumn{1}{C{19mm}|}{\textbf{Column \#}} & \multicolumn{1}{C{30mm}|}{\textbf{Column Generation CPU Time (mins)}} & \multicolumn{1}{C{20mm}|}{\textbf{Total CPU Time (mins)}} & \multicolumn{1}{C{14mm}|}{\textbf{Optimality Gap (\%)}} & \multicolumn{1}{C{20mm}|}{\textbf{Evacuation Percentage (\%)}}\\
			\hline
			HN80-I1.0 & 75  & 11626 & 153  & 1593 & 13.5  & 100.0 \\
			\hline
			HN80-I1.1 & 88  & 13463 & 261  & 1701 & 14.9  & 100.0 \\
			\hline
			HN80-I1.2 & 191 & 18020 & 792  & 2232 & 15.0  & 100.0 \\
			\hline
			HN80-I1.4 & 116 & 17412 & 681  & 2121 & 16.9  & 100.0 \\
			\hline
			HN80-I1.7 & 119 & 19415 & 1275 & 2715 & 100.0 & 96.1 \\
			\hline
			HN80-I2.0 & 163 & 25883 & 5760 & 7200 & 100.0 & 92.5 \\
			\hline
			HN80-I2.5 & 126 & 23185 & 5760 & 7200 & 100.0 & 89.7 \\
			\hline
			HN80-I3.0 & 95  & 23146 & 5760 & 7200 & 63.0  & 81.1 \\
			\hline
		\end{tabular}
	\end{center}
	\caption{Results of the CG method under the deadline setting without contraflow.}
	\label{tab:CG_deadline_noCF}
\end{table}

\begin{table}[!tbhp]
	\begin{center}
		\begin{tabular}{ |C{16mm}|R{12mm}|R{19mm}|R{30mm}|R{20mm}|R{14mm}|R{20mm}| }
			\hline
			\multicolumn{1}{|C{16mm}|}{\textbf{Instance}} & \multicolumn{1}{C{12mm}|}{\textbf{Iter \#}} & \multicolumn{1}{C{19mm}|}{\textbf{Column \#}} & \multicolumn{1}{C{30mm}|}{\textbf{Column Generation CPU Time (mins)}} & \multicolumn{1}{C{20mm}|}{\textbf{Total CPU Time (mins)}} & \multicolumn{1}{C{14mm}|}{\textbf{Optimality Gap (\%)}} & \multicolumn{1}{C{20mm}|}{\textbf{Evacuation Percentage (\%)}}\\
			\hline
			HN80-I1.0 & 60  & 3757  & 25   & 1465 & 4.2  & 100.0 \\
			\hline
			HN80-I1.1 & 69  & 4636  & 55   & 1495 & 4.5  & 100.0 \\
			\hline
			HN80-I1.2 & 105 & 5398  & 137  & 1577 & 5.1  & 100.0 \\
			\hline
			HN80-I1.4 & 152 & 8015  & 304  & 1744 & 5.0  & 100.0 \\
			\hline
			HN80-I1.7 & 226 & 11880 & 829  & 2269 & 7.1  & 100.0 \\
			\hline
			HN80-I2.0 & 310 & 14971 & 1899 & 3339 & 9.9  & 100.0 \\
			\hline
			HN80-I2.5 & 453 & 21474 & 5741 & 7181 & 99.9 & 99.7 \\
			\hline
			HN80-I3.0 & 201 & 17094 & 5760 & 7200 & 99.9 & 99.8 \\
			\hline
		\end{tabular}
	\end{center}
	\caption{Results of the CG method under the deadline setting with contraflow.}
	\label{tab:CG_deadline_withCF}
\end{table}

\paragraph{Comparison of the Evacuation Rates}

\begin{figure}[!tbhp]
	\centering
	\includegraphics[width=0.98\linewidth]{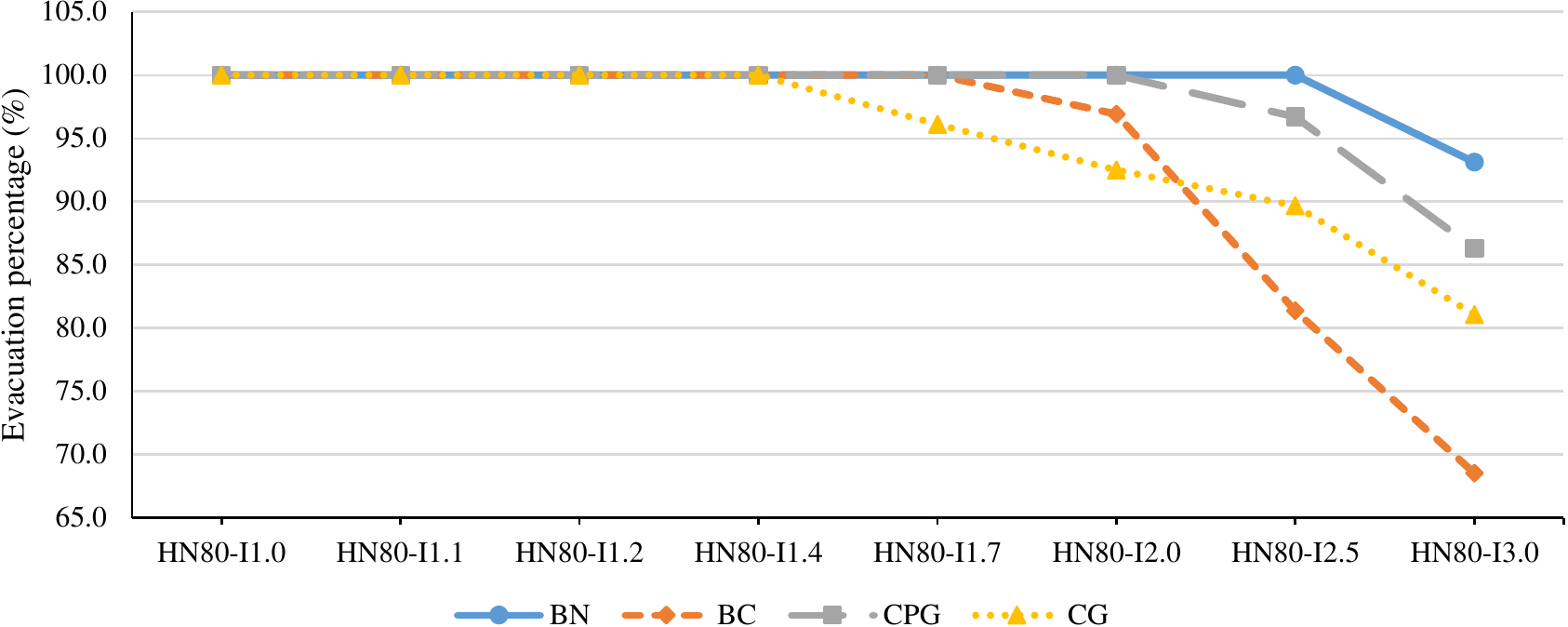}
	\caption{Evacuation percentage under deadline setting without contraflow summary.}
	\label{fig:deadline_summary_nocf}
\end{figure}

\begin{figure}[!tbhp]
	\centering
	\includegraphics[width=0.98\linewidth]{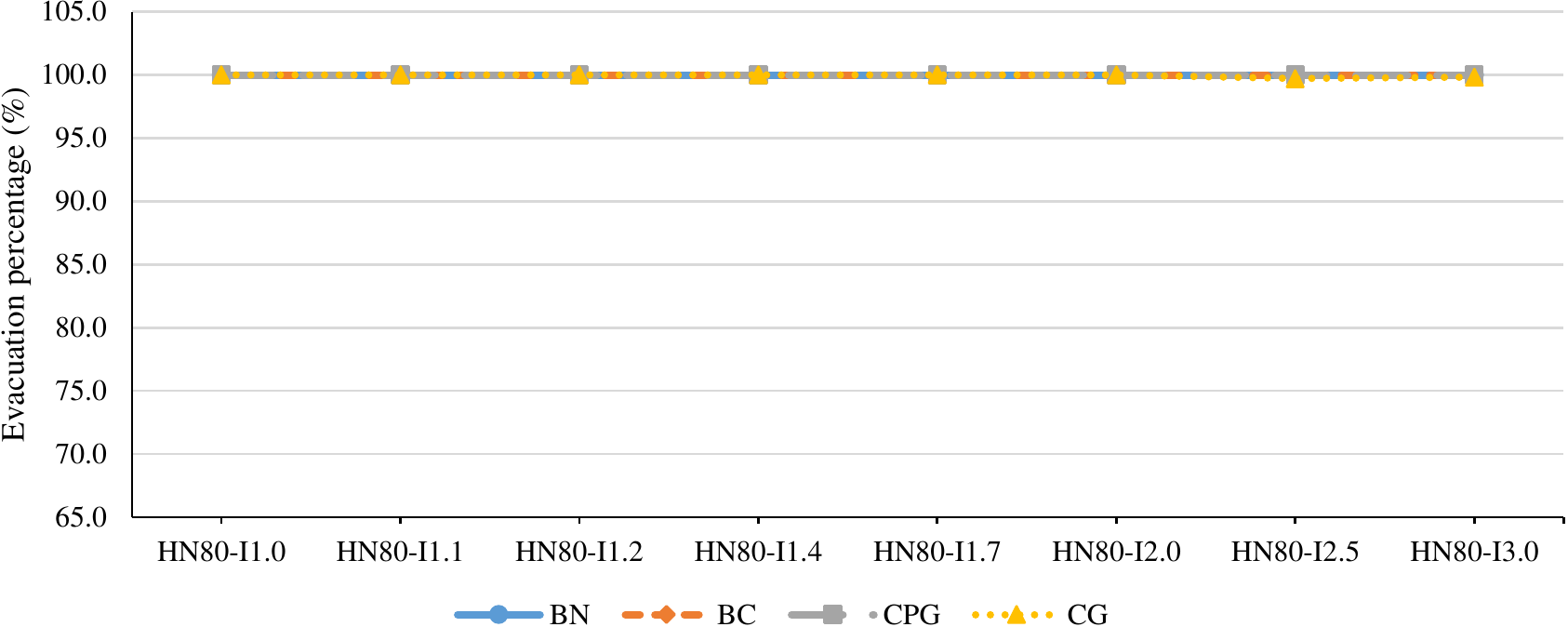}
	\caption{Evacuation percentage under deadline setting with contraflow summary.}
	\label{fig:deadline_summary_withcf}
\end{figure}

Figures \ref{fig:deadline_summary_nocf} and
\ref{fig:deadline_summary_withcf} compare the evacuation percentages
achieved by all four methods. Figure \ref{fig:deadline_summary_withcf}
compares the performance of each method when contraflow is allowed,
and it can be seen that all methods achieve 100\% evacuation for all
instances, except for the CG method which achieves 99.7\% and 99.8\%
evacuation for HN80-I2.5 and HN80-I3.0 respectively. The performance
comparison of each method without contraflow is shown in Figure
\ref{fig:deadline_summary_nocf} which paints a slightly different
picture. All methods achieve 100\% evacuation only up to instance
HN80-I1.4. For instances with larger population scaling factors, the
percentage evacuation starts to drop off significantly for some
methods and slightly for others. A common trend prevalent in these
instances (HN80-Ix with x $\in[1.7,3.0]$) is that the BN method
consistently produces the highest evacuation percentage, followed by
CPG. Although both methods generate evacuation plans with the same
characteristics, their performance disparity could be attributed to
the heuristic nature of CPG: The BN method, although not strictly
optimal, returns optimal results when the integrality constraints on
flows is relaxed.  The lower evacuation rate of the BC and CG methods
can be explained by the additional constraints imposed on their
respective evacuation plans: The BC method produces convergent paths
and the CG method generates non-preemptive evacuation schedules. These
constraints limit their ability of matching the performance of the BN
and CPG methods in a macroscopic analysis. Indeed, the benefits of
these methods cannot be captured in a macroscopic analysis, as they
concern human behavior and the realities of enforcing evacuation
plans.  Note also that the CG method is unique in that the evacuation
rates of its plans are limited to the set of predefined flow rates
$\gamma$ that was specified in Section \ref{sec:casestudy},
whereas the same limitation does not apply to the other methods. The
method's performance is dependent on these preset $\gamma$ values, and
that its performance would change given different sets of $\gamma$
values.

\paragraph{Comparison of the CPU Times}

Tables \ref{tab:deadline_cputime_nocf} and
\ref{tab:deadline_cputime_withcf} compile the total CPU times of all
four methods without and with contraflow respectively. The tables also
show average CPU times for each method across all instances with their
associated standard errors representing uncertainty. The standard
errors are relatively large for most methods due to the spread in CPU
times across the various instances. Nevertheless, a quantitative
comparison can still be made. Regardless of whether contraflow is
allowed or not, the BC and CPG methods consistently consume the
smallest amount of CPU time, followed by the BN method, and finally
the CG method. On top of being the most expensive method, the CG
method reaches its CPU time limit in a few of the more challenging
instances. The difference in run times across the different methods is
not surprising due to the different algorithms employed as well as the
varying constraints imposed on their evacuation plans.

\begin{table}[!tbhp]
	\begin{center}
		\begin{tabular}{|c|r|r|r|r|r|r|r|r|}
			\hline
			\multirow{2}{*}{\textbf{Instance}} & \multicolumn{4}{c|}{\textbf{Total CPU Time (mins)}} \\
			\cline{2-5}
			& \multicolumn{1}{c|}{\textbf{BN}} & \multicolumn{1}{c|}{\textbf{BC}} & \multicolumn{1}{c|}{\textbf{CPG}} & \multicolumn{1}{c|}{\textbf{CG}} \\
			\hline
			HN80-I1.0 & 135 & 0.51 & 0.05 & 1593 \\
			\hline
			HN80-I1.1 & 117 & 10.28 & 0.07 & 1701 \\
			\hline
			HN80-I1.2 & 110 & 10.24 & 0.08 & 2232 \\
			\hline
			HN80-I1.4 & 96 & 1.21 & 0.10 & 2121 \\
			\hline
			HN80-I1.7 & 110 & 10.32 & 0.15 & 2715 \\
			\hline
			HN80-I2.0 & 84 & 1.80 & 0.69 & 7200 \\
			\hline
			HN80-I2.5 & 98 & 3.43 & 11.04 & 7200 \\
			\hline
			HN80-I3.0 & 1449 & 16.22 & 96.48 & 7200 \\
			\hline
			\textbf{Average} & 275 $\pm$ 168 & 6.75 $\pm$ 2.04 & 13.58 $\pm$ 11.92 & 3995 $\pm$ 946 \\
			\hline
		\end{tabular}
	\end{center}
	\caption{CPU times of all methods under deadline setting without contraflow.}
	\label{tab:deadline_cputime_nocf}
\end{table}

\begin{table}[!tbhp]
	\begin{center}
		\begin{tabular}{|c|r|r|r|r|r|r|r|r|}
			\hline
			\multirow{2}{*}{\textbf{Instance}} & \multicolumn{4}{c|}{\textbf{Total CPU Time (mins)}} \\
			\cline{2-5}
			& \multicolumn{1}{c|}{\textbf{BN}} & \multicolumn{1}{c|}{\textbf{BC}} & \multicolumn{1}{c|}{\textbf{CPG}} & \multicolumn{1}{c|}{\textbf{CG}} \\
			\hline
			HN80-I1.0 & 47 & 0.27 & 0.04 & 1465 \\
			\hline
			HN80-I1.1 & 45 & 0.30 & 0.03 & 1495 \\
			\hline
			HN80-I1.2 & 43 & 0.37 & 0.04 & 1577 \\
			\hline
			HN80-I1.4 & 50 & 2.42 & 0.05 & 1744 \\
			\hline
			HN80-I1.7 & 45 & 0.60 & 0.06 & 2269 \\
			\hline
			HN80-I2.0 & 45 & 0.70 & 0.15 & 3339 \\
			\hline
			HN80-I2.5 & 39 & 0.48 & 2.44 & 7181 \\
			\hline
			HN80-I3.0 & 37 & 0.30 & 21.23 & 7200 \\
			\hline
			\textbf{Average} & 44 $\pm$ 2 & 0.68 $\pm$ 0.25 & 3.00 $\pm$ 2.62 & 3284 $\pm$ 880 \\
			\hline
		\end{tabular}
	\end{center}
	\caption{CPU times of all methods under deadline setting with contraflow.}
	\label{tab:deadline_cputime_withcf}
\end{table}

The disparity in CPU times between the BN and BC methods deserves
special attention as their algorithms share a lot of similarities,
both utilizing Benders decomposition. This discrepancy is explained by
the sizes of the RMP and SP of both methods, which are summarized in
Table \ref{tab:bc_bn_constr_var_num}. Although the BC method solves
three problems (the RMP, SP, and DMWP) during each iteration and the
BM method solves only two, the number of variables and constraints in
the latter is significantly larger, leading to correspondingly larger
CPU times. Moreover the nature of the subproblem is fundamentally
different: The Benders subproblems of the BC method are maximum flows,
whereas those of the BN method are multi-commodity flow problems,
where evacuees from an evacuation node can be seen as a single
commodity.

\begin{table}[!tbhp]
	\begin{center}
		\begin{tabular}{|C{43mm}|l|C{42mm}|C{42mm}|}
			\cline{3-4}
			\multicolumn{2}{c|}{} & \multicolumn{1}{c|}{\textbf{BC Method}} & \multicolumn{1}{c|}{\textbf{BN Method}} \\
			\hline
			\multirow{2}{43mm}{\centering \textbf{Restricted Master Problem}} & Constraint \# & $|\mathcal{T}|+|\mathcal{E}\cup\mathcal{T}|+|\mathcal{A}|+|\mathcal{E}|+|\mathcal{C}|+1$ & $|\mathcal{E}|^2+2|\mathcal{T}||\mathcal{E}|+|\mathcal{E}|+|\mathcal{A}||\mathcal{E}|+|\mathcal{A}|+|\mathcal{A}_c|+|\mathcal{C}|+1$ \\
			\cline{2-4}
			& Variable \# & $2|\mathcal{A}|+1$ & $2|\mathcal{A}||\mathcal{E}|+|\mathcal{A}_c|+1$ \\
			\hline
			\multirow{2}{43mm}{\centering \textbf{Subproblem}} & Constraint \# & $|\mathcal{T}||\mathcal{H}|+|\mathcal{A}||\mathcal{H}|+|\mathcal{E}|$ & $|\mathcal{T}||\mathcal{H}||\mathcal{E}|+|\mathcal{E}|+|\mathcal{A}||\mathcal{H}||\mathcal{E}|+|\mathcal{A}||\mathcal{H}|$ \\
			\cline{2-4}
			& Variable \# & $|\mathcal{A}||\mathcal{H}|$ & $|\mathcal{A}||\mathcal{H}||\mathcal{E}|$ \\
			\hline
			\multirow{2}{43mm}{\centering \textbf{Dual of Magnanti-Wong Problem}} & Constraint \# & $|\mathcal{T}||\mathcal{H}|+|\mathcal{A}||\mathcal{H}|+|\mathcal{E}|$ & - \\
			\cline{2-4}
			& Variable \# & $|\mathcal{A}||\mathcal{H}|+1$ & - \\
			\hline
		\end{tabular}
	\end{center}
	\caption{Number of constraints and variables in the problems of the BC and BN methods.}
	\label{tab:bc_bn_constr_var_num}
\end{table}

\subsection{The Minimum Clearance Time Setting}

In the minimum clearance time experiments, each method finds the
smallest amount of time needed to evacuate the entire region, \ie to
achieve 100\% evacuation. A precise definition of minimum clearance
time $h^*$ is given in Equation \eqref{eqn:mincleartime}.

As explained in Section \ref{sec:mincleartime}, the CG method's
multi-objective function simultaneously minimizes the overall
evacuation time and maximizes the evacuation percentage. As such,
$h^*$ can be determined from its deadline setting results for
instances where 100\% evacuation is achieved by identifying the time
at which the last evacuee reaches its safe node. For the BN and BC
methods, the approach outlined in Section \ref{sec:mct_bn_bc} is
applied to find $h^*$. As summarized in Algorithm
\ref{alg:mincleartime}, a binary search procedure is first applied to
find a lower bound $h^\dagger$ to the minimum clearance time. A
sequential search procedure is then applied to find $h^*$. For both
methods, a CPU time limit of 10 minutes is applied when solving each
RMP instance in the binary search procedure. However, in the
subsequent sequential search, a CPU time limit of 10 hours is used for
each search step of the BN method, while a time limit of 2 hours is
used for that of the BC method. The different time limits are selected
to cater for the correspondingly different convergence times of each
method. The approach outlined in Section \ref{sec:mct_cpg} is used to find
$h^*$ for the CPG method. Under the minimum clearance time setting,
the method's RMP and SP are each allocated only 2 minutes of CPU
time.

\begin{table}[!tbhp]
	\begin{center}
		\begin{tabular}{|c|r|r|r|r|r|r|r|r|}
			\hline
			\multirow{2}{*}{\textbf{Instance}} & \multicolumn{2}{c|}{\textbf{BN}} & \multicolumn{2}{c|}{\textbf{BC}} & \multicolumn{2}{c|}{\textbf{CPG}} & \multicolumn{2}{c|}{\textbf{CG}} \\
			\cline{2-9}
			& \multicolumn{1}{C{12mm}|}{\textbf{CPU Time (mins)}} & \multicolumn{1}{C{12mm}|}{\textbf{Min Clear Time (mins)}} & \multicolumn{1}{C{12mm}|}{\textbf{CPU Time (mins)}} & \multicolumn{1}{C{12mm}|}{\textbf{Min Clear Time (mins)}} & \multicolumn{1}{C{12mm}|}{\textbf{CPU Time (mins)}} & \multicolumn{1}{C{12mm}|}{\textbf{Min Clear Time (mins)}} & \multicolumn{1}{C{12mm}|}{\textbf{CPU Time (mins)}} & \multicolumn{1}{C{12mm}|}{\textbf{Min Clear Time (mins)}} \\
			\hline
			HN80-I1.0 & 814   & 260   & 11    & 335   & 29    & 280   & 1593  & 455 \\
			\hline
			HN80-I1.1 & 1936  & 290   & 12    & 365   & 27    & 315   & 1701  & 555 \\
			\hline
			HN80-I1.2 & 650   & 300   & 10    & 395   & 38    & 330   & 2232  & 570 \\
			\hline
			HN80-I1.4 & 1440  & 350   & 13    & 450   & 35    & 380   & 2121  & 600 \\
			\hline
			HN80-I1.7 & 690   & 405   & 10    & 535   & 55    & 445   & -     & - \\
			\hline
			HN80-I2.0 & 705   & 470   & 10    & 630   & 41    & 520   & -     & - \\
			\hline
			HN80-I2.5 & 1374  & 575   & 2     & 770   & 94    & 625   & -     & - \\
			\hline
			HN80-I3.0 & 1128  & 675   & 16    & 925   & 67    & 815   & -     & - \\
			\hline
		\end{tabular}
	\end{center}
	\caption{Minimum clearance times of all methods without contraflow.}
	\label{tab:mincleartime_nocf}
\end{table}

\begin{table}[!tbhp]
	\begin{center}
		\begin{tabular}{|c|r|r|r|r|r|r|r|r|}
			\hline
			\multirow{2}{*}{\textbf{Instance}} & \multicolumn{2}{c|}{\textbf{BN}} & \multicolumn{2}{c|}{\textbf{BC}} & \multicolumn{2}{c|}{\textbf{CPG}} & \multicolumn{2}{c|}{\textbf{CG}} \\
			\cline{2-9}
			& \multicolumn{1}{C{12mm}|}{\textbf{CPU Time (mins)}} & \multicolumn{1}{C{12mm}|}{\textbf{Min Clear Time (mins)}} & \multicolumn{1}{C{12mm}|}{\textbf{CPU Time (mins)}} & \multicolumn{1}{C{12mm}|}{\textbf{Min Clear Time (mins)}} & \multicolumn{1}{C{12mm}|}{\textbf{CPU Time (mins)}} & \multicolumn{1}{C{12mm}|}{\textbf{Min Clear Time (mins)}} & \multicolumn{1}{C{12mm}|}{\textbf{CPU Time (mins)}} & \multicolumn{1}{C{12mm}|}{\textbf{Min Clear Time (mins)}} \\
			\hline
			HN80-I1.0 & 500   & 195   & 26    & 225   & 46    & 200   & 1465  & 235 \\
			\hline
			HN80-I1.1 & 624   & 205   & 27    & 240   & 56    & 210   & 1495  & 255 \\
			\hline
			HN80-I1.2 & 713   & 225   & 47    & 260   & 57    & 235   & 1577  & 285 \\
			\hline
			HN80-I1.4 & 49    & 250   & 3     & 285   & 79    & 260   & 1744  & 320 \\
			\hline
			HN80-I1.7 & 207   & 290   & 1     & 335   & 56    & 320   & 2269  & 575 \\
			\hline
			HN80-I2.0 & 49    & 335   & 10    & 390   & 63    & 385   & 3339  & 535 \\
			\hline
			HN80-I2.5 & 52    & 405   & 10    & 475   & 99    & 440   & -     & - \\
			\hline
			HN80-I3.0 & 59    & 475   & 0     & 560   & 106   & 500   & -     & - \\
			\hline
		\end{tabular}
	\end{center}
	\caption{Minimum clearance times of all methods with contraflow.}
	\label{tab:mincleartime_withcf}
\end{table}

\begin{figure}[!tbhp]
	\centering
	\includegraphics[width=0.98\linewidth]{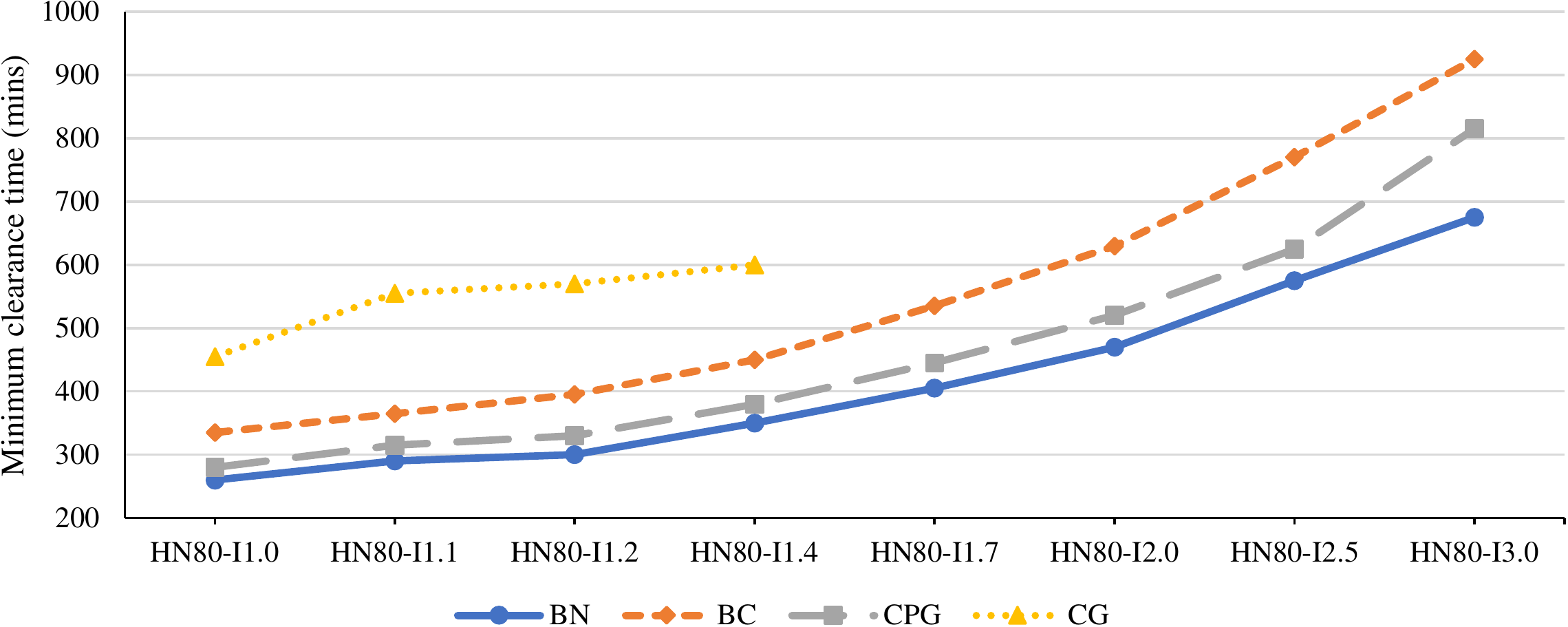}
	\caption{Minimum clearance time without contraflow summary.}
	\label{fig:mincleartime_summary_nocf}
\end{figure}

\begin{figure}[!tbhp]
	\centering
	\includegraphics[width=0.98\linewidth]{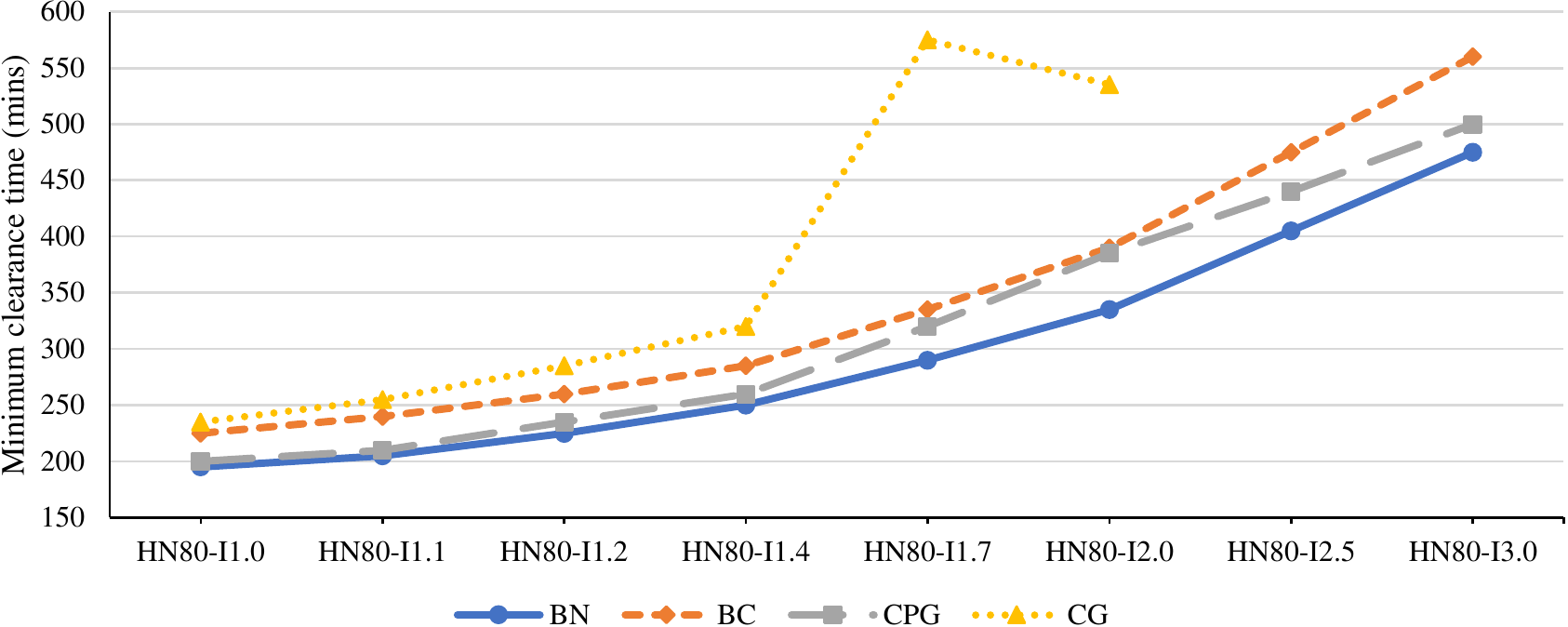}
	\caption{Minimum clearance time with contraflow summary.}
	\label{fig:mincleartime_summary_withcf}
\end{figure}

The total CPU times and the value of $h^*$ for each method are
presented in Tables \ref{tab:mincleartime_nocf} and
\ref{tab:mincleartime_withcf} for cases without and with contraflow
respectively. For the CG method, results are only shown for instances
where the method achieved 100\% evacuation under the deadline
setting. Figures \ref{fig:mincleartime_summary_nocf} and
\ref{fig:mincleartime_summary_withcf} compare each method without and
with contraflow respectively. As expected, the clearance times
increase as the population grows and contraflows help reducing
them. Across all instances, the BN method consistently produces the
smallest clearance time, followed by CPG, BC, and CG. A slight anomaly
is observed in the results of the CG method for instances HN80-I1.7
and HN80-I2.0 when contraflow is allowed, where the clearance time is
smaller for the latter instance.  This result is likely caused by the
early termination of the RMP's last iteration, in which it reached its
24 hours CPU time limit. 

Finally, Table \ref{tab:mincleartime_cputime_summary} shows the
average CPU times consumed by each method under the minimum clearance
time setting. Whether contraflow is allowed or not, the BC method
consistently consumes the least amount of CPU time, followed by CPG,
BN, and CG.

\begin{table}[!tbhp]
	\begin{center}
		\begin{tabular}{|c|c|c|}
			\hline
			\multirow{2}{*}{\textbf{Method}} & \multicolumn{2}{c|}{\textbf{Average CPU Time (mins)}} \\
			\cline{2-3}
			& \multicolumn{1}{c|}{\textbf{Without Contraflow}} & \multicolumn{1}{c|}{\textbf{With Contraflow}} \\
			\hline
			BN  & 1092 $\pm$ 164 & 282 $\pm$ 101  \\
			\hline
			BC  & 11 $\pm$ 1     & 16 $\pm$ 6     \\
			\hline
			CPG & 48 $\pm$ 8     & 70 $\pm$ 8     \\
			\hline
			CG & 1912 $\pm$ 156 & 1981 $\pm$ 297 \\
			\hline
		\end{tabular}
	\end{center}
	\caption{Average CPU times of all methods under minimum clearance time setting.}
	\label{tab:mincleartime_cputime_summary}
\end{table}

\section{Microscopic Evaluation}
\label{sec:microscopic}

The four evacuation methods presented in this paper are macroscopic
and do not capture individual evacuee behaviors, movements, and
interactions, as well as vehicle dynamics such as acceleration and
deceleration, lane changing, and collision avoidance which are all
reflective of what would actually happen in a real world evacuation
scenario. These factors could introduce unanticipated delays or induce
congestion, both of which would negatively affect the performance of
an evacuation plan. Unfortunately, these factors are not captured in
the macroscopic evaluation.

This section presents a microscopic evaluation of the four evacuation
algorithms. Each evacuation plans was simulated, using a road traffic
simulation package SUMO (Simulation of Urban Mobility)
\cite{sumo2012}. SUMO is a full-featured, open-source, microscopic
traffic flow simulation suite developed primarily by researchers of
the Institute of Transportation Systems at the German Aerospace
Center. Its traffic simulator realistically models the traversal
behavior of each vehicle in a road network by computing each vehicle's
instantaneous speed according to the speed limit, a safe distance
to be maintained from a leading vehicle, and the leading vehicle's
speed according to a car-following model described in
\cite{krauss1998} and a lane changing model described in
\cite{erdmann2015}. This simulator not only allows for ascertaining and
evaluating the robustness of the evacuation plans generated in terms of
how they would actually perform in a real world evacuation scenario;
it can also reveal the benefits of convergent and non-preemptive plans. 

The simulations utilize actual road network information of the HN
region, including speed limits, lane counts, and GPS coordinates of
each node to construct the road network for SUMO. They also define the
demand by specifying, for each vehicle, an evacuation path and a
departure time derived from plans generated in the deadline and
minimum clearance time experiments.

\subsection{The Deadline Setting}

\begin{figure}[!tbhp]
	\centering
	\includegraphics[width=0.98\linewidth]{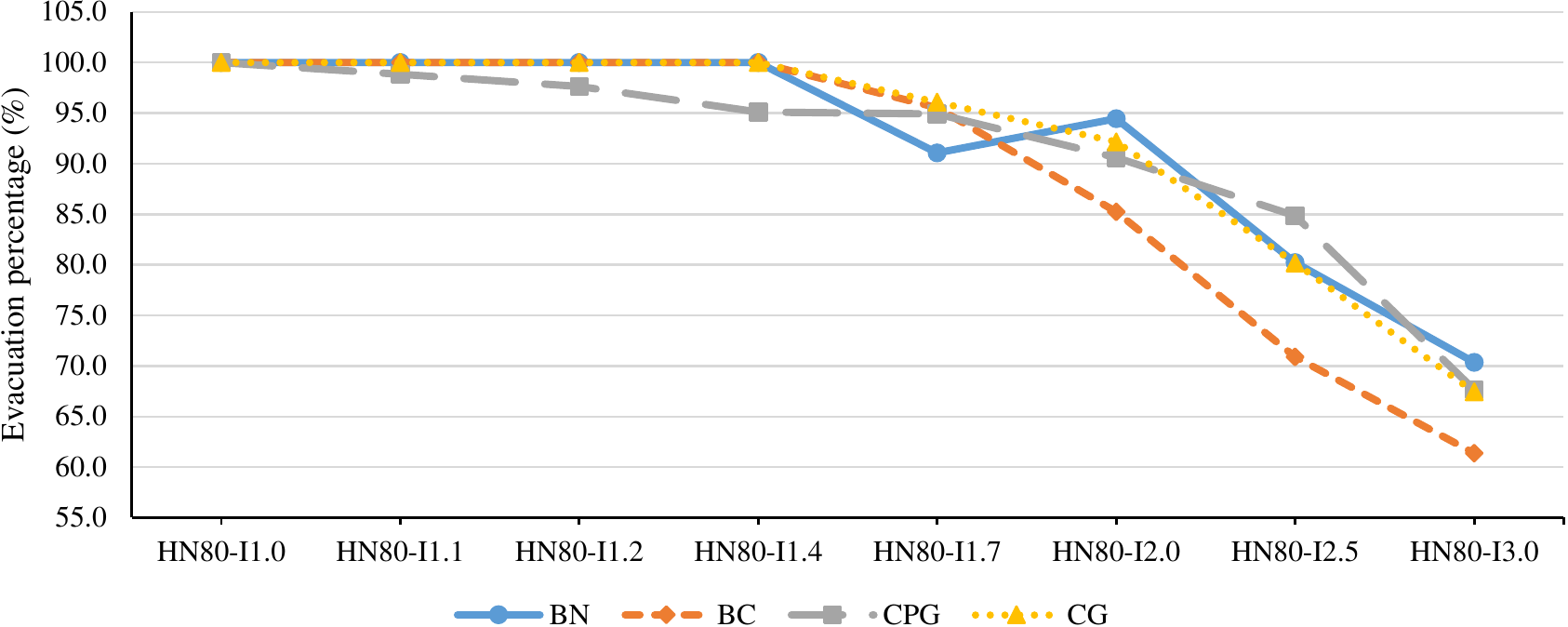}
	\caption{Evacuation percentages obtained from SUMO simulations of deadline setting evacuation plans without contraflow.}
	\label{fig:sumo_deadline_nocf}
\end{figure}

\begin{figure}[!tbhp]
	\centering
	\includegraphics[width=0.98\linewidth]{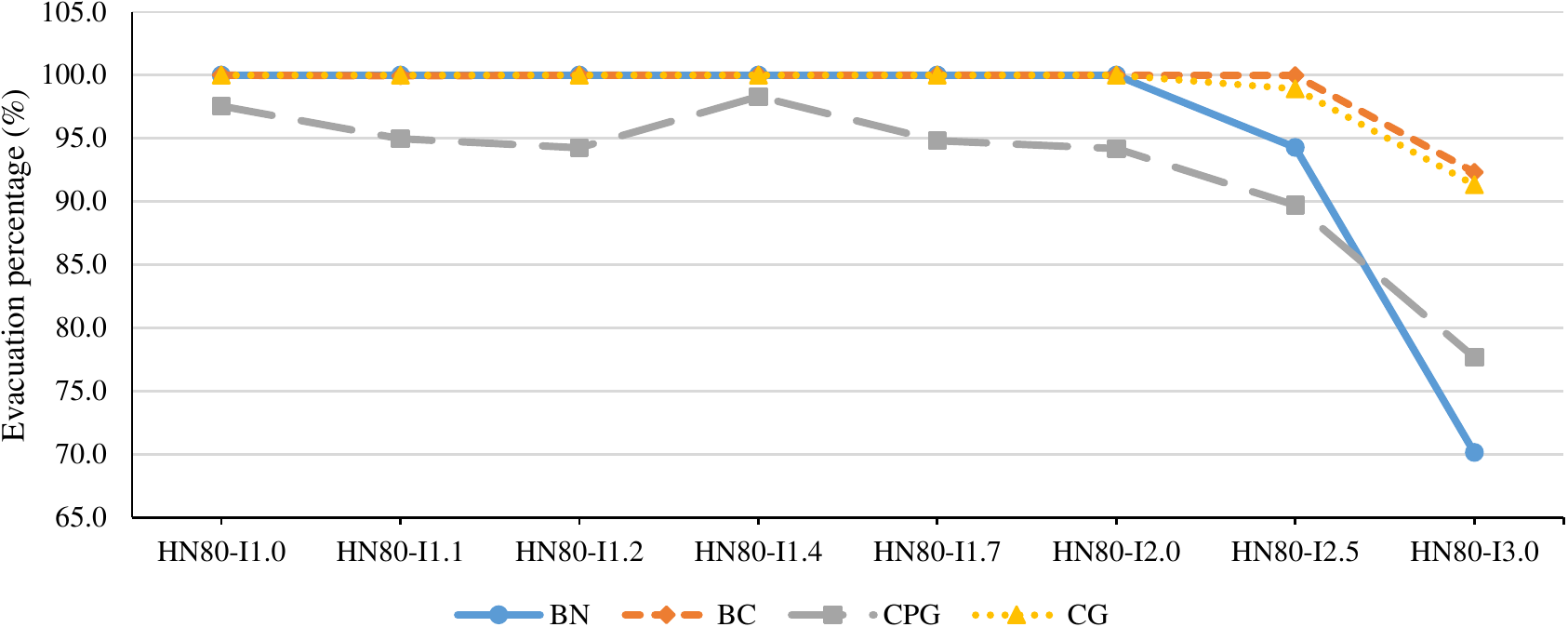}
	\caption{Evacuation percentages obtained from SUMO simulations of deadline setting evacuation plans with contraflow.}
	\label{fig:sumo_deadline_withcf}
\end{figure}

Table \ref{tab:sumo_deadline_results} shows the evacuation percentages
obtained from simulating the evacuation plans generated by the
deadline setting experiments. Figures \ref{fig:sumo_deadline_nocf} and
\ref{fig:sumo_deadline_withcf} present the same information in
graphical form for settings without and with contraflow
respectively. To get a better perspective on the results shown in
these figures, they should be compared with their counterparts from
the deadline experiments, Figures \ref{fig:deadline_summary_nocf} and
\ref{fig:deadline_summary_withcf}, which show corresponding evacuation
percentages produced by the four algorithms. When contraflows are not
allowed, the evacuation percentages start to decrease sooner as the
population scaling factor increases. More importantly, the clear
performance advantage of the BN and CPG methods are not preserved in
these simulations: In some instances, they are outperformed by the
other two methods while, for some other instances, they could barely
outperform the CG method. When contraflows are allowed, the methods
are no longer able to achieve a 100\% evacuation rate. In fact, the
CPG method cannot achieve a 100\% evacuation on any instance, while
the evacuation percentage of other methods start decreasing after
instance HN80-I2.0., with the BN method having the steepest drop.

\begin{table}[!tbhp]
	\begin{center}
		\begin{tabular}{|c|r|r|r|r|r|r|r|r|}
			\hline
			\multirow{3}{*}{\textbf{Instance}} & \multicolumn{8}{c|}{\textbf{Evacuation Percentage (\%)}} \\
			\cline{2-9}
			& \multicolumn{4}{c|}{\textbf{Without Contraflow}} & \multicolumn{4}{c|}{\textbf{With Contraflow}} \\
			\cline{2-9}
			& \multicolumn{1}{c|}{\textbf{BN}} & \multicolumn{1}{c|}{\textbf{BC}} & \multicolumn{1}{c|}{\textbf{CPG}} & \multicolumn{1}{c|}{\textbf{CG}} & \multicolumn{1}{c|}{\textbf{BN}} & \multicolumn{1}{c|}{\textbf{BC}} & \multicolumn{1}{c|}{\textbf{CPG}} & \multicolumn{1}{c|}{\textbf{CG}} \\
			\hline
			HN80-I1.0 & 100.0 & 100.0 & 100.0 & 100.0 & 100.0 & 100.0 & 97.6  & 100.0 \\
			\hline
			HN80-I1.1 & 100.0 & 100.0 & 98.8  & 100.0 & 100.0 & 100.0 & 95.0  & 100.0 \\
			\hline
			HN80-I1.2 & 100.0 & 100.0 & 97.6  & 100.0 & 100.0 & 100.0 & 94.3  & 100.0 \\
			\hline
			HN80-I1.4 & 100.0 & 100.0 & 95.1  & 100.0 & 100.0 & 100.0 & 98.3  & 100.0 \\
			\hline
			HN80-I1.7 & 91.1  & 95.6  & 94.9  & 96.1  & 100.0 & 100.0 & 94.8  & 100.0 \\
			\hline
			HN80-I2.0 & 94.4  & 85.2  & 90.6  & 92.2  & 100.0 & 100.0 & 94.2  & 100.0 \\
			\hline
<			HN80-I2.5 & 80.2  & 70.9  & 84.8  & 80.2  & 94.3  & 100.0 & 89.7  & 98.9 \\
			\hline
			HN80-I3.0 & 70.4  & 61.4  & 67.7  & 67.4  & 70.2  & 92.3  & 77.7  & 91.3 \\
			\hline
		\end{tabular}
	\end{center}
	\caption{Evacuation percentages obtained from SUMO simulations of deadline setting evacuation plans.}
	\label{tab:sumo_deadline_results}
\end{table}

In order to obtain a better understanding of how closely the simulated
results match those of deadline experiments, it is interesting to
normalize the evacuation percentage obtained from the simulations
relative to those from the deadline settings. More precisely, the
normalized evacuation percentage is calculated as follows:
\begin{equation}
\text{Normalized evacuation percentage}=\frac{\text{microscopic evacuation percentage}}{\text{macroscopic evacuation percentage}}
\end{equation}
\noindent
These normalized evacuation percentages are summarized in Table
\ref{tab:sumo_deadline_normalized}. Average values for each method are
also shown with their corresponding uncertainties calculated using
standard errors. A normalized value of 1.0 is achieved by 
most methods for the smaller population instances. However, as
the population increases, the normalized values start to decrease, in some
instances down to 0.70 in the case of the BN method with
contraflow. The general observation is that, regardless of what method
is used and whether contraflows are allowed, microscopic results tend
to diverge from macroscopic results as the population grows. {\em What
  is most interesting however is how the microscopic results change
  the ranking of the methods. The BC and CG methods not only have the
  highest normalized ratios: They also have the highest evacuation
  percentages, especially when contraflows are allowed.}  By using
convergent plans and non-preemptive evacuations, BC and CG produce
plans whose objectives are realistic from a microscopic standpoint.

\begin{table}[!tbhp]
	\begin{center}
		\begin{tabular}{|c|c|c|c|c|c|c|c|c|}
			\hline
			\multirow{3}{*}{\textbf{Instance}} & \multicolumn{8}{c|}{\textbf{Normalized Evacuation Percentage}} \\
			\cline{2-9}
			& \multicolumn{4}{c|}{\textbf{Without Contraflow}} & \multicolumn{4}{c|}{\textbf{With Contraflow}} \\
			\cline{2-9}
			& \multicolumn{1}{c|}{\textbf{BN}} & \multicolumn{1}{c|}{\textbf{BC}} & \multicolumn{1}{c|}{\textbf{CPG}} & \multicolumn{1}{c|}{\textbf{CG}} & \multicolumn{1}{c|}{\textbf{BN}} & \multicolumn{1}{c|}{\textbf{BC}} & \multicolumn{1}{c|}{\textbf{CPG}} & \multicolumn{1}{c|}{\textbf{CG}} \\
			\hline
			HN80-I1.0 & 1.00  & 1.00  & 1.00  & 1.00  & 1.00  & 1.00  & 0.98  & 1.00 \\
			\hline
			HN80-I1.1 & 1.00  & 1.00  & 0.99  & 1.00  & 1.00  & 1.00  & 0.95  & 1.00 \\
			\hline
			HN80-I1.2 & 1.00  & 1.00  & 0.98  & 1.00  & 1.00  & 1.00  & 0.94  & 1.00 \\
			\hline
			HN80-I1.4 & 1.00  & 1.00  & 0.95  & 1.00  & 1.00  & 1.00  & 0.98  & 1.00 \\
			\hline
			HN80-I1.7 & 0.91  & 0.96  & 0.95  & 1.00  & 1.00  & 1.00  & 0.95  & 1.00 \\
			\hline
			HN80-I2.0 & 0.94  & 0.88  & 0.91  & 1.00  & 1.00  & 1.00  & 0.94  & 1.00 \\
			\hline
			HN80-I2.5 & 0.80  & 0.87  & 0.88  & 0.89  & 0.94  & 1.00  & 0.90  & 0.99 \\
			\hline
			HN80-I3.0 & 0.76  & 0.90  & 0.78  & 0.83  & 0.70  & 0.92  & 0.78  & 0.91 \\
			\hline
			\textbf{Average} & \multicolumn{1}{C{11mm}|}{0.93 $\pm$ 0.03} & \multicolumn{1}{C{11mm}|}{0.95 $\pm$ 0.02} & \multicolumn{1}{C{11mm}|}{0.93 $\pm$ 0.03} & \multicolumn{1}{C{11mm}|}{0.97 $\pm$ 0.02} & \multicolumn{1}{C{11mm}|}{0.96 $\pm$ 0.04} & \multicolumn{1}{C{11mm}|}{0.99 $\pm$ 0.01} & \multicolumn{1}{C{11mm}|}{0.93 $\pm$ 0.02} & \multicolumn{1}{C{11mm}|}{0.99 $\pm$ 0.01} \\
			\hline
		\end{tabular}
	\end{center}
	\caption{SUMO evacuation percentages normalized relative to deadline setting evacuation percentages.}
	\label{tab:sumo_deadline_normalized}
\end{table}

\subsection{The Minimum Clearance Time Setting}

\begin{figure}[!tbhp]
	\centering
	\includegraphics[width=0.98\linewidth]{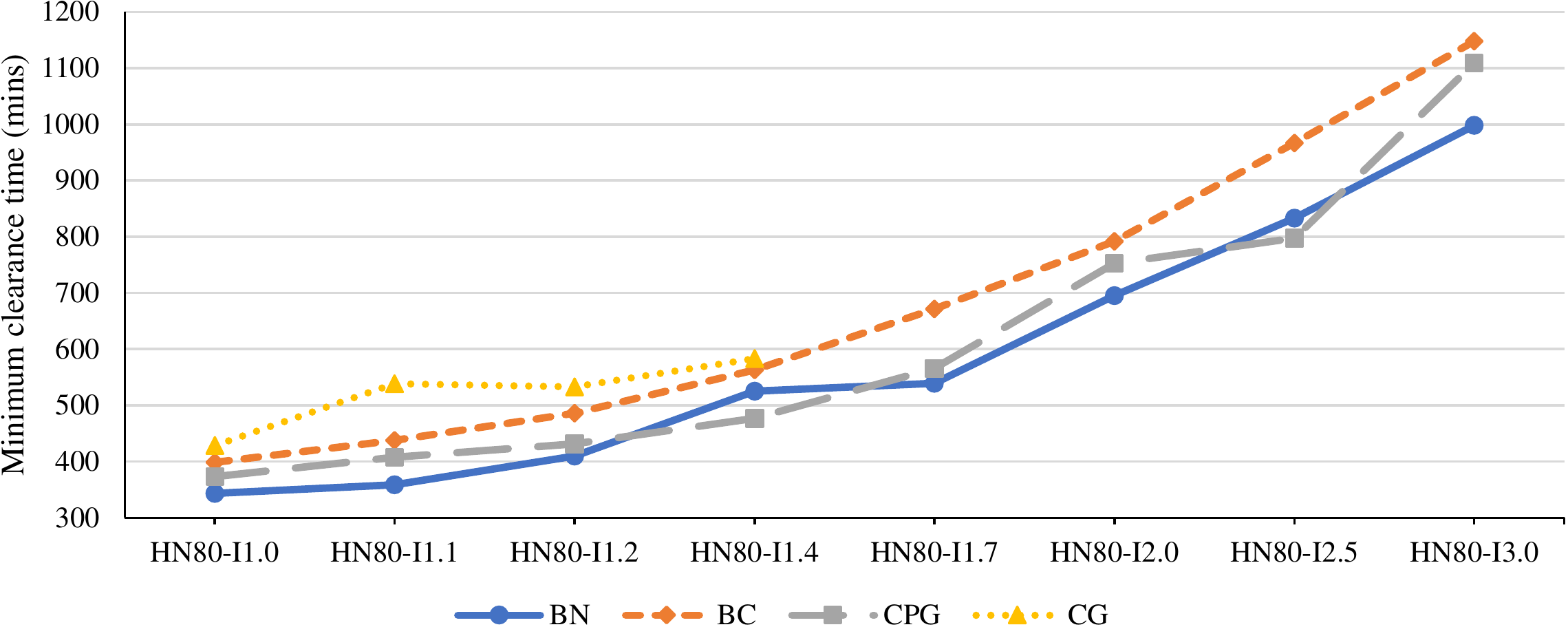}
	\caption{SUMO minimum clearance time without contraflow.}
	\label{fig:sumo_mincleartime_nocf}
\end{figure}

\begin{figure}[!tbhp]
	\centering
	\includegraphics[width=0.98\linewidth]{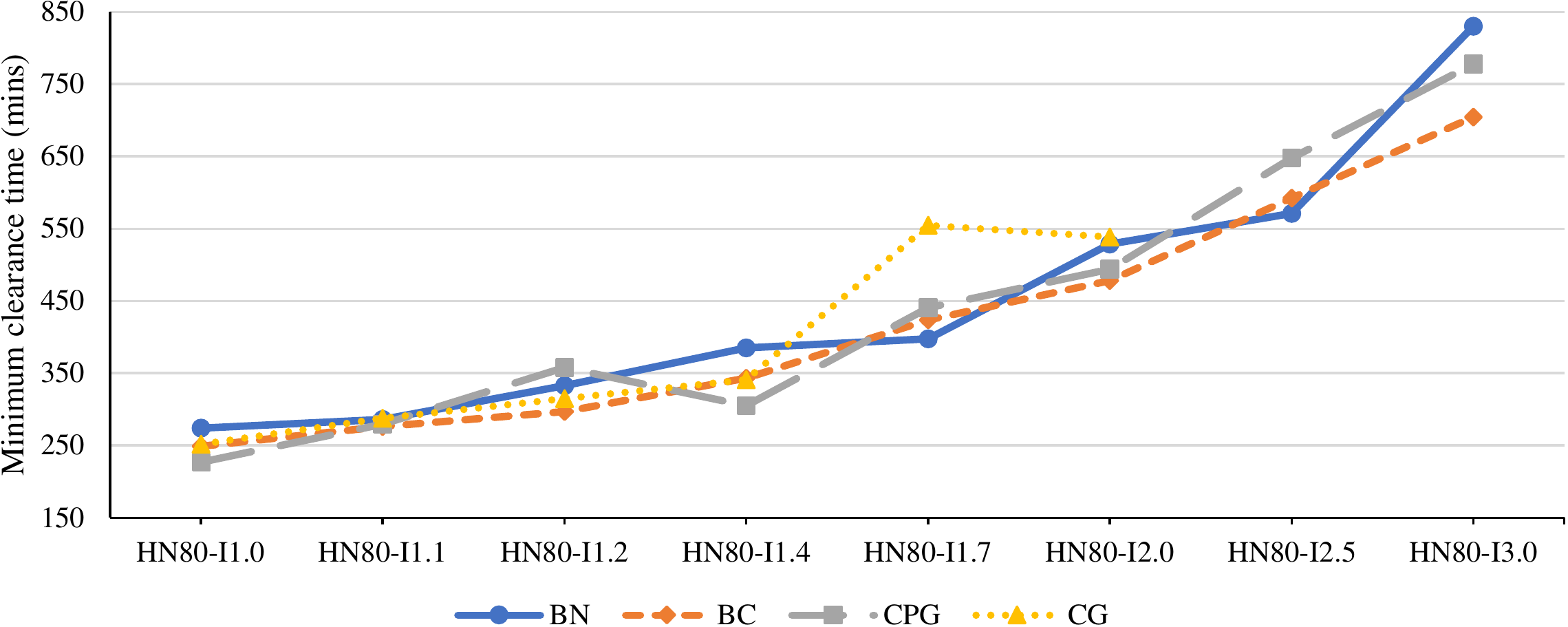}
	\caption{SUMO minimum clearance time with contraflow.}
	\label{fig:sumo_mincleartime_withcf}
\end{figure}

This section presents the microscopic results for the minimum
clearance time setting. Table \ref{tab:sumo_mincleartime_results}
shows minimum clearance times obtained by the simulations, while
Figures \ref{fig:sumo_mincleartime_nocf} and
\ref{fig:sumo_mincleartime_withcf} show the same information in
graphical form without and with contraflows respectively. To place the
results into context, Figure \ref{fig:sumo_mincleartime_nocf} is
compared with Figure \ref{fig:mincleartime_summary_nocf} which contain
corresponding algorithmic results, while Figure
\ref{fig:sumo_mincleartime_withcf} is compared to Figure
\ref{fig:mincleartime_summary_withcf}. The main take-away from these
results is that the gaps between the various methods either decrease
or disappear entirely. For instance, looking at Figure
\ref{fig:sumo_mincleartime_withcf}, there is no clear winner in terms
of minimum clearance time from the simulations.

\begin{table}[!tbhp]
	\begin{center}
		\begin{tabular}{|c|r|r|r|r|r|r|r|r|}
			\hline
			\multirow{3}{*}{\textbf{Instance}} & \multicolumn{8}{c|}{\textbf{Minimum Clearance Time (mins)}} \\
			\cline{2-9}
			& \multicolumn{4}{c|}{\textbf{Without Contraflow}} & \multicolumn{4}{c|}{\textbf{With Contraflow}} \\
			\cline{2-9}
			& \multicolumn{1}{c|}{\textbf{BN}} & \multicolumn{1}{c|}{\textbf{BC}} & \multicolumn{1}{c|}{\textbf{CPG}} & \multicolumn{1}{c|}{\textbf{CG}} & \multicolumn{1}{c|}{\textbf{BN}} & \multicolumn{1}{c|}{\textbf{BC}} & \multicolumn{1}{c|}{\textbf{CPG}} & \multicolumn{1}{c|}{\textbf{CG}} \\
			\hline
			HN80-I1.0 & 343   & 399   & 373   & 428   & 274   & 249   & 227   & 252 \\
			\hline
			HN80-I1.1 & 359   & 438   & 408   & 539   & 286   & 276   & 280   & 288 \\
			\hline
			HN80-I1.2 & 410   & 486   & 432   & 533   & 333   & 297   & 358   & 315 \\
			\hline
			HN80-I1.4 & 526   & 563   & 477   & 583   & 385   & 343   & 305   & 341 \\
			\hline
			HN80-I1.7 & 539   & 672   & 566   & -     & 398   & 424   & 440   & 555 \\
			\hline
			HN80-I2.0 & 695   & 792   & 753   & -     & 529   & 478   & 494   & 539 \\
			\hline
			HN80-I2.5 & 833   & 967   & 797   & -     & 571   & 592   & 648   & - \\
			\hline
			HN80-I3.0 & 998   & 1148  & 1109  & -     & 830   & 705   & 778   & - \\
			\hline
		\end{tabular}
	\end{center}
	\caption{SUMO minimum clearance time results.}
	\label{tab:sumo_mincleartime_results}
\end{table}

\begin{table}[!tbhp]
	\begin{center}
		\begin{tabular}{|c|c|c|c|c|c|c|c|c|}
			\hline
			\multirow{3}{*}{\textbf{Instance}} & \multicolumn{8}{c|}{\textbf{Normalized Minimum Clearance Time}} \\
			\cline{2-9}
			& \multicolumn{4}{c|}{\textbf{Without Contraflow}} & \multicolumn{4}{c|}{\textbf{With Contraflow}} \\
			\cline{2-9}
			& \multicolumn{1}{c|}{\textbf{BN}} & \multicolumn{1}{c|}{\textbf{BC}} & \multicolumn{1}{c|}{\textbf{CPG}} & \multicolumn{1}{c|}{\textbf{CG}} & \multicolumn{1}{c|}{\textbf{BN}} & \multicolumn{1}{c|}{\textbf{BC}} & \multicolumn{1}{c|}{\textbf{CPG}} & \multicolumn{1}{c|}{\textbf{CG}} \\
			\hline
			HN80-I1.0 & 1.32  & 1.19  & 1.33  & 0.94  & 1.41  & 1.11  & 1.13  & 1.07 \\
			\hline
			HN80-I1.1 & 1.24  & 1.20  & 1.29  & 0.97  & 1.40  & 1.15  & 1.33  & 1.13 \\
			\hline
			HN80-I1.2 & 1.37  & 1.23  & 1.31  & 0.94  & 1.48  & 1.14  & 1.52  & 1.10 \\
			\hline
			HN80-I1.4 & 1.50  & 1.25  & 1.26  & 0.97  & 1.54  & 1.20  & 1.17  & 1.06 \\
			\hline
			HN80-I1.7 & 1.33  & 1.26  & 1.27  & -     & 1.37  & 1.26  & 1.38  & 0.96 \\
			\hline
			HN80-I2.0 & 1.48  & 1.26  & 1.45  & -     & 1.58  & 1.23  & 1.28  & 1.01 \\
			\hline
			HN80-I2.5 & 1.45  & 1.26  & 1.28  & -     & 1.41  & 1.25  & 1.47  & - \\
			\hline
			HN80-I3.0 & 1.48  & 1.24  & 1.36  & -     & 1.75  & 1.26  & 1.56  & - \\
			\hline
			\textbf{Average} & \multicolumn{1}{C{11mm}|}{1.40 $\pm$ 0.03} & \multicolumn{1}{C{11mm}|}{1.23 $\pm$ 0.01} & \multicolumn{1}{C{11mm}|}{1.32 $\pm$ 0.02} & \multicolumn{1}{C{11mm}|}{0.95 $\pm$ 0.01} & \multicolumn{1}{C{11mm}|}{1.49 $\pm$ 0.05} & \multicolumn{1}{C{11mm}|}{1.20 $\pm$ 0.02} & \multicolumn{1}{C{11mm}|}{1.36 $\pm$ 0.06} & \multicolumn{1}{C{11mm}|}{1.06 $\pm$ 0.03} \\
			\hline
		\end{tabular}
	\end{center}
	\caption{SUMO minimum clearance time normalized results.}
	\label{tab:sumo_mincleartime_normalized}
\end{table}

Table \ref{tab:sumo_mincleartime_normalized} summarizes the normalized
values for this setting. The table reveals statistically significant
differences in the average normalized values of each method.

The CG method has normalized values which are closest to 1.0, followed
by the BC, CPG, and BN methods. This indicates that the CG method
produces minimum clearance times that are most reproducible in
simulations, \ie their times are the least optimistic and most
accurately reflects what could be achieved in a real world
setting. This is interesting since CG method is the only method that
utilizes non-preemptive evacuation schedules. A possible rationale for
this result is that the step response curves representing evacuation
flow rates over time prevent the transportation network from being
utilized to its full capacity systematically. Indeed, the flow rate is
fixed to a constant value and cannot be increased to saturate
available arc capacities at any given time unlike the flow rates of
preemptive schedules. While this characteristic of the CG hampers its
macroscopic performance, its conservative utilization of the
transportation network has a positive side benefit in that it is less
likely to cause congestion. The fact that the CG method's normalized
values are occasionally less than 1.0 warrants further
explanation. This result is due to the discretization into 5 minute
intervals in the time-expanded graph. This discretization forces
travel time along each arc to be rounded up to the nearest 5 minute
multiple when the arc is represented in the time-expanded graph, even
though the actual travel time maybe less than this multiple. This
quantization causes total travel times along an evacuation path to be
a slight overestimate of actual times. In the simulation of other
methods, this overestimation gets canceled out by congestion induced
delays. However, the CG method produces flows that tend to under
utilize the transportation network, which result in less congestion
and delays. As such, the travel time overestimation does not get
canceled out and is reflected in the normalized clearance times.

\begin{figure}[!tbhp]        
	\centering
	\begin{subfigure}{.5\textwidth}
		\centering
		\includegraphics[width=.9\linewidth]{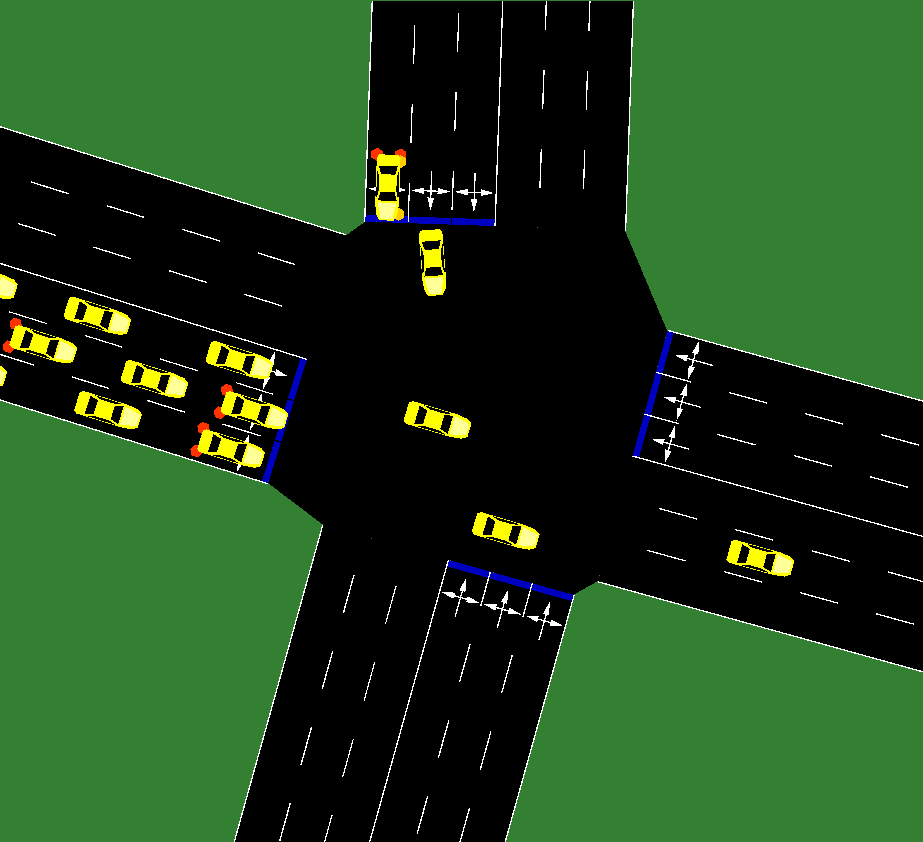}
		\caption{Convergent traffic flow of BC method}
		\label{fig:sumo_bc}
	\end{subfigure}%
	\begin{subfigure}{.5\textwidth}
		\centering
		\includegraphics[width=.9\linewidth]{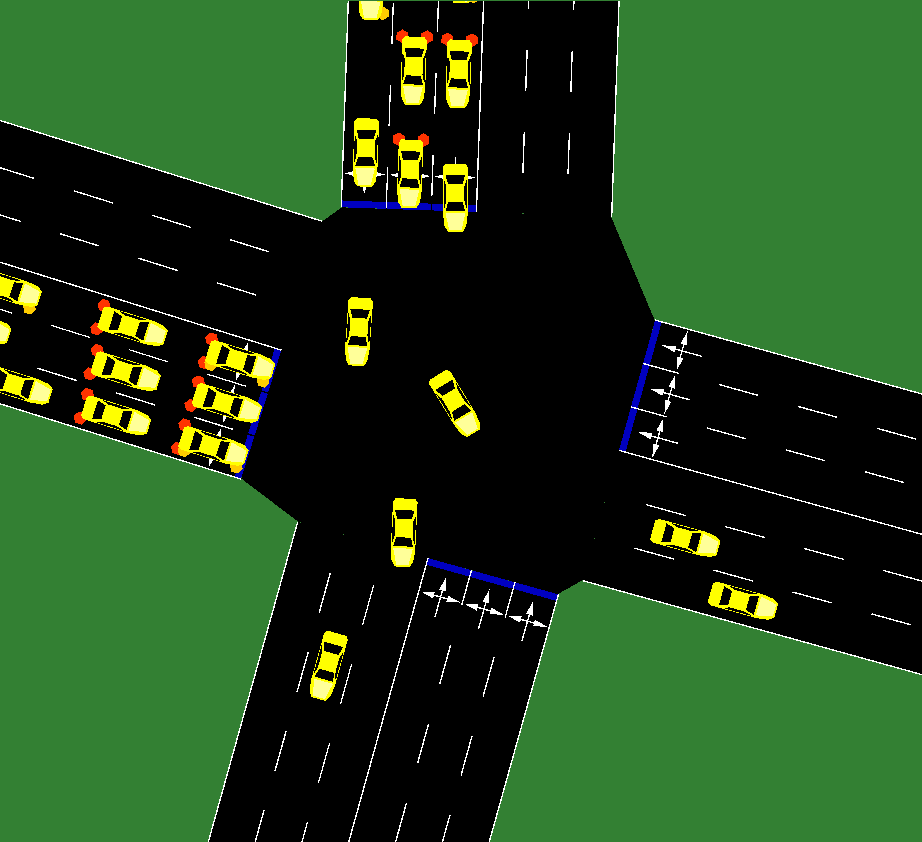}
		\caption{Non-convergent traffic flow of BN method}
		\label{fig:sumo_bn}
	\end{subfigure}
	\caption{Visualization of simulated traffic flow at an intersection produced by evacuaton plans of the BC and BN methods.}
	\label{fig:sumo_bc_bn}
\end{figure}

Method BC also produces normalized values close to 1.0.  To understand
why this was the case, it is useful to take a closer look at the
visualizations of the simulations and, in particular, busy road
intersections and to compare the congestion severity at these
intersections for the the BN and CPG methods. Figures
\ref{fig:sumo_bc} and \ref{fig:sumo_bn} show visualizations of the
same intersection taken from simulations of the BC and BN methods
respectively. The inspections reveal that even though congestion is
present at these intersections in the BC simulations, it was
consistently less severe than those of the other two methods. The BC
method forces traffic flows at intersections to converge onto a single
outgoing road. This limitation is not present in the other methods,
which allow for convergent, divergent, and sometimes even crossing
paths at any given intersection while often lead to more severe
congestion. These observations support the normalized values obtained
for the BC method which indicate its microscopic and macroscopic
results are much closer.

\section{Perspectives}
\label{sec:perspective}

\paragraph{Extensions of the Optimization Algorithms}
The macroscopic and microscopic evaluations contain two important
findings:
\begin{enumerate}
\item Adding constraints on the evacuations, e.g., path convergence
  and non-preemption, increases the fidelity and overall effectiveness
  of the evacuation plans. In particular, for large population sizes,
  problems C-ZEPP and NP-ZEPP produce solutions which are dominated in
  macroscopic evaluations but become superior when validated in
  microscopic simulators.

\item Adding constraints on the evacuations may lead to more elegant
  optimization approaches. Indeed, problems C-ZEPP and NP-ZEPP are
  solved using exact Benders decomposition (method BC) and column
  generation (method CG), which is not the case of Problem ZEPP, which
  is approximated by methods BN and CPG. This elegance may, or may not,
  lead to improved computational performance: Method BC is the most efficient
  method studied in this paper but method CG is the most demanding.
\end{enumerate}

\noindent
An obvious direction for future work is to study convergent,
non-preemptive zone-based evacuation planning. Enhancing method BC
with non-preemption is not an easy task however. The subproblem
becomes a complex scheduling problem that was studied in
\cite{CP2015evacuation}. The key technical issue is to find effective
combinatorial Benders cuts to link the subproblem and the restricted
master problem. Method CG seems easier to extend with convergent
paths, since they only affect the master problem. It remains to
evaluate whether the resulting restricted master problem can be solved
effectively.

\paragraph{Fidelity of the Macroscopic Approaches} This study
demonstrated that macroscopic approaches such as methods BC and CG are
of high fidelity: Their evacuation plans do not degrade in any
significant way when evaluated with a microscopic simulator. This is
an interesting finding on its own. Indeed, it indicates that, for the
case study, it is not necessary to consider more complex models of the
transportation network, such as the Cell Transmission Model (CTM)
\cite{Daganzo1994} which has raised significant interest due to its
elegance and practicality in various settings. However, the CTM is a
mixed integer formulation and hence it adds additional computational
complexity to the models. Some work have considered the linear
relaxation of the CTM: However, this linearization makes it possible
to delay vehicles everywhere in the network, which is not practical in
evacuations and makes the model overly optimistic.

Methods BN and CPG have been shown overly optimistic by the
microscopic simulations. This is an important finding as well: It
suggests that flow-based modeling of evacuations will also be
significantly optimistic. Indeed, flow-based methods decide the route
and timing of each individual evacuee and hence they have much more
flexibility than zone-based evacuations. Similarly, methods not based
on time-expanded graphs and not reasoning about time are unlikely to
be of high fidelity in practice. 

\paragraph{Evacuations Under Uncertainty}
{\em Moving from deterministic evacuation planning algorithms to
  optimization under uncertainty is a critical extension to the
  methods presented in this study}. There are many sources of
uncertainty in evacuations, from the natural or human-made disaster
itself to evacuee compliance and potential accidents. Large-scale
zone-based evacuation planning under uncertainty is largely
unexplored. The work of Andreas and Smith \cite{andreas2009} is a
notable exception: They study a stochastic evacuation planning
problems, where scenarios specify the uncertainty (e.g., the loss of
an arc or an increased travel time). They assign penalties to each arc
in the network and define the cost of a plan as the sum of the
penalties on the arcs used in the evacuation. The resulting problem is
to find a convergent evacuation plan that minimizes the expected
cost. Andreas and Smith apply a Benders decomposition but their
algorithm does not scale beyond a dozen nodes.

It is important to emphasize that, in evacuations, the uncertainty is
not only exogenous (e.g., the flood may be more severe than expected):
It also depends on the evacuation plan itself. Indeed, the probability
of an accident along an arc increases with the flow being scheduled on
the arc. This type of endogenous uncertainty is particularly difficult
to handle efficiently. 

\paragraph{Reliable Evacuation Plans} An interesting direction in evacuation
planning is to borrow some ideas from power systems engineering and
defines various notions of reliability. For instance, the dispatches
of transmission systems satisfy a property known as n-1 reliability
which means that the network can sustain the loss of a line or a
generator and redispatch flows appropriately. Defining and finding n-1
reliable evacuation plans is an interesting avenue for further
research. Within this framework, it would be typical to have a master
problem to find a robust plan and a subproblem for each
contingency. The subproblem will then aim at rerouting the flow of
evacuees affected by the contingency without altering the rest of the
evacuation significantly.

\paragraph{Evacuating Low-Mobility Population and Multi-Modal Evacuations}
The problems studied in this paper assume that each evacuee (or each
household) has a vehicle available for evacuation and do not consider
low-mobility population. Bus evacuations (e.g.,
\cite{Bish2011,doi:10.1002/atr.1224,Goerigk2014,Lee2019,Dhamala2018})
have been studied rather extensively for evacuating low-mobility
population. It is important however to consider the multi-modal
evacuations that integrate both vehicle owners and low-mobility
population, since they share the same transportation network and may
affect the design of contraflows. On-demand multimodal transit systems
\cite{TS2019BusPlus} may bring some significant benefits in that
context, since these novel transit systems solve the first/last mile
problem: They will be able to pick up riders at their location.

\paragraph{Additional Functionalities} There are a number of additional
and critical issues that also deserve to be mentioned. The problem of
reversed evacuations, i.e., bringing residents back to the evacuation area,
is also of great importance and needs to take into account various
cleaning efforts and priorities. The evacuation process itself should
also take into account lodging and other amenities when scheduling an
evacuation for a large area that will not be livable for a significant
period of time. Taking into account these considerations may
significantly reduce the financial burden of many families.

\section{Related Work}
\label{sec:related}

A comprehensive survey of evacuation planning is available in
\cite{bayram2016}. The goal of this section is to not duplicate this
survey but to provide some context for this study. The importance of
prescriptive evacuations, in which evacuation is orchestrated by a
central authority as opposed to self-evacuation, is well recognized in
the field, making it the focal point of various studies.

As mentioned in the introduction, Hamacher and Tjandra
\cite{hamacher2002} distinguishes between microscopic and macroscopic
approaches to evacuation modeling. Microscopic approaches model
individual characteristics of evacuees, their interactions with each
other, and how these factors influence their movement. On the other
hand, macroscopic approaches aggregates evacuees and models their
movements as flow in a network. While all of the algorithms described
This paper has presented macroscopic approaches that were evaluated
using microscopic methods.

{\em The majority of macroscopic approaches are flow-based: They solve
  the evacuation planning problem as a flow on a time-expanded graph.}
For instance, \cite{Lu2003,Lu2005} propose three heuristics to design
an evacuation plan with multiple evacuation routes per evacuation
node, minimizing the time of the last evacuation.  They show that, in
the best case, the proposed heuristic is able to solve randomly
generated instances of up to 50,000 nodes and 150,000 edges in under 6
minutes. Liu \etal \cite{Liu2007a} propose a Heuristic Algorithm for
Staged Traffic Evacuation (HASTE), whose main difference is the use of
the CTM to capture more accurately the flow of evacuees. Lim \etal
\cite{Lim2012} consider a short-notice regional evacuation maximizing
the number of evacuees reaching safety weighted by the severity of the
threat.  They propose two solution approaches to solve the problem,
and present computational experiments on instances derived from the
Houston-Galveston region (USA) with up to 66 nodes, 187 edges, and a
horizon of 192 time steps. Bretschneider and Kimms
\cite{Bretschneider2011, Bretschneider2012} focus on modeling more
accurately the transportation network; they present a free-flow
mathematical model over a detailed street network and include
computational experiments on generated instances with a grid topology
of up to 240 nodes, 330 edges, and considering 150 times steps.

Flow-based approaches are not zone-based and hence do not provide
evacuation plans that comply with the procedures in place in most
emergency services and local authorities. Assigning a routing and
timing to every indivdual evacuee is unlikely to result in actionable
plans. To the best of our knowledge, only a handful of studies design
zone-based evacuation plans that produce both a set of evacuation
routes and an evacuation schedule.  Huibregtse \etal
\cite{Huibregtse2011} propose a two-stage algorithm that first
generates a set of evacuation routes and feasible evacuation times,
and then assigns a route and time to each evacuated area using an ant
colony optimization algorithm.  The main difference with the present
work is that the approach does not explicitly schedule the evacuation
but relies on a third party simulator (EVAQ) to simulate the departure
time of evacuees depending on the evacuation time decided for each
area and evaluate the quality of the solution.  In later work, the
authors studied the robustness of the produced solution
\cite{Huibregtse2010}, and strategies to improve the compliance of
evacuees \cite{Huibregtse2012}. Bish and Sherali \cite{Bish2013}
present a model based on a CTM that assigns a single evacuation path
to each evacuation node but fixes a response curve for each zone.
Computational results include instances with up to 13 evacuation
nodes, 2 safe nodes, and 72 edges. The CPG and CG methods were
originally proposed in \cite{pillac2016} and \cite{pillac2015}
respectively. Even \etal \cite{even2015} takes a different approach by
introducing the Convergent Evacuation Planning Problem (CEPP) which
aims at removing forks from all evacuation routes combined to
eliminate delays caused by diverging and crossing routes. They propose
a two-stage approach to solve the problem which separates route design
and evacuation scheduling. Their work served as the foundation of the
algorithm by Romanski and Van Hentenryck \cite{romanski2016}, one of
the algorithms described in detail in this paper (method BC). As
mentioned earlier, Andreas and Smith \cite{andreas2009} considers a
stochastic evacuation planning with convergent paths. Although not
directly related, it is interesting to mention the work of Chen and
Miller-Hooks \cite{chen2008} who provide an exact technique based on
Benders decomposition to solve the Building Evacuation Problem with
Shared Information. The problem is formulated as a MIP with shared
information constraints to ensure evacuees departing from the same
location at the same departure time receive common instructions and
seeks to route evacuees from multiple locations in a building to exits
such that total evacuation time is minimized.

\section{Conclusion}
\label{sec:conclusion}

This paper presented a systematic study large-scale zone-based
evacuation planning, both from an effectiveness and a computational
standpoint. In zone-based evacuation planning, the region to evacuate
is divided into zones and each zone must be assigned a path to safety
and departure times along the path. Zone-based evacuations are highly
desirable in practice because they allow emergency services to
communicate evacuation orders and to control the evacuation more
accurately.

The paper reviewed several existing optimization algorithms, and
presented new ones, and evaluated them, on a real, large-scale case
study, both from a macroscopic standpoint and through microscopic
simulations under a variety of assumptions. In particular, the paper
evaluated the impact of contraflow, preemption, and convergent
evacuations.

The main take-away is probably the success of mathematical programming
for designing macroscopic evacuation planning algorithms that are
effective when evaluated by microscopic simulators on large-scale
scenarios. The algorithms thus provide both primal solutions and a
performance guarantee. A second key take-away is the importance of
constraining the optimization algorithms to produce realistic plans.
By imposing the non-preemption constraints, the optimization algorithm
produces evacuation plans that are easy to enforce and of high
fidelity.  By imposing convergence properties, the optimization
algorithm avoid congestion and driver hesitations. The optimization
algorithms have been shown to be highly practical and some of them can
even be used in real-time settings.

The papers also listed a number of directions for future work and
perspectives on the field. Perhaps the most pressing issues center
around delivering plans that are robust under uncertainty and to
support multimodal evacuations that evacuate both vehicle owners and
low-mobility populations. 

\section*{Acknowledgements}

We would like to express our deep gratitude to Peter Cinque (New
South-Wales State Emergency Services) and Peter Liehn (then at NICTA)
for their leadership in the NICTA evacuation project, the foundation
of this work. Manuel Cebrian, Caroline Even, Victor Pillac, and
Andreas Schutt drove the development of many of the algorithms
presented here and have been truly amazing collaborators.

\end{document}